\documentclass[a4paper,12pt]{article}
\usepackage{amsthm}
\usepackage{amsmath}
\usepackage{amssymb}
\usepackage{tabularx}
\usepackage{indentfirst}
\usepackage{enumerate}
\usepackage{dsfont}
\usepackage{graphicx}
\usepackage{subfigure}
\usepackage{epsfig}
\usepackage{graphics}
\usepackage{cases}
\usepackage[compress]{cite}
\usepackage{txfonts}
\usepackage{geometry}
 \usepackage{epstopdf}
 \usepackage{lineno}
 \usepackage[justification=centering]{caption}
 \usepackage{color}

 \usepackage{hyperref}
\hypersetup{
colorlinks=true,
linkcolor=blue,
anchorcolor=blue,
citecolor=blue}

\makeatletter
\newsavebox\myboxA
\newsavebox\myboxB
\newlength\mylenA
\newcommand*\xoverline[2][0.8]{%
\sbox{\myboxA}{$\m@th#2$}%
\setbox\myboxB\null
\ht\myboxB=\ht\myboxA%
\dp\myboxB=\dp\myboxA%
\wd\myboxB=#1\wd\myboxA
\sbox\myboxB{$\m@th\overline{\copy\myboxB}$}
\setlength\mylenA{\the\wd\myboxA}
\addtolength\mylenA{-\the\wd\myboxB}%
\ifdim\wd\myboxB<\wd\myboxA%
\rlap{\hskip 0.5\mylenA\usebox\myboxB}{\usebox\myboxA}%
\else
\hskip -0.5\mylenA\rlap{\usebox\myboxA}{\hskip 0.5\mylenA\usebox\myboxB}%
\fi}
\makeatother
\makeatletter
\newcommand\dlmu[2][0.2cm]{\hskip0.5pt\underline{\hb@xt@ #1{\hss#2\hss}}\hskip0.5pt}
\makeatother
\newtheorem{theorem}{Theorem}[section]
\newtheorem{lemma}{Lemma}[section]

\newtheorem{exm}{Example}[section]

\newtheorem{Remark}{Remark}[section]
\newtheorem{definition}{Definition}[section]

\numberwithin{equation}{section}

\newcommand{\be}{\begin{equation}}
\newcommand{\ee}{\end{equation}}
\newcommand\bes{\begin{eqnarray}}
\newcommand\ees{\end{eqnarray}}
\newcommand{\bess}{\begin{eqnarray*}}
\newcommand{\eess}{\end{eqnarray*}}

\begin{document}

\thispagestyle{empty}

\title{On an impulsive faecal-oral model in a periodically evolving environment \thanks{The second author is supported by the National Natural Science Foundation of China (No. 12271470) and the third author acknowledges the support of CNPq/Brazil Proc. $N^{o}$ $311562/2020-5$ and FAPDF grant 00193.00001133/2021-80.}}

\date{\empty}

\author{Qi Zhou$^1$, Zhigui Lin$^1 \thanks{Corresponding author. Email: zglin@yzu.edu.cn (Z. Lin).}$ and Carlos Alberto Santos$^{2}$
\\
{\small $^1$ School of Mathematical Science, Yangzhou University, Yangzhou 225002, China}\\
{\small $^2$ Department of Mathematics, University of Brasilia, BR-70910900 Brasilia, DF, Brazil}
}
 \maketitle
\begin{quote}
\noindent
{\bf Abstract.}{\footnotesize\small~To understand how impulsive intervention and regional evolution jointly influence the spread of faecal-oral diseases, this paper develops an impulsive faecal-oral model in a periodically evolving environment. The well-posedness of the model is first checked. Then, the existence of the principal eigenvalue dependent on impulse intensity and evolving rate is proved based on Krein-Rutman theorem. With the help of this value, the threshold dynamical behaviours of the model are explored. More importantly, this paper also derives the monotonicity of the principal eigenvalue with respect to initial region and impulse intensity and estimates the principal eigenvalue in some special cases. Finally, numerical simulations are used to verify the correctness of the theoretical results and to explore the impact of regional evolution rate on the spread of the diseases. Our research shows that large impulsive intensity $1-g'(0)$ and small evolving rate $\rho(t)$ play a positive role in the prevention and control of the diseases.
}

\noindent {\bf MSC:} 35K57, 
35R12, 
92B05 

\medskip
\noindent {\bf Keywords:} Faecal-oral model; Impulsive intervention; Evolving environment; Dynamical behaviours; Principal eigenvalue
\end{quote}

\section{Introduction}
Faecal-oral pathway is a particular route of transmission of a disease wherein pathogenic bacteria are excreted in the faeces of patients or carriers, and then ingested by infected persons through various means \cite{1}. This kind of pathway can be used to describe the spread of infectious diseases such as hand-foot-mouth diseases, cholera and so on \cite{2}, and of epidemic diseases such as typhoid fever, infectious hepatitis, poliomyelitis and the like \cite{3}. In addition, it has been suggested from \cite{4, 5} that COVID-19 may also be transmitted by this route. Therefore, faecal-oral pathway has been widely concerned by some epidemiologists, mathematicians, and medical scientists.

Mathematical models have become potentially powerful tools for understanding the intrinsic transmission routes, anticipating the trends and designing control measures of disease \cite{6}. In the past several decades, some faecal-oral models have been proposed  to study the evolution of faecally-orally transmitted diseases, see, for instance, \cite{7, 8, 9, 10}. To the best of authors' knowledge, to investigate the cholera epidemic that spread in the European Mediterranean regions in 1973, Capasso and Paveri-Fontana \cite{7} first proposed the following faecal-oral model
\begin{eqnarray}\label{1-1}
\left\{
\begin{array}{ll}
\frac{d u}{d t}=-a_{11}u+a_{12}v,\\
\frac{d v}{d t}=-a_{22}v+f(u)
\end{array}
\right.
\end{eqnarray}
with suitable initial values $u_{0}$ and $v_{0}$, where $u(t)$ and $v(t)$ denote the average densities of pathogenic bacteria and infected individuals at time $t$, respectively. $f(u)$ is the growth function of infected individuals caused by the bacteria. Constants $a_{11}$, $a_{12}$ and $a_{22}$ are all positive. For a more detailed biological explanation of model \eqref{1-1}, the interested reader can refer to \cite{7}. Taking into account the influence of spatial diffusion, Capasso and Maddalena modified model \eqref{1-1} into the reaction-diffusion models with Neumann boundary conditions \cite{8, 9} and Robin boundary conditions \cite{10}, and explored the dynamical behaviours of the models.

Notice that all already mentioned works on faecal-oral models were done on a fixed domain, that is, the domain $\Omega$ is independent of time. In fact, the habitat of biological organisms may exhibit periodic change with time, that is, $\Omega_{t}\in\mathds{R}^{n}(n\geq 1)$ is a time-varying region with period $\tau$. For example, seasonal shifts cause changes in the depths of lakes \cite{11}, and climate variations may contribute to changes in the region of infection for mosquito borne diseases \cite{12}. In recent years, some mathematical models on periodically evolving region have been established in order to understand as clearly as possible the influence of altered region on the spread of diseases, see, for instance, \cite{13, 14, 15} and the references therein. Among these works, Zhou, Zhang and Ling \cite{13} developed the following faecal-oral model
\begin{eqnarray}\label{1-2}
\left\{
\begin{array}{ll}
\frac{\partial u}{\partial t}=\frac{d_{1}}{\rho^{2}(t)}\Delta u-\frac{n\dot{\rho}(t)}{\rho(t)}u-a_{11}u+a_{12}v,\; &\, t> 0, x\in\Omega_{0}, \\[2mm]
\frac{\partial v}{\partial t}=\frac{d_{2}}{\rho^{2}(t)}\Delta v-\frac{n\dot{\rho}(t)}{\rho(t)}v-a_{22}v+f(u),\; &\, t> 0, x\in\Omega_{0}, \\[2mm]
u(0,x)=u_{0}(x),v(0,x)=v_{0}(x),\; &\, x\in \bar{\Omega}_{0}
\end{array}
\right.
\end{eqnarray}
with Neumann boundary conditions, where habitat $\Omega_{0}\in\mathds{R}^{n}(n\geq1)$ be a bounded smooth domain with smooth boundary $\partial\Omega_0$. $u(t,x)$ and $v(t,x)$ respectively represent the spatial densities of pathogenic bacteria and infected individuals at time $t$ and location $x$, and $\rho(t)$ is a positive periodic function and satisfies $\rho(0)=1$. $d_{1}$ and $d_{2}$ are positive constants that denote the diffusion coefficients of pathogenic bacteria and infected individuals, respectively. Here, the meanings of other symbols are the same as in model \eqref{1-1}. With regard to the more detailed establishment procedure of model \eqref{1-2}, interested readers can refer to \cite{13}.

For diseases that are transmitted by faecal-oral route, environmental disinfection has become one of the most effective measures for control and prevention.
Note that this intervention usually leads to a sharp reduction in the number of pathogenic bacteria in a short period of time, and has an important influence on the spread of the diseases. Classical differential equation model cannot capture this instantaneous phenomenon. Scholars have thus turned their interest to impulse differential equation models that can capture this phenomenon.

Theoretical research on impulse differential equation started by Mil'man and Myshkis \cite{17}, and has been greatly developed after 1980s. In the field of biomathematics, a pioneering work on impulsive reaction-diffusion model has been given by Lewis and Li \cite{18}. Based on work \cite{18}, there has been a series of impulsive researches. We refer to \cite{19} for a diffusive logistic model, \cite{20} for a juvenile-adult model, \cite{21} for a competition model in stream environments, \cite{22} for a model with a nonlocal discrete time map, \cite{23} for a reaction-diffusion model with non-monotone birth function and age structure, \cite{23-1} for a model with shifting environments, \cite{23-2} for a predator-prey model with seasonality and fear effect, \cite{23-3} for a pulsed-precipitation model of dryland vegetation pattern formation, \cite{23-4} for a stage-structured continuous-discrete-time population model, and references therein. However, to the best of the authors' knowledge, no research has been carried out on impulsive faecal-oral model in a periodically evolving environment or in a fixed environment.

In this paper, we develop an impulsive faecal-oral model in a periodically evolving environment by incorporating the pulsed intervention into model \eqref{1-2}.
It is natural to ask how the combination of pulsed intervention and regional evolution can influence or even change the dynamical behaviours of the diseases? This is the research motivation of this paper. The novelty and contribution of this work are summarized as follows:
\begin{itemize}
\item{An impulsive faecal-oral model in a periodically evolving environment is developed via taking into account the pulsed intervention.}
\item{The existence of the principal eigenvalue is proved. Based on this, the global threshold dynamics of the model are explored.}
\item{The effects of habitat and impulse intensity on the dynamical behaviours are given.}
\end{itemize}

The rest of the paper is organized as follows. Model derivation and well-posedness are presented in \autoref{Section-2}. \textcolor{blue}{Section}~\ref{Section-3} proves the existence of the principal eigenvalue of the periodic linear problem corresponding to the impulse model. The global threshold dynamics of the model are explored in \autoref{Section-4}. \textcolor{blue}{Section}~\ref{Section-5} shows the monotonicity of the principal eigenvalue with respect to initial region and impulse intensity, and estimates the principal eigenvalue in some special cases. Through numerical simulations, \autoref{Section-6} verifies the correctness of the theoretical results obtained in the previous sections, and explores the influence of evolutionary rate on the spread of the diseases.  Finally, \autoref{Section-7} concludes the paper.
\section{Model derivation and well-posedness}\label{Section-2}
In this section, we first introduce the model and then prove that this model possesses a unique nonnegative solution.

As stated in introduction, impulsive phenomenon are often present during the spread of faecal-oral diseases. By incorporating this phenomenon into model \eqref{1-2}, the impulsive faecal-oral model in a periodically evolving environment is given by the following equations:
\begin{eqnarray}\label{1-2-20}
\left\{
\begin{array}{ll}
\frac{\partial u}{\partial t}=\frac{d_{1}}{\rho^{2}(t)}\Delta u-\frac{n\dot{\rho}(t)}{\rho(t)}u-a_{11}u+a_{12}v,\; &\, t\in((k\tau)^{+}, (k+1)\tau], x\in\Omega_{0}, \\[2mm]
\frac{\partial v}{\partial t}=\frac{d_{2}}{\rho^{2}(t)}\Delta v-\frac{n\dot{\rho}(t)}{\rho(t)}v-a_{22}v+f(u),\; &\, t\in((k\tau)^{+}, (k+1)\tau], x\in\Omega_{0}, \\[2mm]
u=0, v=0,\; &\, t\in\mathds{R}^{+}, x\in\partial\Omega_{0},\\[2mm]
u(0,x)=u_{0}(x), v(0,x)=v_{0}(x),\; &\,x\in\bar{\Omega}_{0},\\[2mm]
u((k\tau)^{+},x)=g(u(k\tau,x)), v((k\tau)^{+},x)=v((k\tau),x), \; &\, x\in\Omega_{0}, k=0,1,2,\cdots,
\end{array} \right.
\end{eqnarray}
where $\tau$ denotes the time gap of two neighboring pulse interventions, and $(k\tau)^{+}$ stands for the right limit at $k\tau$. For the convenience of writing, we always omit $k=0,1,2,\cdots$, unless otherwise specified. Here, the meanings of other symbols are the same as in model \eqref{1-2}. In the last line, $u((k\tau)^{+},x)=g(u(k\tau,x))$ represents the spatial density of the pathogenic bacteria at position $x$ after impulsive intervention is applied at time $k\tau$, and $v((k\tau)^{+},x)=v((k\tau),x)$ denotes that there is no relevant control measure applied to infected individual. The first two lines denote that $u$ and $v$ satisfy the equations for $(t,x)\in((k\tau)^{+}, (k+1)\tau]\times\Omega_{0}$ and take the initial value $u((k\tau)^{+},x)$ and $v((k\tau)^{+},x)$. The third line means that the model has Dirichlet boundary condition. This boundary condition can describe the fact that infected individuals and pathogenic bacteria vanish at the boundary of habitat \cite{24}.

In this paper, to guarantee the existence and uniqueness of the classical solution of model \eqref{1-2-20}, assume that $u_{0}(x), v_{0}(x)\in\mathbb{C}^{1}(\bar{\Omega}_{0})$ and $u_{0}(x)$, $v_{0}(x)\geq 0(\not \equiv 0)$, and that $\rho(t)\in\mathbb{C}^{1}([0,\infty])$. For the convenience of writing, let $h_{1}=f$ and $h_{2}=g$. Moreover, we make the following assumptions about growth function $h_{1}$ and impulsive function $h_{2}$:
\begin{itemize}
\item[(A1)] For each $i\in\{1,2\}$, $h_{i}(u)\in\mathbb{C}^{1}([0,\infty))$, $h'_{i}(0)>0$ and $h'_{i}(u)\geq 0$ for $u>0$, and $h_{i}(u)\geq 0$ for $u\geq0$ and $h_{i}(u)=0$ if and only if $u=0$.
\item[(A2)] $\frac{h_{1}(u)}{u}$ is decreasing for $u>0$ and $\lim\limits_{u\rightarrow +\infty}\frac{h_{1}(u)}{u}<\frac{(a_{11}+\check{m})(a_{22}+\check{m})}{a_{12}}$, where $\check{m}=\min\limits_{t\in [0, \tau]}\frac{n\dot{\rho}(t)}{\rho(t)}$ and $a_{11}, a_{22}>|\check{m}|$.
\item[(A3)] $\frac{h_{2}(u)}{u}$ is decreasing for $u>0$ and $0<\frac{h_{2}(u)}{u}\leq1$ for $u>0$.
\item[(A4)] For each $i\in\{1,2\}$, there exist positive constants $H_{i}$, $\varpi_{i}$ and $\kappa_{i}>1$ such that $h_{i}(u)\geq h'_{i}(0)u-H_{i}u^{\kappa_{i}}$ for $0\leq u\leq\varpi_{i}$.
\end{itemize}
In the sequel, the assumptions about initial values $u_{0}(x), v_{0}(x)$, growth function $f$ and pulse function $g$ always hold without further explanation.
\begin{Remark}
For Assumptions (A1)-(A4), some satisfied functions are as follows:
\begin{itemize}
\item[(1)] linear functions $f(u)=c_{1}u$, $g(x)=c_{2}u$, where $c_{1}\in (0, \frac{(a_{11}+\check{m})(a_{22}+\check{m})}{a_{12}})$ and $c_{2}\in(0,1]$ \cite{20};
\item[(2)] Beverton-Holt functions $f(u)=\frac{m_{1}u}{a_{1}+u}$, $g(x)=\frac{m_{2}u}{a_{2}+u}$, where $a_{1}, a_{2}, m_{1}\in(0, \infty)$ and $m_{2}\in(0,a_{2}]$ \cite{24};
\item[(3)] exponential functions $f(u)=b_{1}e^{-r_{1}t}u$, $g(u)=b_{2}e^{-r_{2}t}u$, where $r_{1}, r_{2}\in(0, \infty)$, $b_{1}\in (0, \frac{(a_{11}+\check{m})(a_{22}+\check{m})}{a_{12}})$ and $b_{2}\in(0,1]$ \cite{25}.
\end{itemize}
\end{Remark}
Now, an impulsive faecal-oral model is introduced. In order to make the model more realistic, the solution of the model must be existent, unique, and nonnegative. We will study this issue below.

For the convenience of stating the solution of impulse differential equation, we introduce the following notation:
\begin{equation}
\begin{aligned}
\mathds{PC}:=\mathds{PC}\big([0,\infty)\times\bar{\Omega}_{0}\big)&:=\Big\{u:u\in \cap_{k\in\mathbb{N}}\mathbb{C}\big((k\tau, (k+1)\tau]\times \bar{\Omega}_{0} \big)\Big\};\\
\mathds{PC}^{2}:=\mathds{PC}^{2}\big([0,&\infty)\times\bar{\Omega}_{0}\big):=\Big\{(u, v):u, v\in \mathds{PC}\Big\};\\
\mathds{PC}_{1,2}:=\mathds{PC}_{1,2}\big((0,\infty)\times\Omega_{0}\big)&:=\Big\{u:u\in \cap_{k\in\mathbb{N}}\mathbb{C}^{1,2}\big((k\tau, (k+1)\tau]\times \Omega_{0} \big)\Big\};\\
\mathds{PC}_{1,2}^{2}:=\mathds{PC}_{1,2}^{2}\big(&[0,\infty)\times\bar{\Omega}_{0}\big):=\Big\{(u, v):u, v\in \mathds{PC}_{1,2}\Big\}.
\nonumber
\end{aligned}
\end{equation}
Inspired by work \cite{25}, we next attempt to give the well-posedness of the solution of model \eqref{1-2-20}.
\begin{theorem}\label{theorem 1-1}
For any initial value $(u_{0}(x), v_{0}(x))$, there exists a unique nonnegative solution $\big(u(t,x), v(t,\\x)\big)\in\mathds{PC}_{1,2}^{2}\cap\mathds{PC}^{2}$ for the model \eqref{1-2-20}. Furthermore, the solution of model \eqref{1-2-20} satisfies that $u(t,x), v(t,x)>0$ for $t>0$ and $x\in\Omega_{0}$.
\begin{proof}
The conclusion that model \eqref{1-2-20} admits a unique locally nonnegative solution $(u(t,x), v(t,x))$ on the maximal existence interval $[0, t_{\max})$ can be obtained by using a similar approach to \cite[Theorem 3.1]{25}. The proof of the boundedness of $u(t,x)$ and $v(t,x)$ is given below. By Assumption (A2), it follows that for a fixed sufficiently small $\epsilon_{0}>0$, there exists a sufficiently large constant
\begin{equation*}
G>\max\bigg\{\max\limits_{\bar{\Omega}_{0}}u_{0}(x), \frac{a_{12}+\epsilon_{0}}{a_{11}+\check{m}}\max\limits_{\bar{\Omega}_{0}}v_{0}(x)\bigg\}
\end{equation*}
such that when $G>G_{0}$, we have
\begin{equation*}
f(G)< \frac{(a_{11}+\check{m})(a_{22}+\check{m})}{a_{12}+\epsilon_{0}}G.
\end{equation*}
Define $C_{1}:=G$, $C_{2}:=\frac{a_{11}+\check{m}}{a_{12}+\epsilon_{0}}G$, $N_{1}:=\max\limits_{[0, t_{\max})\times\Omega_{0}}u(t,x)$, and $N_{2}:=\max\limits_{0\leq u\leq\max(C_{1},N_{1})}f'(u)$.

Let
\begin{equation}\label{WS-1}
P=(C_{1}-u)e^{-(a_{12}+N_{2})t}~\text{and}~Q=(C_{2}-v)e^{-(a_{12}+N_{2})t}.
\end{equation}
By regular calculation, it follows that
\begin{equation}\label{WS-2}
\begin{array}{ll}
&P_{t}-\frac{d_{1}}{\rho^{2}(t)}\Delta P+(a_{11}+\frac{n\dot{\rho}(t)}{\rho(t)}+a_{12}+N_{2})P-a_{12}Q\\
=&\Big[-u_{t}-(C_{1}-u)(a_{12}+N_{2})+\frac{d_{1}}{\rho^{2}(t)}\Delta u-a_{12}(C_{2}-v)\\
&+(a_{11}+\frac{n\dot{\rho}(t)}{\rho(t)}+a_{12}+N_{2})(C_{1}-u)\Big]e^{-(a_{12}+N_{2})t}\\
=&\Big[(a_{11}+\frac{n\dot{\rho}(t)}{\rho(t)})C_{1}-a_{12}C_{2}\Big]e^{-(a_{12}+N_{2})t}>0
\end{array}
\end{equation}
and
\begin{equation}\label{WS-3}
\begin{array}{ll}
&Q_{t}-\frac{d_{2}}{\rho^{2}(t)}\Delta Q+(a_{22}+\frac{n\dot{\rho}(t)}{\rho(t)}+a_{12}+N_{2})Q-f'(\zeta)P\\
=&\Big[-v_{t}-(C_{2}-v)(a_{12}+N_{2})+\frac{d_{2}}{\rho^{2}(t)}\Delta v-f'(\zeta)(C_{1}-u)\\
&+(a_{22}+\frac{n\dot{\rho}(t)}{\rho(t)}+a_{12}+N_{2})(C_{2}-v)\Big]e^{-(a_{12}+N_{2})t}\\
=&\Big[(a_{22}+\frac{n\dot{\rho}(t)}{\rho(t)})C_{2}-f(C_{1})\Big]e^{-(a_{12}+N_{2})t}>0,
\end{array}
\end{equation}
where $\zeta$ is between $C_{1}$ and $u$. Therefore, through the assumption that $g(u)>0$ for $u>0$ in (A1), and \eqref{WS-1}-\eqref{WS-3}, we have that
\begin{eqnarray}\label{abb}
\left\{
\begin{array}{ll}
P_{t}>\frac{d_{1}}{\rho^{2}(t)}\Delta P-(a_{11}+\frac{n\dot{\rho}(t)}{\rho(t)}+a_{12}+N_{2})P+a_{12}Q,\; &\, t\in(0^{+}, \tau], x\in\Omega_{0}, \\[2mm]
Q_{t}>\frac{d_{2}}{\rho^{2}(t)}\Delta Q-(a_{22}+\frac{n\dot{\rho}(t)}{\rho(t)}+a_{12}+N_{2})Q+f'(\zeta)P,\; &\, t\in(0^{+}, \tau], x\in\Omega_{0}, \\[2mm]
P(t,x)>0, Q(t,x)>0,\; &\, t\in(0, \tau], x\in\partial\Omega_{0},\\[2mm]
P(0,x)>0, Q(0,x)>0,\; &\,x\in\xoverline{\Omega}_{0},\\[2mm]
P(0^{+},x)>0, Q(0^{+},x)>0, \; &\, x\in\Omega_{0}.
\end{array} \right.
\end{eqnarray}
We claim that
\begin{equation*}
\min\bigg\{\min\limits_{[0, \tau]\times\bar{\Omega}_{0}}P(t,x), \min\limits_{[0, \tau]\times\bar{\Omega}_{0}}Q(t,x)\bigg\}>0.
\end{equation*}
Otherwise, there would exist $(x_{0}, t_{0})\in(0^{+}, \tau]\times\Omega_{0}$ such that
\begin{equation*}
\min\Big\{Q(t_{0},x_{0}), P(t_{0},x_{0})\Big\}\leq0.
\end{equation*}
Without loss of generality, suppose $\min\Big\{P(t_{0},x_{0}), Q(t_{0},x_{0})\Big\}=P(t_{0},x_{0})$. Then, it follows that
\begin{equation}\label{WS-4}
P_{t}(t_{0},x_{0})-\frac{d_{1}}{\rho^{2}(t)}\Delta P(t_{0},x_{0})\leq 0.
\end{equation}
On the other hand, we have that
\begin{equation*}
-(a_{11}+\frac{n\dot{\rho}(t)}{\rho(t)}+a_{12}+N_{2})P(t_{0},x_{0})+a_{12}Q(t_{0},x_{0})\geq0,
\end{equation*}
which shows a contradiction with \eqref{WS-4} and the first inequality in \eqref{abb}. We thus derive that $u(t,x)< C_{1}$ and $v(t,x)< C_{2}$ in $[0, \tau]\times\bar{\Omega}_{0}$.
Obviously, $u(\tau,x)$ and $v(\tau, x)$ as new initial values satisfy the conditions $\min\limits_{\bar{\Omega}_{0}}u(\tau,x)<C_{1}$ and $\min\limits_{\bar{\Omega}_{0}}v(\tau,x)<C_{2}$. Repeating the above procedure yields that $u(t,x)< C_{1}$ and $v(t,x)< C_{2}$ in $[\tau, 2\tau]\times\bar{\Omega}_{0}$. Step by step, we can get that $u(t,x)< C_{1}$ and $v(t,x)< C_{2}$ in $[0, \infty)\times\bar{\Omega}_{0}$.
As a result, $t_{max}=\infty$, that is, model \eqref{1-2-20} exists a unique global nonnegative solution.

In the end, we prove that the solution of model \eqref{1-2-20} satisfies $u, v>0$ for $(t,x)\in\mathds{R}^{+} \times\Omega_{0}$. The nonnegativity of the solution, the assumption of $f(u)\geq0$ in (A1), and the assumption of $0<\frac{g(u)}{u}\leq1$ in (A3) give that
\begin{eqnarray*}
\left\{
\begin{array}{ll}
\frac{\partial u}{\partial t}\geq\frac{d_{1}}{\rho^{2}(t)}\Delta u-\frac{n\dot{\rho}(t)}{\rho(t)}u-a_{11}u,\; &\, t\in(0^{+}, \tau], x\in\Omega_{0}, \\[2mm]
\frac{\partial v}{\partial t}\geq\frac{d_{2}}{\rho^{2}(t)}\Delta v-\frac{n\dot{\rho}(t)}{\rho(t)}v-a_{22}v,\; &\, t\in(0^{+}, \tau], x\in\Omega_{0}, \\[2mm]
u=0, v=0,\; &\, t\in (0, \tau], x\in\partial\Omega_{0},\\[2mm]
u(0,x)=u_{0}(x), v(0,x)=v_{0}(x),\; &\,x\in\bar{\Omega}_{0},\\[2mm]
u(0^{+},x)=g(u(0,x)), v(0^{+},x)=v(0,x), \; &\, x\in\Omega_{0}.
\end{array} \right.
\end{eqnarray*}
This together with strong maximum principle yields $u, v>0$ for $(t,x)\in\mathds{R}^{+} \times\Omega_{0}$. Similarly, one can obtain that $u, v>0$ for $t\in(\tau, 2\tau]$ and $x\in\Omega_{0}$. Step by step, we can get that $u, v>0$ for $t\in(0, \infty]$ and $x\in\Omega_{0}$. This ends the proof.
\end{proof}
\end{theorem}
\section{Existence of the principal eigenvalue}\label{Section-3}
The principal eigenvalue is fundamental to the analysis of the dynamical behaviors of model \eqref{1-2-20}. In \cite{2-1}, Ant\'{o}n et al. presented the theory of the principal eigenvalue for cooperative periodic-parabolic systems. In this section, we are going to overcome the obstacles arising from pulses and give the existence of the principal eigenvalue for problem with pulses.

For convenience of writing, we first introduce the notations as follows:
\begin{equation*}
\begin{aligned}
\mathds{X}:=\bigg\{ U=(u_{1},u_{2})\in\big[\mathbb{C}^{0,1}([0,\tau]\times\xoverline{\Omega}_{0})\big]^{2}: U=(0,0)~\forall (t.x)\in [0,\tau]\times\partial{\Omega}_{0},~~~~~~~~~~~~~~\\
~~~~~~~~~~~~~~~~~~~~~~~~~~~~~~~~~~u_{1}(0,x)=g'(0)u_{1}(\tau,x),~u_{2}(0,x)=u_{2}(\tau,x)~\forall x\in \xoverline{\Omega}_{0}\bigg\};
\end{aligned}
\end{equation*}
\begin{equation*}
\begin{aligned}
\mathds{X}^{+}:=\text{closure}\bigg\{ U\in\mathds{X}: U\gg 0~\forall (t.x)\in [0,\tau]\times{\Omega}_{0}, \frac{\partial U}{\partial n}\ll 0~\forall (t.x)\in [0,\tau]\times\partial\Omega_{0}\bigg\};
\end{aligned}
\end{equation*}
\begin{equation*}
\begin{aligned}
\overset{\text{o}} {\mathds{X}^{+}}:=\text{int}\bigg\{ U\in\mathds{X}: U\gg 0~\forall (t.x)\in [0,\tau]\times{\Omega}_{0}, \frac{\partial U}{\partial n}\ll 0~\forall (t.x)\in [0,\tau]\times\partial\Omega_{0}\bigg\},
\end{aligned}
\end{equation*}
where $U\gg 0$ and $\frac{\partial U}{\partial n}\ll 0$ denote $u_{1}, u_{2}>0$ and $\frac{\partial u_{1}}{\partial n}, \frac{\partial u_{2}}{\partial n}<0$, respectively. It is obvious that $\mathds{X}$ is a Banach space, and $\mathds{X}^{+}$ is the positive cone of $\mathds{X}$. In addition, the interior $\overset{\text{o}} {\mathds{X}^{+}}$ of $\mathds{X}^{+}$ is nonempty. Linearizing model \eqref{1-2-20} at equilibrium point $(0,0)$ gives the following linear model
\begin{eqnarray}\label{2-1}
\left\{
\begin{array}{ll}
\frac{\partial u}{\partial t}=\frac{d_{1}}{\rho^{2}(t)}\Delta u-\frac{n\dot{\rho}(t)}{\rho(t)}u-a_{11}u+a_{12}v,\; &\, t\in((k\tau)^{+}, (k+1)\tau], x\in\Omega_{0}, \\[2mm]
\frac{\partial v}{\partial t}=\frac{d_{2}}{\rho^{2}(t)}\Delta v-\frac{n\dot{\rho}(t)}{\rho(t)}v-a_{22}v+f'(0)u,\; &\, t\in((k\tau)^{+}, (k+1)\tau], x\in\Omega_{0}, \\[2mm]
u=0, v=0,\; &\, t\in\mathds{R}^{+}, x\in\partial\Omega_{0},\\[2mm]
u(0,x)=u_{0}(x), v(0,x)=v_{0}(x),\; &\,x\in\xoverline{\Omega}_{0},\\[2mm]
u((k\tau)^{+},x)=g'(0)u(k\tau,x), v((k\tau)^{+},x)=v((k\tau),x), \; &\, x\in\Omega_{0}, k=0,1,2,\cdots.
\end{array} \right.
\end{eqnarray}
Now, we take into account the following periodic eigenvalue problem
\begin{eqnarray}\label{2-2}
\left\{
\begin{array}{ll}
\frac{\partial \phi}{\partial t}=\frac{d_{1}}{\rho^{2}(t)}\Delta \phi-\frac{n\dot{\rho}(t)}{\rho(t)}\phi-a_{11}\phi+a_{12}\psi+\lambda \phi,\; &\, t\in(0^{+}, \tau], x\in\Omega_{0}, \\[2mm]
\frac{\partial \psi}{\partial t}=\frac{d_{2}}{\rho^{2}(t)}\Delta \psi-\frac{n\dot{\rho}(t)}{\rho(t)}\psi-a_{22}\psi+f'(0)\phi+\lambda \psi,\; &\, t\in(0^{+}, \tau], x\in\Omega_{0}, \\[2mm]
\phi=0, \psi=0,\; &\, t\in[0, \tau], x\in\partial\Omega_{0},\\[2mm]
\phi(0,x)=\phi(\tau,x), \psi(0,x)=\psi(\tau,x),\; &\,x\in\xoverline{\Omega}_{0},\\[2mm]
\phi(0^{+},x)=g'(0)\phi(0,x), \psi(0^{+},x)=\psi(0,x), \; &\, x\in\Omega_{0}.
\end{array} \right.
\end{eqnarray}

In order to remove the obstacles arising from the pulse, we present the following periodic eigenvalue problem:
\begin{eqnarray}\label{2-3}
\left\{
\begin{array}{ll}
\frac{\partial u_{1}}{\partial t}=\frac{d_{1}}{\rho^{2}(t)}\Delta u_{1}-\frac{n\dot{\rho}(t)}{\rho(t)}u_{1}-a_{11}u_{1}+a_{12}u_{2}+\lambda u_{1},\; &\, t\in(0, \tau], x\in\Omega_{0}, \\[2mm]
\frac{\partial u_{2}}{\partial t}=\frac{d_{2}}{\rho^{2}(t)}\Delta u_{2}-\frac{n\dot{\rho}(t)}{\rho(t)}u_{2}-a_{22}u_{2}+f'(0)u_{1}+\lambda u_{2},\; &\, t\in(0, \tau], x\in\Omega_{0}, \\[2mm]
u_{1}=0, u_{2}=0,\; &\, t\in[0, \tau], x\in\partial\Omega_{0},\\[2mm]
u_{1}(0,x)=g'(0)u_{1}(\tau,x), u_{2}(0,x)=u_{2}(\tau,x), \; &\, x\in\xoverline{\Omega}_{0}.
\end{array} \right.
\end{eqnarray}
In fact, problem \eqref{2-2} is equivalent to problem \eqref{2-3}, which is because we can take
\begin{equation*}
\phi(t,x)=u_{1}(t,x),~~~\forall (t,x)\in (0,\tau]\times\xoverline{\Omega}_{0},
\end{equation*}
\begin{equation*}
\psi(t,x)=u_{2}(t,x),~~~\forall (t,x)\in [0,\tau]\times\xoverline{\Omega}_{0},
\end{equation*}
\begin{equation*}
\phi(0,x)=\phi(\tau,x)=u_{1}(\tau,x)=1/g'(0)u_{1}(0,x),~~~\forall x \in\xoverline{\Omega}_{0}.
\end{equation*}
For problem \eqref{2-2}, at the pulse point $0^{+}$, we have that $\phi(0^{+},x)=g'(0)\phi(0,x)=u_{1}(0,x)~\forall x\in \xoverline{\Omega}_{0}$ and $\psi(0^{+},x)=u_{2}(0,x)~\forall x\in \xoverline{\Omega}_{0}$.

Based on the above preparation work, we will give the existence of the principal eigenvalue of problem \eqref{2-3}.
\begin{theorem}\label{theorem 2-1}
The periodic eigenvalue problem \eqref{2-3} has the principal eigenvalue $\lambda_{1}$ with an eigenfunction $(\phi(t,x), \psi(t,x))\gg 0$ in $[0,\tau]\times\Omega_{0}$.
\end{theorem}
\begin{proof}
For an arbitrary given vector $U(t,x)=\big(u_{1}(t,x), u_{2}(t,x)\big)\in \mathds{X}$, let
\begin{equation*}
\bar{u}_{1}=\sup_{[0,\tau]\times\bar{\Omega}_{0}}u_{1}(t,x),~~
\bar{u}_{2}=\sup_{[0,\tau]\times\bar{\Omega}_{0}}u_{2}(t,x),~~
\xoverline{\mathcal{C}}=\frac{\bar{u}_{1}+\bar{u}_{2}}{2},
\end{equation*}
and
\begin{equation*}
\begin{aligned}
\mathds{Y}=\Big\{ V=\big(v_{1}(x),v_{2}(x)\big)\in\big[\mathbb{C}(\xoverline{\Omega}_{0})\big]^{2}: 0\ll V\ll \xoverline{\mathcal{C}}, v_{1}(x)|_{x\in\partial\Omega_{0}}=v_{2}(x)|_{x\in\partial\Omega_{0}}=0\Big\},
\end{aligned}
\end{equation*}
where $0\ll V\ll \xoverline{\mathcal{C}}$ denotes $0<v_{1}(x),v_{2}(x)<\xoverline{\mathcal{C}}$. Now, for an arbitrary given vector $V\in \mathds{Y}$, we define
\begin{equation*}
\kappa(V)=(g'(0)z(\tau,x), w(\tau,x)),
\end{equation*}
where $(z(t,x),w(t,x))$ is the solution of the following problem
\begin{eqnarray}\label{2-4}
\left\{
\begin{array}{ll}
\frac{\partial z}{\partial t}=\frac{d_{1}}{\rho^{2}(t)}\Delta z-\Big(\frac{n\dot{\rho}(t)}{\rho(t)}+a_{11}+\varpi\Big)z+a_{12}w+u_{1},~~~~~~t\in(0, \tau], x\in\Omega_{0}, \\[2mm]
\frac{\partial w}{\partial t}=\frac{d_{2}}{\rho^{2}(t)}\Delta w-\Big(\frac{n\dot{\rho}(t)}{\rho(t)}+a_{22}+\varpi\Big)w+f'(0)z+u_{2},~t\in(0, \tau], x\in\Omega_{0}, \\[2mm]
\varpi=a_{12}\max\big\{g'(0), 1\big\}-\inf \limits_{t\in(0,\tau)}\frac{n\dot{\rho}(t)}{\rho(t)}+f'(0),\\[2mm]
z=0, w=0,~~~~~~~~~~~~~~~~~~~~~~~~~~~~~~~~~~~~~~~~~~~~~~~~~~~~~~~~ t\in[0, \tau], x\in\partial\Omega_{0},\\[2mm]
z(0,x)=v_{1}(x), w(0,x)=v_{2}(x), ~~~~~~~~~~~~~~~~~~~~~~~~~~~~x\in\xoverline{\Omega}_{0}.
\end{array} \right.
\end{eqnarray}
It is not hard to obtain that $(z(t,x),w(t,x))\in \mathds{X}^{1,2}_{p}\big((0,\tau]\times \Omega_{0}\big)\cap\mathbb{C}^{(1+\alpha)/2, 1+\alpha}\big([0,\tau]\times \xoverline{\Omega}_{0}\big)$ by classical $L^p$ theory and embedding theorem for any $p>0$ and $0<\alpha<1$. With the help of the comparison principle, we have that $\kappa(V)\in \mathds{Y}$ for $V\in \mathds{Y}$. Additionally, Sobolev embedding theorem yields that $\kappa$ is compact. This together with Schauder fixed point theorem give that $\kappa$ admits at least one fixed point $V\in \mathds{Y}$, and a solution $(z(t,x),w(t,x))$ of problem \eqref{2-4} with $V=\big(v_{1}(x),v_{2}(x)\big)$ as initial value. Hence, $(z(t,x),w(t,x))$ is a solution of problem
\begin{eqnarray}\label{2-5}
\left\{
\begin{array}{ll}
\frac{\partial \xi}{\partial t}=\frac{d_{1}}{\rho^{2}(t)}\Delta \xi-\Big(\frac{n\dot{\rho}(t)}{\rho(t)}+a_{11}+\varpi\Big)\xi+a_{12}\eta+u_{1},~~~~~~t\in(0, \tau], x\in\Omega_{0}, \\[2mm]
\frac{\partial \eta}{\partial t}=\frac{d_{2}}{\rho^{2}(t)}\Delta \eta-\Big(\frac{n\dot{\rho}(t)}{\rho(t)}+a_{22}+\varpi\Big)\eta+f'(0)\xi+u_{2},~~~t\in(0, \tau], x\in\Omega_{0}, \\[2mm]
\varpi=a_{12}\max\big\{g'(0), 1\big\}-\inf \limits_{t\in(0,\tau)}\frac{n\dot{\rho}(t)}{\rho(t)}+f'(0),\\[2mm]
\xi=0, \eta=0,~~~~~~~~~~~~~~~~~~~~~~~~~~~~~~~~~~~~~~~~~~~~~~~~~~~~~~~~ t\in[0, \tau], x\in\partial\Omega_{0},\\[2mm]
\xi(0,x)=g'(0)\xi(\tau,x), \eta(0,x)=\eta(\tau,x), ~~~~~~~~~~~~~~~~x\in\xoverline{\Omega}_{0},
\end{array} \right.
\end{eqnarray}
and $\big(z(t,x),w(t,x)\big)\in \mathds{X}\cap\mathbb{C}^{(1+\alpha)/2, 1+\alpha}\big([0,\tau]\times \xoverline{\Omega}_{0}\big)$.

Next, we will prove that the solution of problem \eqref{2-5} is unique. In order to prove this, we only need to prove that there is only the zero solution $(0,0)$ to the problem
\begin{eqnarray}\label{2-6}
\left\{
\begin{array}{ll}
\frac{\partial \xi}{\partial t}=\frac{d_{1}}{\rho^{2}(t)}\Delta \xi-\Big(\frac{n\dot{\rho}(t)}{\rho(t)}+a_{11}+\varpi\Big)\xi+a_{12}\eta,~~~~~~~~~~~~~~~t\in(0, \tau], x\in\Omega_{0}, \\[2mm]
\frac{\partial \eta}{\partial t}=\frac{d_{2}}{\rho^{2}(t)}\Delta \eta-\Big(\frac{n\dot{\rho}(t)}{\rho(t)}+a_{22}+\varpi\Big)\eta+f'(0)\xi,~~~~~~~~~~~~t\in(0, \tau], x\in\Omega_{0}, \\[2mm]
\varpi=a_{12}\max\big\{g'(0), 1\big\}-\inf \limits_{t\in(0,\tau)}\frac{n\dot{\rho}(t)}{\rho(t)}+f'(0),\\[2mm]
\xi(t,x)=0, \eta(t,x)=0,~~~~~~~~~~~~~~~~~~~~~~~~~~~~~~~~~~~~~~~~~~ t\in[0, \tau], x\in\partial\Omega_{0},\\[2mm]
\xi(0,x)=g'(0)\xi(\tau,x), \eta(0,x)=\eta(\tau,x), ~~~~~~~~~~~~~~~~~x\in\xoverline{\Omega}_{0},
\end{array} \right.
\end{eqnarray}
which is the homogeneous linear problem corresponding to problem \eqref{2-5}. On the contrary, assume that
\begin{equation*}
\min\bigg\{\min\limits_{[0,\tau]\times\bar{\Omega}_{0} }\xi(t,x), \min\limits_{[0,\tau]\times\bar{\Omega}_{0} }\eta(t,x)\bigg\}<0
\end{equation*}
or
\begin{equation*}
\max\bigg\{\max\limits_{[0,\tau]\times\bar{\Omega}_{0} }\xi(t,x), \max\limits_{[0,\tau]\times\bar{\Omega}_{0} }\eta(t,x)\bigg\}>0
\end{equation*}
and show a contradiction. Without loss of generality, suppose that there exists $(t_{0},x_{0})\in[0,\tau]\times\Omega_{0}$ such that
\begin{equation*}
\min\bigg\{\min\limits_{[0,\tau]\times\bar{\Omega}_{0} }\xi(t,x), \min\limits_{[0,\tau]\times\bar{\Omega}_{0} }\eta(t,x)\bigg\}=\xi(t_{0},x_{0})<0.
\end{equation*}
If $0<t_{0}\leq \tau$, then
\begin{equation*}
a_{12}\eta(t_{0},x_{0})-\Big(\frac{n\dot{\rho}(t_{0})}{\rho(t_{0})}+a_{11}+\varpi\Big)\xi(t_{0},x_{0})>0.
\end{equation*}
On the other hand, it is immediate to get that
\begin{equation*}
\frac{\partial \xi(t_{0},x_{0})}{\partial t}-\frac{d_{1}}{\rho^{2}(t_{0})}\Delta \xi(t_{0},x_{0})\leq 0,
\end{equation*}
which leads to a contradiction to the equation
\begin{equation}\label{2-7}
\frac{\partial \xi(\tau, x_{0})}{\partial t}=\frac{d_{1}}{\rho^{2}(\tau)}\Delta \xi(\tau, x_{0})-\Big(\frac{n\dot{\rho}(\tau)}{\rho(\tau)}+a_{11}+\varpi\Big)\xi(\tau, x_{0})+a_{12}\eta(\tau, x_{0}).
\end{equation}
If $t_{0}=0$, then \eqref{2-7} together with the boundary condition $\xi(0,x)=g'(0)\xi(\tau,x)$ and $\eta(0,x)=\eta(\tau,x)$ give that
\begin{equation*}
\frac{\partial \xi(0, x_{0})}{\partial t}-\frac{d_{1}}{\rho^{2}(\tau)}\Delta \xi(0, x_{0})=-\Big(\frac{n\dot{\rho}(\tau)}{\rho(\tau)}+a_{11}+\varpi\Big)\xi(0, x_{0})+a_{12}g'(0)\eta(0,x).
\end{equation*}
Since $\varpi>a_{12}\max\big\{g'(0), 1\big\}-\inf \limits_{t\in(0,\tau)}\frac{n\dot{\rho}(t)}{\rho(t)}$, it follows that
\begin{equation*}
-\Big(\frac{n\dot{\rho}(\tau)}{\rho(\tau)}+a_{11}+\varpi\Big)\xi(0, x_{0})+a_{12}g'(0)\eta(0,x)>0,
\end{equation*}
On the other hand, it follows that
\begin{equation*}
\frac{\partial \xi(0, x_{0})}{\partial t}-\frac{d_{1}}{\rho^{2}(\tau)}\Delta \xi(0, x_{0})\leq 0,
\end{equation*}
which leads to a contradiction.  Hence, the solution of problem \eqref{2-5} is unique.

Now, define the operator $\mathcal{A}U=V$. By the strong maximum principle and Hopf boundary lemma, we have that $\mathcal{A}$ is strongly positive with respect to $\mathds{X}$. Additionally, $\mathcal{A}$ is a linear compact operator for the reason that the imbedding $\mathbb{C}^{(1+\alpha)/2, 1+\alpha}\hookrightarrow \mathbb{C}^{0,1}$ is compact. According to a strong version of the Kerin-Rutman theorem \cite[Theorem 19.3]{2-2},  $r(\mathcal{A})$ is an algebraically simple eigenvalue of $\mathcal{A}$ with an eigenvector $V\in \text{Int}({\mathds{X}^{+}})$, and there is no other eigenvalue with a positive eigenvector. Therefore, problem \eqref{2-3} has the principal eigenvalue
\begin{equation*}
\lambda_{1}=\frac{1}{r(\mathcal{A})}-\varpi
\end{equation*}
 with an eigenfunction $V=\big(\phi(t,x), \psi(t,x)\big)\gg 0$ in $[0,\tau]\times\Omega_{0}$. This ends the proof.
\end{proof}
\begin{Remark}
For model \eqref{1-2-20}, if we also take into account the impact of the instantaneous phenomenon on the number of infected individuals, that is, $v\big((k\tau)^{+},x\big)=v\big((k\tau),x\big)$ is replaced by $v\big((k\tau)^{+},x\big)\\=h\big(v(k\tau,x)\big)$, the above procedure can also prove the existence of the corresponding principal eigenvalue.
\end{Remark}

\section{Dynamical behaviours of the model}\label{Section-4}
This section will explore the long-time dynamical behaviours of model \eqref{1-2-20} in terms of the principal eigenvalue. To this end, we begin with the
definition of upper and lower solutions of model \eqref{1-2-20} with periodic pulses, and then prove that the comparison principle is applicable for this model.
\begin{definition}\label{definition 3-1}
Functions $\big(\overline{u}(t,x), \overline{v}(t,x)\big)$, $\big(\underline{u}(t,x), \underline{v}(t,x)\big)\in\mathds{PC}_{1,2}^{2}\cap\mathds{PC}^{2}$ satisfying
\begin{equation*}
(0,0)\leq \big(\underline{u}(t,x), \underline{v}(t,x)\big)\leq \big(\overline{u}(t,x), \overline{v}(t,x)\big)
\end{equation*}
are called ordered upper and lower solutions of model \eqref{1-2-20}, respectively, if $\big(\overline{u}(t,x), \overline{v}(t,x)\big)$ and\\ $\big(\underline{u}(t,x), \underline{v}(t,x)\big)$ satisfy
\begin{eqnarray*}
\left\{
\begin{array}{ll}
\frac{\partial \overline{u}}{\partial t}\geq\frac{d_{1}}{\rho^{2}(t)}\Delta \overline{u}-\frac{n\dot{\rho}(t)}{\rho(t)}\overline{u}-a_{11}\overline{u}+a_{12}\overline{v},\; &\, t\in\big((k\tau)^{+}, (k+1)\tau\big], x\in\Omega_{0}, \\[2mm]
\frac{\partial \overline{v}}{\partial t}\geq\frac{d_{2}}{\rho^{2}(t)}\Delta \overline{v}-\frac{n\dot{\rho}(t)}{\rho(t)}\overline{v}-a_{22}\overline{v}+f(\overline{u}),\; &\, t\in\big((k\tau)^{+}, (k+1)\tau\big], x\in\Omega_{0}, \\[2mm]
\overline{u}\geq0, \overline{v}\geq0,\; &\, t>0, x\in\partial\Omega_{0},\\[2mm]
\overline{u}(0,x)\geq u_{0}(x), \overline{v}(0,x)\geq v_{0}(x),\; &\,x\in\xoverline{\Omega}_{0},\\[2mm]
\overline{u}\big((k\tau)^{+},x\big)\geq g\big(\overline{u}(k\tau,x)\big), \overline{v}\big((k\tau)^{+},x\big)\geq \overline{v}\big((k\tau),x\big), \; &\, x\in\Omega_{0}, k=0,1,2,\cdots,
\end{array} \right.
\end{eqnarray*}
and
\begin{eqnarray*}
\left\{
\begin{array}{ll}
\frac{\partial \underline{u}}{\partial t}\leq\frac{d_{1}}{\rho^{2}(t)}\Delta \underline{u}-\frac{n\dot{\rho}(t)}{\rho(t)}\underline{u}-a_{11}\underline{u}+a_{12}\underline{v},\; &\, t\in\big((k\tau)^{+}, (k+1)\tau\big], x\in\Omega_{0}, \\[2mm]
\frac{\partial \underline{v}}{\partial t}\leq\frac{d_{2}}{\rho^{2}(t)}\Delta \underline{v}-\frac{n\dot{\rho}(t)}{\rho(t)}\underline{v}-a_{22}\underline{v}+f(\underline{u}),\; &\, t\in\big((k\tau)^{+}, (k+1)\tau\big], x\in\Omega_{0}, \\[2mm]
\underline{u}=0, \underline{v}=0,\; &\, t>0, x\in\partial\Omega_{0},\\[2mm]
\underline{u}(0,x)\leq u_{0}(x), \underline{v}(0,x)\leq v_{0}(x),\; &\,x\in\xoverline{\Omega}_{0},\\[2mm]
\underline{u}\big((k\tau)^{+},x\big)\leq g\big(\underline{u}(k\tau,x)\big), \underline{v}\big((k\tau)^{+},x\big)\leq \underline{v}\big((k\tau),x\big), \; &\, x\in\Omega_{0}, k=0,1,2,\cdots,
\end{array} \right.
\end{eqnarray*}
respectively.
\end{definition}
\begin{lemma}\label{lemma 3-1}
If $\big(\overline{u}(t,x), \overline{v}(t,x)\big)$ and $\big(\underline{u}(t,x), \underline{v}(t,x)\big)$ are the ordered upper and lower solutions of model \eqref{1-2-20}, then the unique solution $\big(u(t,x), v(t,x)\big)$ of model \eqref{1-2-20} satisfies
\begin{equation*}
\big(\underline{u}(t,x), \underline{v}(t,x)\big)\leq \big(u(t,x), v(t,x)\big) \leq\big(\overline{u}(t,x), \overline{v}(t,x)\big), ~(t,x)\in[0,+\infty)\times\overline{\Omega}_{0}.
\end{equation*}
\end{lemma}
\begin{proof}
Let $(U, V)=\big(\xoverline{u}-u, \xoverline{v}-v\big)$. With the help of Assumption (A1), \autoref{definition 3-1}, and Lagrange mean value theorem, it follows that
\begin{eqnarray*}
\left\{
\begin{array}{ll}
\frac{\partial U}{\partial t}\geq\frac{d_{1}}{\rho^{2}(t)}\Delta U-\frac{n\dot{\rho}(t)}{\rho(t)}U-a_{11}U+a_{12}V,\; &\, t\in(0^{+}, \tau], x\in\Omega_{0}, \\[2mm]
\frac{\partial V}{\partial t}\geq\frac{d_{2}}{\rho^{2}(t)}\Delta V-\frac{n\dot{\rho}(t)}{\rho(t)}V-a_{22}V+f'(\xi)U,\; &\, t\in(0^{+}, \tau], x\in\Omega_{0}, \\[2mm]
U\geq0, V\geq0,\; &\, t\in\mathds{R}^{+}, x\in\partial\Omega_{0},\\[2mm]
U(0,x)\geq 0, V(0,x)\geq 0,\; &\,x\in\xoverline{\Omega}_{0},\\[2mm]
U\big(0^{+},x\big)\geq0, V\big(0^{+},x\big)\geq 0, \; &\, x\in\Omega_{0},
\end{array} \right.
\end{eqnarray*}
where for an arbitrary given $(t,x)\in(0^{+}, \tau]\times \Omega_{0}$,
\begin{equation*}
\min\big(\xoverline{u}, u\big)\leq\xi \leq  \max\big(\xoverline{u}, u\big).
\end{equation*}
The maximum principle for cooperative systems (see, e.g.,  \cite[Theorem 13 in Chapter 3]{2-3}) yields that $(U, V)\geq (0,0)$ for $(t,x)\in(0^{+}, \tau]\times\xoverline{\Omega}_{0}$. Next, we take $(U, V)$ at time $\tau^{+}$ as a new initial value for the period $(\tau^{+}, 2\tau]$, which satisfies
$U(\tau^{+},x)\geq0$ and $V(\tau^{+},x)\geq0$. The same procedure gives that $(U, V)\geq (0,0)$ for $(t,x)\in(\tau^{+}, 2\tau]\times\xoverline{\Omega}_{0}$. Step by step, we have that $(U, V)\geq (0,0)$ for $(t,x)\in[0,+\infty)\times\xoverline{\Omega}_{0}$, that is, $(u, v) \leq(\overline{u}, \overline{v})$ for $(t,x)\in[0,+\infty)\times\overline{\Omega}$. By the same method, the result that $(\underline{u}, \underline{v})\leq(u, v) $ for $(t,x)\in[0,+\infty)\times\overline{\Omega}$ can be also proved. This ends the proof.
\end{proof}

In what follows, with the help of \autoref{lemma 3-1}, the long-time dynamical behaviours of model \eqref{1-2-20} will be explored.
\begin{theorem}\label{theorem 3-1}
If $\lambda_{1}>0$, then the solution of model \eqref{1-2-20} satisfies
\begin{equation*}
\lim\limits_{t\rightarrow+\infty}(u(t,x),v(t,x))=(0,0) \text{~uniformly~for~}x\in\xoverline{\Omega}_{0}.
\end{equation*}
\begin{proof}
We first normalize the eigenfunction $\big(\phi(t,x), \psi(t,x)\big)$ corresponding to the principal eigenvalue $\lambda_{1}$, that is,
\begin{equation}\label{3-3}
\bar{\phi}(t,x)=\frac{1}{C}\phi(t,x), ~~\bar{\psi}(t,x)=\frac{1}{C}\psi(t,x),
\end{equation}
where
\begin{equation*}
C=\max\bigg\{\max\limits_{[0,\infty)\times\Omega_{0}} \phi(t,x), \max\limits_{[0,\infty)\times\Omega_{0} }\psi(t,x)\bigg\}.
\end{equation*}
For an arbitrary given initial value $(u_{0}(x), v_{0}(x) )$, there exists a sufficiently large $M$ such that for any $x\in\Omega_{0}$,
\begin{equation*}
M\bar{\phi}(0,x)\geq u_{0}(x)~\text{and}~M\bar{\psi}(0,x)\geq v_{0}(x).
\end{equation*}
Now, construct the functions
\begin{equation*}
\overline{u}(t,x)=Me^{-\lambda_{1}t}\bar{\phi}(t,x), ~~\overline{v}(t,x)=Me^{-\lambda_{1}t}\bar{\psi}(t,x).
\end{equation*}

We next will verify that $\big(\overline{u}(t,x), \overline{v}(t,x)\big)$ is the upper solution of model \eqref{1-2-20}. In virtue of Assumptions (A1) and (A2) and periodic eigenvalue problem \eqref{2-2}, the direct calculation yields that
\begin{equation*}
\begin{aligned}
&\frac{\partial \xoverline{u}}{\partial t}-\frac{d_{1}}{\rho^{2}(t)}\Delta \xoverline{u}+\frac{n\dot{\rho}(t)}{\rho(t)}\xoverline{u}+a_{11}\xoverline{u}-a_{12}\xoverline{v}\\
\geq&~Me^{-\lambda_{1}t}\bigg[\bar{\phi}_{t}-\lambda_{1}\bar{\phi}-\frac{d_{1}}{\rho^{2}(t)}\Delta\bar{\phi}+\frac{n\dot{\rho}(t)}{\rho(t)}\bar{\phi}+a_{11}\bar{\phi}-a_{12}\bar{\psi}\bigg]\\
\geq&~0
\end{aligned}
\end{equation*}
and
\begin{equation*}
\begin{aligned}
&\frac{\partial \xoverline{v}}{\partial t}-\frac{d_{2}}{\rho^{2}(t)}\Delta \xoverline{v}+\frac{n\dot{\rho}(t)}{\rho(t)}\xoverline{v}+a_{22}\xoverline{v}-f(\xoverline{u})\\
\geq &~\frac{\partial \xoverline{v}}{\partial t}-\frac{d_{2}}{\rho^{2}(t)}\Delta \xoverline{v}+\frac{n\dot{\rho}(t)}{\rho(t)}\xoverline{v}+a_{22}\xoverline{v}-f'(0)\xoverline{u}\\
\geq&~Me^{-\lambda_{1}t}\bigg[\bar{\psi}_{t}-\lambda_{1}\bar{\psi}-\frac{d_{2}}{\rho^{2}(t)}\Delta\bar{\psi}+\frac{n\dot{\rho}(t)}{\rho(t)}\bar{\psi}+a_{22}\bar{\psi}-f'(0)\bar{\phi}\bigg]\\
\geq &~0.
\end{aligned}
\end{equation*}
Additionally, at the pulse point, we have from Assumptions (A1) and (A3) that
\begin{equation*}
\begin{aligned}
&\xoverline{u}\big((n\tau)^{+}, x\big)\geq Me^{-\lambda_{1}n\tau}\bar{\phi}\big((n\tau)^{+},x\big)\\
\geq &~Mg'(0)e^{-\lambda_{1}n\tau}\bar{\phi}\big(n\tau,x\big)\\
\geq &~g'(0)\xoverline{u}\big(n\tau,x\big)\geq g\big(\xoverline{u}(n\tau,x)\big)
\end{aligned}
\end{equation*}
and
\begin{equation*}
\begin{aligned}
&\xoverline{v}\big((n\tau)^{+}, x\big)\geq Me^{-\lambda_{1}n\tau}\bar{\psi}\big((n\tau)^{+},x\big)\\
\geq&~ Me^{-\lambda_{1}n\tau}\bar{\psi}\big(n\tau,x\big)\geq \xoverline{v}(n\tau,x).
\end{aligned}
\end{equation*}
Therefore, it follows from \autoref{definition 3-1} that $\big(\overline{u}(t,x), \overline{v}(t,x)\big)$ is the upper solution of model \eqref{1-2-20}.

In virtue of \autoref{lemma 3-1}, it is not hard to get that
\begin{equation*}
(u, v) \leq(\overline{u}, \overline{v}), ~(t,x)\in[0,+\infty)\times\overline{\Omega}_{0}.
\end{equation*}
This together with the condition $\lambda_{1}>0$ yield that
\begin{equation*}
\lim\limits_{t\rightarrow+\infty}\big(u(t,x),v(t,x)\big)=(0,0), \text{~for~all~}x\in\xoverline{\Omega}_{0}.
\end{equation*}
This completes the proof.
\end{proof}
\end{theorem}
The above theorem presents the long-time dynamical behaviour of model \eqref{1-2-20} when $\lambda_{1}>0$. Next, we will discuss the case of when $\lambda_{1}<0$. To this end, we first introduce the following auxiliary periodic problem:
\begin{eqnarray}\label{3-1}
\left\{
\begin{array}{ll}
\frac{\partial u}{\partial t}=\frac{d_{1}}{\rho^{2}(t)}\Delta u-\frac{n\dot{\rho}(t)}{\rho(t)}u-a_{11}u+a_{12}v,\; &\, t\in(0^{+}, \tau], x\in\Omega_{0}, \\[2mm]
\frac{\partial v}{\partial t}=\frac{d_{2}}{\rho^{2}(t)}\Delta v-\frac{n\dot{\rho}(t)}{\rho(t)}v-a_{22}v+f(u),\; &\, t\in(0^{+}, \tau], x\in\Omega_{0}, \\[2mm]
u(t,x)=0, v(t,x)=0,\; &\, t\in[0, \tau], x\in\partial\Omega_{0}, \\[2mm]
u(0,x)=u(\tau,x), v(0,x)=v(\tau,x),\; &\,x\in\xoverline{\Omega}_{0},\\[2mm]
u(0^{+},x)=g(u(0,x)), v(0^{+},x)=v(0,x), \; &\, x\in\Omega_{0}.
\end{array} \right.
\end{eqnarray}

Similarly to \autoref{definition 3-1}, the definition of the ordered upper and lower solutions of auxiliary periodic problem \eqref{3-1} is given below.
\begin{definition}\label{definition 3-2}
Vector-valued functions $\big(\overline{u}(t,x), \overline{v}(t,x)\big)$, $\big(\underline{u}(t,x), \underline{v}(t,x)\big)\in\mathds{PC}_{1,2}^{2}\cap\mathds{PC}^{2}$ satisfying
\begin{equation*}
(0,0)\leq \big(\underline{u}(t,x), \underline{v}(t,x)\big)\leq \big(\overline{u}(t,x), \overline{v}(t,x)\big)
\end{equation*}
are the ordered upper solution and lower solution of problem \eqref{3-1}, respectively, if $\big(\overline{u}(t,x), \overline{v}(t,x)\big)$ and $\big(\underline{u}(t,x), \underline{v}(t,x)\big)$ satisfy
\begin{eqnarray*}
\left\{
\begin{array}{ll}
\frac{\partial \overline{u}}{\partial t}\geq\frac{d_{1}}{\rho^{2}(t)}\Delta \overline{u}-\frac{n\dot{\rho}(t)}{\rho(t)}\overline{u}-a_{11}\overline{u}+a_{12}\overline{v},\; &\, t\in\big(0^{+}, \tau\big], x\in\Omega_{0}, \\[2mm]
\frac{\partial \overline{v}}{\partial t}\geq\frac{d_{2}}{\rho^{2}(t)}\Delta \overline{v}-\frac{n\dot{\rho}(t)}{\rho(t)}\overline{v}-a_{22}\overline{v}+f(\overline{u}),\; &\, t\in\big(0^{+}, \tau\big], x\in\Omega_{0}, \\[2mm]
\overline{u}\geq0, \overline{v}\geq0,\; &\, t\in\big[0, \tau\big], x\in\partial\Omega_{0},\\[2mm]
\overline{u}(0,x)\geq \overline{u}(\tau,x), \overline{v}(0,x)\geq \overline{v}(\tau,x),\; &\,x\in\xoverline{\Omega}_{0},\\[2mm]
\overline{u}\big(0^{+},x\big)\geq g\big(\overline{u}(0,x)\big), \overline{v}\big(0^{+},x\big)\geq \overline{v}\big(0,x\big), \; &\, x\in\Omega_{0},
\end{array} \right.
\end{eqnarray*}
and
\begin{eqnarray*}
\left\{
\begin{array}{ll}
\frac{\partial \underline{u}}{\partial t}\leq\frac{d_{1}}{\rho^{2}(t)}\Delta \underline{u}-\frac{n\dot{\rho}(t)}{\rho(t)}\underline{u}-a_{11}\underline{u}+a_{12}\underline{v},\; &\, t\in\big(0^{+}, \tau\big], x\in\Omega_{0}, \\[2mm]
\frac{\partial \underline{v}}{\partial t}\leq\frac{d_{2}}{\rho^{2}(t)}\Delta \underline{v}-\frac{n\dot{\rho}(t)}{\rho(t)}\underline{v}-a_{22}\underline{v}+f(\underline{u}),\; &\, t\in\big(0^{+}, \tau\big], x\in\Omega_{0}, \\[2mm]
\underline{u}=0, \underline{v}=0,\; &\, t\in\big[0, \tau\big], x\in\partial\Omega_{0},\\[2mm]
\underline{u}(0,x)\leq \underline{u}(\tau,x), \underline{v}(0,x)\leq \underline{v}(\tau,x),\; &\,x\in\xoverline{\Omega}_{0},\\[2mm]
\underline{u}\big(0^{+},x\big)\leq g\big(\underline{u}(0,x)\big), \underline{v}\big(0^{+},x\big)\leq \underline{v}\big(0, x\big), \; &\, x\in\Omega_{0},
\end{array} \right.
\end{eqnarray*}
respectively.
\end{definition}

Based on \autoref{definition 3-2}, we next will construct corresponding monotone iteration sequence, which plays an important role in checking that problem \eqref{3-1} has a unique positive periodic solution when $\lambda_{1}<0$. For the convenience of writing, let
\begin{equation*}
H_{1}(t, u, v):=h_{1}(t, u, v)+m_{1}u, ~~H_{2}(t, u, v):=h_{2}(t, u, v)+m_{2}v,
\end{equation*}
where
\begin{eqnarray*}
\begin{array}{ll}
h_{1}(t, u, v):=a_{12}v-\frac{n\dot{\rho}(t)}{\rho(t)}u-a_{11}u,\; &\, m_{1}:=\sup \limits_{t\in(0,\tau)}\frac{n\dot{\rho}(t)}{\rho(t)}+\frac{3}{2}a_{11},\\[2mm]
h_{1}(t, u, v):=f(u)-\frac{n\dot{\rho}(t)}{\rho(t)}v-a_{22}v,\; &\, m_{2}:=\sup \limits_{t\in(0,\tau)}\frac{n\dot{\rho}(t)}{\rho(t)}+\frac{3}{2}a_{22}.
\end{array}
\end{eqnarray*}
Note that $H_{1}$ and $H_{2}$ are increasing about $u$ and $v$, respectively. Therefore, by choosing the upper solution $(\xoverline{u}, \xoverline{v})$ of auxiliary periodic problem \eqref{3-1} as an initial iteration $(\xoverline{u}^{(0)}, \xoverline{v}^{(0)})$, we can construct a sequence $\big\{(\xoverline{u}^{(i)}, \xoverline{v}^{(i)})\big\}_{i=0}^{n}$ from periodic problem
\begin{eqnarray}\label{3-2}
\left\{
\begin{array}{ll}
\frac{\partial \xoverline{u}^{(i)}}{\partial t}-\frac{d_{1}}{\rho^{2}(t)}\Delta \xoverline{u}^{(i)}+m_{1}\xoverline{u}^{(i)}=H_{1}(t, \xoverline{u}^{(i-1)}, \xoverline{v}^{(i-1)}),\; &\, t\in(0^{+}, \tau], x\in\Omega_{0}, \\[2mm]
\frac{\partial \xoverline{v}^{(i)}}{\partial t}-\frac{d_{2}}{\rho^{2}(t)}\Delta \xoverline{v}^{(i)}+m_{2}\xoverline{v}^{(i)}=H_{2}(t, \xoverline{u}^{(i-1)}, \xoverline{v}^{(i-1)}),\; &\, t\in(0^{+}, \tau], x\in\Omega_{0}, \\[2mm]
\xoverline{u}^{(i)}(t,x)=0, \xoverline{v}^{(i)}(t,x)=0,\; &\, t\in[0, \tau], x\in\partial\Omega_{0}, \\[2mm]
\xoverline{u}^{(i)}(0,x)=\xoverline{u}^{(i-1)}(\tau,x), \xoverline{v}^{(i)}(0,x)=\xoverline{v}^{(i-1)}(\tau,x),\; &\,x\in\xoverline{\Omega}_{0},\\[2mm]
\xoverline{u}^{(i)}(0^{+},x)=g(\xoverline{u}^{(i-1)}(\tau,x)), \xoverline{v}^{(i)}(0^{+},x)=\xoverline{v}^{(i-1)}(\tau,x), \; &\, x\in\Omega_{0}.
\end{array} \right.
\end{eqnarray}
Similarly, a sequence $\big\{(\underline{u}^{(i)}, \underline{v}^{(i)})\big\}_{i=0}^{n}$ can also be constructed by choosing the lower solution $(\underline{u}, \underline{v})$ of auxiliary periodic problem \eqref{3-1} as an initial iteration $(\underline{u}^{(0)}, \underline{v}^{(0)})$. In virtue of \cite[Lemma 3.1]{2-4}, the constructed sequences have the following property.
\begin{lemma}\label{lemma 3-2}
Let $\big\{(\xoverline{u}^{(i)}, \xoverline{v}^{(i)})\big\}_{i=0}^{n}$ and $\big\{(\underline{u}^{(i)}, \underline{v}^{(i)})\big\}_{i=0}^{n}$ be the sequences of iteration with $(\xoverline{u}, \xoverline{v})$ and $(\underline{u}, \underline{v})$ as the initial iteration, respectively. Then
\begin{equation*}
(\underline{u}, \underline{v})\leq(\underline{u}^{(i)}, \underline{v}^{(i)})\leq (\underline{u}^{(i+1)}, \underline{v}^{(i+1)})\leq (\xoverline{u}^{(i+1)}, \xoverline{v}^{(i+1)})\leq (\xoverline{u}^{(i)}, \xoverline{v}^{(i)})\leq(\xoverline{u}, \xoverline{v})
\end{equation*}
holds for all $(t,x)\in[0,\tau]\times\xoverline{\Omega}_{0}$.
\end{lemma}

Based on \autoref{lemma 3-2}, we now check that problem \eqref{3-1} has a unique positive periodic solution when $\lambda_{1}<0$.
\begin{theorem}\label{theorem 3-2}
If $\lambda_{1}<0$, then problem \eqref{3-1} has a unique positive periodic solution $(U(t,x), V(t,x))$.
\begin{proof}
Because the proof of this theorem is long, it is divided into three parts for readability.

\textbf{(1) The upper and lower solutions}

We first construct the upper solution of problem \eqref{3-1}. In virtue of Assumption (A2), it follows that there exists a sufficiently large positive constant $G$ such that
\begin{equation*}
f(G)\leq \frac{(a_{11}+\check{m})(a_{22}+\check{m})}{a_{12}}G.
\end{equation*}
Then, let $\xoverline{u}(t,x)=G$ and $\xoverline{v}(t,x)=\frac{a_{11}+\check{m}}{a_{12}}G$ for all $(t,x)\in[0,\tau]\times\xoverline{\Omega}_{0}$. Now, we will verify that $\big(\overline{u}(t,x), \overline{v}(t,x)\big)$ is the upper solution of problem \eqref{3-1}. The direct calculation yields that
\begin{equation*}
\begin{aligned}
\frac{\partial \xoverline{u}}{\partial t}-\frac{d_{1}}{\rho^{2}(t)}\Delta \xoverline{u}+\frac{n\dot{\rho}(t)}{\rho(t)}\xoverline{u}+a_{11}\xoverline{u}-a_{12}\xoverline{v}\geq0
\end{aligned}
\end{equation*}
and
\begin{equation*}
\begin{aligned}
\frac{\partial \xoverline{v}}{\partial t}-\frac{d_{2}}{\rho^{2}(t)}\Delta \xoverline{v}+\frac{n\dot{\rho}(t)}{\rho(t)}\xoverline{v}+a_{22}\xoverline{v}-f(\xoverline{u})\geq 0.
\end{aligned}
\end{equation*}
At the pulse point, Assumptions (A1) and (A3) give that
\begin{equation*}
\begin{aligned}
\xoverline{u}\big(0^{+}, x\big)-g\big(\xoverline{u}(0,x)\big)=G-g\big(G\big)\geq G\big(1-g'(0)\big)\geq 0
\end{aligned}
\end{equation*}
and
\begin{equation*}
\begin{aligned}
&\xoverline{v}\big(0^{+}, x\big)- \xoverline{v}(0,x)=0.
\end{aligned}
\end{equation*}
Additionally, it is not hard to see that both the boundary condition and the period condition that are required in \autoref{definition 3-2} are satisfied. Therefore, $\big(\overline{u}(t,x), \overline{v}(t,x)\big)$ is the upper solution of problem \eqref{3-1}.

Next, we construct the lower solution of problem \eqref{3-1}. For an arbitrary given negative number $\lambda_{1}$, there exists a sufficiently small positive number $\alpha$ such that $\lambda_{1}+\alpha<0$. Let
\begin{eqnarray*}
\underline{u}(t,x)=
\left\{
\begin{array}{ll}
\epsilon \bar{\phi}(t,x),\; &\, t=0, x\in\xoverline{\Omega}_{0}, \\[2mm]
\epsilon e^{(\lambda_{1}+\alpha)\tau} \xoverline{\phi}(t,x),\; &\, t=0^{+}, x\in\xoverline{\Omega}_{0}, \\[2mm]
\epsilon e^{(\lambda_{1}+\alpha)(\tau-t)}\xoverline{\phi}(t,x),\; &\, t\in(0^{+}, \tau], x\in\xoverline{\Omega}_{0},
\end{array} \right.
\end{eqnarray*}
and
\begin{eqnarray*}
\underline{v}(t,x)=
\left\{
\begin{array}{ll}
\epsilon \bar{\psi}(t,x),\; &\, t=0, x\in\xoverline{\Omega}_{0}, \\[2mm]
\epsilon e^{(\lambda_{1}+\alpha)\tau} \xoverline{\psi}(t,x),\; &\, t=0^{+}, x\in\xoverline{\Omega}_{0}, \\[2mm]
\epsilon e^{(\lambda_{1}+\alpha)(\tau-t)}\xoverline{\psi}(t,x),\; &\, t\in(0^{+}, \tau], x\in\xoverline{\Omega}_{0},
\end{array} \right.
\end{eqnarray*}
where $\epsilon$ is a positive constant that is small sufficiently. In virtue of Assumption (A4) and problem \eqref{2-2}, the direct calculation yields that
\begin{equation*}
\begin{aligned}
&\frac{\partial \underline{u}}{\partial t}-\frac{d_{1}}{\rho^{2}(t)}\Delta \underline{u}+\frac{n\dot{\rho}(t)}{\rho(t)}\underline{u}+a_{11}\underline{u}-a_{12}\underline{v}\\
\leq&~\epsilon e^{(\lambda_{1}+\alpha)(\tau-t)}\Bigg[\xoverline{\phi}_{t}-(\lambda_{1}+\alpha) \xoverline{\phi}-\frac{d_{1}}{\rho^{2}(t)}\Delta \xoverline{\phi} +\frac{n\dot{\rho}(t)}{\rho(t)} \xoverline{\phi}+a_{11} \xoverline{\phi}-a_{12}\bar{\psi}\Bigg] \\
\leq&-\epsilon\alpha e^{(\lambda_{1}+\alpha)(\tau-t)}\bar{\phi}\leq0
\end{aligned}
\end{equation*}
and
\begin{equation*}
\begin{aligned}
&\frac{\partial \underline{v}}{\partial t}-\frac{d_{2}}{\rho^{2}(t)}\Delta \underline{v}+\frac{n\dot{\rho}(t)}{\rho(t)}\underline{v}+a_{22}\underline{v}-f(\underline{u})\\
\leq&~\frac{\partial \underline{v}}{\partial t}-\frac{d_{2}}{\rho^{2}(t)}\Delta \underline{v}+\frac{n\dot{\rho}(t)}{\rho(t)}\underline{v}+a_{22}\underline{v}-f'(0)\underline{u}+H_{1}(\underline{u})^{\kappa_{1}}\\
\leq&~\epsilon e^{(\lambda_{1}+\alpha)(\tau-t)}\bigg[\xoverline{\psi}_{t}-(\lambda_{1}+\alpha) \xoverline{\psi}-\frac{d_{1}}{\rho^{2}(t)}\Delta \xoverline{\psi} +\frac{n\dot{\rho}(t)}{\rho(t)} \xoverline{\psi}+a_{22} \xoverline{\psi}-f'(0)\xoverline{\phi}\\
&+H_{1}\Big(\epsilon e^{(\lambda_{1}+\alpha)(\tau-t)}\Big)^{\kappa_{1}-1}\xoverline{\phi}^{\kappa_{1}}\bigg] \\
\leq&~\epsilon e^{(\lambda_{1}+\alpha)(\tau-t)}\bigg[H_{1}\Big(\epsilon e^{(\lambda_{1}+\alpha)(\tau-t)}\Big)^{\kappa_{1}-1}\xoverline{\phi}^{\kappa_{1}}-\alpha\xoverline{\psi}\bigg]\leq0.
\end{aligned}
\end{equation*}
At the pulse point, Assumption (A4) gives that
\begin{equation*}
\begin{aligned}
&\underline{u}\big(0^{+}, x\big)- g\big(\underline{u}(0,x)\big)\leq\underline{u}\big(0^{+}, x\big)-g'(0)\underline{u}(0,x)+H_{2}\big(\underline{u}(0,x)\big)^{\kappa_{2}}\\
\leq&~\epsilon\xoverline{\phi}(0,x)\Big[g'(0)\Big(e^{(\lambda_{1}+\alpha)\tau}-1\Big)+H_{2}\big(\epsilon\xoverline{\phi}(0,x)\big)^{\kappa_{2}-1}\Big]\leq 0
\end{aligned}
\end{equation*}
and
\begin{equation*}
\begin{aligned}
&\underline{v}\big(0^{+}, x\big)-\underline{v}(0,x)=\epsilon e^{(\lambda_{1}+\alpha)\tau}\xoverline{\psi}(0^{+},x)-\epsilon \bar{\psi}(0,x)\\
\leq&~\epsilon\xoverline{\psi}(0,x)\big[e^{(\lambda_{1}+\alpha)\tau}-1\big]\leq 0.
\end{aligned}
\end{equation*}
Moreover, both the boundary condition and the period condition that are required in \autoref{definition 3-2} are satisfied. Hence, $\big(\underline{u}(t,x), \underline{v}(t,x)\big)$ is the lower solution of problem \eqref{3-1}. As a summary, $\big(\overline{u}(t,x), \overline{v}(t,x)\big)$ and $\big(\underline{u}(t,x), \underline{v}(t,x)\big)$ are the upper solution and lower solution of problem \eqref{3-1}.

\textbf{(2) The existence of the positive periodic solution}

By choosing the upper and lower solutions constructed in \textbf{(1)} as the initial iterations, respectively, two iteration sequences $\big\{(\xoverline{u}^{(i)}, \xoverline{v}^{(i)})\big\}_{i=0}^{n}$ and $\big\{(\underline{u}^{(i)}, \underline{v}^{(i)})\big\}_{i=0}^{n}$ can be obtained. Then, \autoref{lemma 3-2} yields that
\begin{equation*}
(\underline{u}, \underline{v})\leq(\underline{u}^{(n-1)}, \underline{v}^{(n-1)})\leq (\underline{u}^{(n)}, \underline{v}^{(n)})\leq (\xoverline{u}^{(n)}, \xoverline{v}^{(n)})\leq (\xoverline{u}^{(n-1)}, \xoverline{v}^{(n-1)})\leq(\xoverline{u}, \xoverline{v}).
\end{equation*}
This together with the monotone bounded convergence theorem give that
\begin{equation}\label{3-5}
\lim\limits_{n\rightarrow+\infty}(\xoverline{u}^{(n)}, \xoverline{v}^{(n)})= (\xoverline{U}, \xoverline{V})\text{~~and~~}\lim\limits_{n\rightarrow+\infty}(\underline{u}^{(n)}, \underline{v}^{(n)})= (\underline{U}, \underline{V}),
\end{equation}
where $(\xoverline{U}, \xoverline{V})$ and $(\underline{U}, \underline{V})$ be the solutions of auxiliary periodic problem \eqref{3-1}.

Next, we prove that $(\underline{U}, \underline{V})$ and $(\xoverline{U}, \xoverline{V})$ are, respectively, the minimal and maximal positive periodic solutions of auxiliary periodic problem \eqref{3-1}. For the sake of distinction, denote by $\big(\tilde{u}(t,x), \tilde{v}(t,x)\big)$ the positive periodic solution of problem \eqref{3-1}. Based on periodic problem \eqref{3-2}, iteration sequence $\big\{(\tilde{u}^{(i)}, \tilde{v}^{(i)})\big\}_{i=0}^{n}$ can be obtained by choosing the upper solution $\big(\tilde{u}(t,x), \tilde{v}(t,x)\big)$ of auxiliary periodic problem \eqref{3-1} as the initial iteration. With the help of \autoref{lemma 3-2} and \eqref{3-5}, it follows that
\begin{equation*}
(\underline{u}, \underline{v})\leq(\underline{u}^{(n-1)}, \underline{v}^{(n-1)})\leq (\underline{u}^{(n)}, \underline{v}^{(n)})\leq (\underline{U}, \underline{V})\leq (\tilde{u}^{(n)}, \tilde{v}^{(n)})\leq(\tilde{u}, \tilde{v}).
\end{equation*}
This implies that $(\tilde{u}, \tilde{v})\geq(\underline{U}, \underline{V})$, that is, $(\underline{U}, \underline{V})$ is the minimal solution of problem \eqref{3-1}. Similarly, the conclusion that $(\xoverline{U}, \xoverline{V})$ is the maximal periodic solution of problem \eqref{3-1} can also be proved by regarding $\big(\tilde{u}(t,x), \tilde{v}(t,x)\big)$ as an lower solution of problem \eqref{3-1}. As a summary, $(\underline{U}, \underline{V})$ and $(\xoverline{U}, \xoverline{V})$ respectively be the minimal and maximal positive periodic solutions of auxiliary periodic problem \eqref{3-1}.

\textbf{(3) The uniqueness of the positive periodic solution}

Assume that $(u, v)$ and $(z, w)$ are two distinct solutions of problem \eqref{3-1}. Define a set
\begin{equation*}
\Gamma:=\Big\{r\in[0,1]:r(u, v)\leq (z, w), (t,x)\in[0,\tau] \times \xoverline{\Omega}_{0}\Big\}.
\end{equation*}
Now, we prove that $1\in \Gamma$. On the contrary, suppose that $r_{0}:=\sup\Gamma<1$. Without loss of generality, suppose that
\begin{equation}\label{3-6}
r_{0}=\sup \Big\{r\in[0,1]:ru\leq z, (t,x)\in[0,\tau] \times \xoverline{\Omega}_{0}\Big\}.
\end{equation}
This suppose means that
\begin{equation}\label{3-7}
r_{0}\leq \sup \Big\{r\in[0,1]:rv\leq w, (t,x)\in[0,\tau] \times \xoverline{\Omega}_{0}\Big\}.
\end{equation}
Based on auxiliary periodic problem \eqref{3-1}, the direct calculation yields that
\begin{equation}\label{3-8}
\begin{aligned}
&\frac{\partial }{\partial t}(z-r_{0}u)-\frac{d_{1}}{\rho^{2}(t)}\Delta (z-r_{0}u) +m_{1}(z-r_{0}u)\\
=&~\bigg(m_{1}-a_{11}-\frac{n\dot{\rho}(t)}{\rho(t)}\bigg)(z-r_{0}u)+a_{12}(w-r_{0}v)\geq 0.
\end{aligned}
\end{equation}
At the pulse point, Assumptions (A1) and (A2) give that
\begin{equation*}
\begin{aligned}
&z\big(0^{+}, x\big)-r_{0}u\big(0^{+}, x\big)=~g\Big(z\big(0, x\big)\Big)-r_{0}g\Big(u\big(0, x\big)\Big)\\
\geq&~g\Big(r_{0}u\big(0, x\big)\Big)-r_{0}g\Big(u\big(0, x\big)\Big)\geq 0.
\end{aligned}
\end{equation*}
Moreover, it is not hard to get that
\begin{equation*}
z(t, x)-r_{0}u(t, x)=0\text{~for~all~}(t,x)\in[0,\tau]\times\partial\Omega_{0}.
\end{equation*}
With the help of the strong maximum principle  (see, e.g.,  \cite[Theorem 4 in Chapter 3]{2-3}), one of the following two conclusions must hold:

(1) $z-r_{0}u>0$ holds for $(t,x)\in\big[0^{+}\cup(0^{+},\tau]\big]\times\Omega_{0}$. Note that $z$ and $u$ are $\tau$-periodic solutions, that is, $z(0,x)=z(\tau,x)$ and $u(0,x)=u(\tau,x)$ for $x\in \Omega_{0}$. This implies that $z-r_{0}u>0$ holds for $(t,x)\in[0,\tau]\times\Omega_{0}$. Additionally, according to the Hopf boundary lemma (see, e.g.,  \cite[Theorem 4 in Chapter 3]{2-3}), it follows that
\begin{equation*}
\frac{\partial }{\partial n}(z-r_{0}u)<0\text{~for~all~}(t,x)\in[0,\tau]\times\partial\Omega_{0}.
\end{equation*}
Therefore, there exists a sufficiently small $\epsilon_{0}$ such that $z-r_{0}u>\epsilon_{0} u$. This shows a  contradiction with \eqref{3-6}.

(2) $z-r_{0}u\equiv0$ holds for $(t,x)\in\big[0^{+}\cup(0^{+},\tau]\big]\times\Omega_{0}$. Then, it follows from \eqref{3-8} that $r_{0}v=w$ for $(t,x)\in[0,\tau]\times\Omega_{0}$. Moreover, it is not hard to get that
\begin{equation*}
w(t, x)-r_{0}v(t, x)=0\text{~for~all~}(t,x)\in[0,\tau]\times\partial\Omega_{0}.
\end{equation*}
Hence, we have that $r_{0}v=w$ for $(t,x)\in[0,\tau]\times\xoverline\Omega_{0}$. This shows a  contradiction with \eqref{3-7}.

As a summary, $1\in \Gamma$, that is, $(u, v)\leq (z, w)$ for $(t,x)\in[0,\tau] \times \xoverline{\Omega}_{0}$. Similarly, the conclude that $(z, w)\leq(u, v)$ for $(t,x)\in[0,\tau] \times \xoverline{\Omega}_{0}$ can also be proved. Therefore, $(u, v)= (z, w)$ for $(t,x)\in[0,\tau] \times \xoverline{\Omega}_{0}$, that is, the solution of  problem \eqref{3-1} is unique. This ends the proof.
\end{proof}
\end{theorem}

\autoref{theorem 3-2} tells us that problem \eqref{3-1} has a unique positive periodic solution $(U(t,x), V(t,x))$. In the following, we will show that the solution to model \eqref{1-2-20} converges to $(U(t,x), V(t,x))$.

\begin{theorem}\label{theorem 3-3}
If $\lambda_{1}<0$, then for any initial value $\big(u_{0}(x), v_{0}(x)\big)$ the solution of model \eqref{1-2-20} satisfies
\begin{equation}\label{3-13}
\lim\limits_{m\rightarrow+\infty}\big(u(t+m\tau,x),v(t+m\tau,x)\big)=\big(U(t,x), V(t,x)\big),~(t,x)\in[0,\infty)\times\xoverline\Omega_{0}.
\end{equation}
\begin{proof}
First, search for the upper and lower solutions of model \eqref{1-2-20}. Without loss of generality, suppose that $\big(u_{0}(x), v_{0}(x)\big)>(0,0)$ for all $x\in\Omega_{0}$. Otherwise, it follows from the strong maximum principle that there exists $t_{0}\in[0,\tau]$ such that $\big(u_{0}(x), v_{0}(x)\big)>(0,0)$ for all $x\in\Omega_{0}$. This supposition together with the Hopf boundary lemma yield that there exists a sufficiently small $\epsilon$ and a sufficiently large $G$ such that
\begin{equation}\label{3-12}
\big(\underline{u}(0,x), \underline{v}(0,x)\big)\leq\big(u_{0}(x), v_{0}(x)\big)\leq\big(\xoverline{u}(0,x), \xoverline{v}(0,x)\big),~x\in\xoverline\Omega_{0},
\end{equation}
where $\big(\underline{u}, \underline{v}\big)$ and $\big(\xoverline{u}, \xoverline{v}\big)$ are defined in the proof of \autoref{theorem 3-2}. With the help of \eqref{3-12} and  Assumptions (A1), it follows that
\begin{equation*}
\big(\underline{u}(0^{+},x), \underline{v}(0^{+},x)\big)\leq\big(u(0^{+}, x), v(0^{+}, x)\big)\leq\big(\xoverline{u}(0^{+},x), \xoverline{v}(0^{+},x)\big),~x\in\xoverline\Omega_{0}.
\end{equation*}
Then, \autoref{lemma 3-1} gives that
\begin{equation*}
\big(\underline{u}, \underline{v}\big)\leq\big(u, v\big)\leq\big(\xoverline{u}, \xoverline{v}\big),~(t,x)\in[0, \tau]\times\xoverline\Omega_{0}.
\end{equation*}
Based on the mathematical induction, it follows that
\begin{equation}\label{3-9}
\big(\underline{u}^{(0)}, \underline{v}^{(0)}\big)\leq\big(u, v\big)\leq\big(\xoverline{u}^{(0)}, \xoverline{v}^{(0)}\big),~(t,x)\in[0, \infty)\times\xoverline\Omega_{0}.
\end{equation}
This means that $\big(\underline{u}^{(0)}, \underline{v}^{(0)}\big)$ and $\big(\xoverline{u}^{(0)}, \xoverline{v}^{(0)}\big)$ are the upper and lower solutions of model \eqref{1-2-20}.

Next, prove that
\begin{equation}\label{3-10}
\big(\underline{u}^{(m)}(t,x), \underline{v}^{(m)}(t,x)\big)\leq\big(u(t+m\tau,x), v(t+m\tau,x)\big)\leq\big(\xoverline{u}^{(m)}(t,x), \xoverline{v}^{(m)}(t,x)\big)
\end{equation}
for $t\geq0$ and $x\in\xoverline\Omega_{0}$, where $m$ is an arbitrary positive integer. According to \eqref{3-9}, it follows that \eqref{3-10} holds for $m=0$. Now, assume that \eqref{3-10} holds for $m=n$, that is,
 \begin{equation}\label{3-11}
\big(\underline{u}^{(n)}(t,x), \underline{v}^{(n)}(t,x)\big)\leq\big(u(t+n\tau,x), v(t+n\tau,x)\big)\leq\big(\xoverline{u}^{(n)}(t,x), \xoverline{v}^{(n)}(t,x)\big)
\end{equation}
for $t\geq0$ and $x\in\xoverline\Omega_{0}$. Next, we will prove that \eqref{3-10} holds for $m=n+1$. In virtue of \eqref{3-11}, one can get that
 \begin{equation*}
\big(\underline{u}^{(n+1)}(0,x), \underline{v}^{(n+1)}(0,x)\big)\leq\big(u((n+1)\tau,x), v((n+1)\tau,x)\big)\leq\big(\xoverline{u}^{(n+1)}(0,x), \xoverline{v}^{(n+1)}(0,x)\big).
\end{equation*}
This together with Assumption (A1) yield that
 \begin{equation*}
\big(\underline{u}^{(n+1)}(0^{+},x), \underline{v}^{(n+1)}(0^{+},x)\big)\leq\big(u(((n+1)\tau)^{+},x), v(((n+1)\tau)^{+},x)\big)\leq\big(\xoverline{u}^{(n+1)}(0^{+},x), \xoverline{v}^{(n+1)}(0^{+},x)\big).
\end{equation*}
Then, \autoref{lemma 3-1} gives that
\begin{equation*}
\big(\underline{u}^{(n+1)}(t,x), \underline{v}^{(n+1)}(t,x)\big)\leq\big(u(t+(n+1)\tau,x), v(t+(n+1)\tau,x)\big)\leq\big(\xoverline{u}^{(n+1)}(t,x), \xoverline{v}^{(n+1)}(t,x)\big)
\end{equation*}
for $t\in[0,\tau]$ and $x\in\xoverline\Omega_{0}$. Using a method similar to obtaining \eqref{3-9}, it follows that
\begin{equation*}
\big(\underline{u}^{(n+1)}(t,x), \underline{v}^{(n+1)}(t,x)\big)\leq\big(u(t+(n+1)\tau,x), v(t+(n+1)\tau,x)\big)\leq\big(\xoverline{u}^{(n+1)}(t,x), \xoverline{v}^{(n+1)}(t,x)\big)
\end{equation*}
for $t\geq0$ and $x\in\xoverline\Omega_{0}$. Therefore, \eqref{3-10} holds for $m=n+1$. With the help of mathematical induction, the conclusion holds.

Finally, show that \eqref{3-13} holds. \autoref{theorem 3-2} gives that
\begin{equation}\label{3-14}
\lim\limits_{m\rightarrow+\infty}\big(\underline{u}^{(m)}(t,x), \underline{v}^{(m)}(t,x)\big)=\lim\limits_{m\rightarrow+\infty}\big(\xoverline{u}^{(n+1)}(t,x), \xoverline{v}^{(n+1)}(t,x)\big)=\big(U(t,x), V(t,x)\big)
\end{equation}
for $t\in[0,\tau]$ and $x\in\xoverline\Omega_{0}$. Using mathematical induction again, it follows that \eqref{3-14} holds for all $(t,x)\in[0, \infty]\times\xoverline\Omega_{0}$.
Therefore, we have that \eqref{3-13} holds. This completes the proof.
\end{proof}
\end{theorem}
The long-time dynamical behaviours of model \eqref{1-2-20} are investigated in this section by using the principal eigenvalue as a judgement condition. However, the relevant conclusions do not tell how the period pulses and evolving region will affect the dynamical behaviours of model \eqref{1-2-20}. This issue is investigated in the following section.
\section{Effects of initial region and impulsive intensity}\label{Section-5}
This section will consider the effects on the principal eigenvalue $\lambda_{1}$ of the initial region $\Omega_{0}$ and the impulsive intensity $g'(0)$, and then give some estimations on the principal eigenvalue $\lambda_{1}$ in special situations. To this end, we first give the adjoint problem of problem \eqref{2-2}, i.e.,
\begin{eqnarray}\label{4-1}
\left\{
\begin{array}{ll}
-\frac{\partial \zeta}{\partial t}=\frac{d_{1}}{\rho^{2}(t)}\Delta \zeta-\frac{n\dot{\rho}(t)}{\rho(t)}\zeta-a_{11}\zeta+f'(0)\eta+\mu \zeta,\; &\, t\in(0^{+}, \tau], x\in\Omega_{0}, \\[2mm]
-\frac{\partial \eta}{\partial t}=\frac{d_{2}}{\rho^{2}(t)}\Delta \eta-\frac{n\dot{\rho}(t)}{\rho(t)}\eta-a_{22}\eta+a_{12}\zeta+\mu \eta,\; &\, t\in(0^{+}, \tau], x\in\Omega_{0}, \\[2mm]
\zeta(t,x)=0, \eta(t,x)=0,\; &\, t\in[0, \tau], x\in\partial\Omega_{0},\\[2mm]
\zeta(0,x)=\zeta(\tau,x), \eta(0,x)=\eta(\tau,x),\; &\,x\in\xoverline{\Omega}_{0},\\[2mm]
\zeta(0^{+},x)=\frac{1}{g'(0)}\zeta(0,x), \eta(0^{+},x)=\eta(0,x), \; &\, x\in\Omega_{0}.
\end{array} \right.
\end{eqnarray}

For the convenience of later discussion, the association between the principal eigenvalue $\lambda_{1}$ of problem \eqref{2-2} and the principal eigenvalue $\mu_{1}$ of problem \eqref{4-1} is given next.
\begin{theorem}\label{theorem 4-1}
The principal eigenvalues of problem \eqref{2-2} and problem \eqref{4-1} satisfy $\lambda_{1}=\mu_{1}$.
\begin{proof}
By multiplying the first equations of problem \eqref{2-2} and problem \eqref{4-1} by $\zeta$ and  $\phi$, respectively, it follows that
\begin{eqnarray}\label{4-2}
\left\{
\begin{array}{ll}
\frac{\partial \phi}{\partial t}\zeta=\frac{d_{1}}{\rho^{2}(t)}\Delta \phi\zeta-\frac{n\dot{\rho}(t)}{\rho(t)}\phi\zeta-a_{11}\phi\zeta+a_{12}\psi\zeta+\lambda_{1} \phi \zeta, \\[2mm]
\frac{\partial \zeta}{\partial t}\phi=-\frac{d_{1}}{\rho^{2}(t)}\Delta \zeta\phi+\frac{n\dot{\rho}(t)}{\rho(t)}\zeta\phi+a_{11}\zeta\phi-f'(0)\eta\phi-\mu_{1} \zeta\phi ,
\end{array}
\right.
\end{eqnarray}
Then, integrating both sides of the first and second equations of \eqref{4-2} on $(0^{+}, \tau]\times\Omega_{0}$, and adding the obtained results give that
\begin{equation}\label{4-3}
\begin{aligned}
(\lambda_{1}-\mu_{1})\int_{\Omega_{0}}\int^{\tau}_{0^{+}}\phi \zeta dtdx=&\int_{\Omega_{0}}\int^{\tau}_{0^{+}}(f'(0)\eta\phi-a_{12}\psi\zeta)dtdx+\int_{\Omega_{0}}\int^{\tau}_{0^{+}}(\frac{\partial \phi}{\partial t}\zeta+\frac{\partial \zeta}{\partial t}\phi)dtdx\\
&+\int_{\Omega_{0}}\int^{\tau}_{0^{+}}\frac{d_{1}}{\rho^{2}(t)}(\Delta \zeta\phi- \Delta\phi\zeta)dtdx.
\end{aligned}
\end{equation}
The direct calculation further yields that
\begin{equation}\label{4-4}
\begin{aligned}
\int_{\Omega_{0}}\int^{\tau}_{0^{+}}(\frac{\partial \phi}{\partial t}\zeta+\frac{\partial \zeta}{\partial t}\phi)dtdx=&\int_{\Omega_{0}}\bigg(\zeta\phi|^{\tau}_{0^{+}}-\int^{\tau}_{0^{+}}\frac{\partial \zeta}{\partial t}\phi dt\bigg)dx+\int_{\Omega_{0}}\int^{\tau}_{0^{+}}\frac{\partial \zeta}{\partial t}\phi dtdx\\
=&\int_{\Omega_{0}}\zeta\phi|^{\tau}_{0^{+}}dx=0,
\end{aligned}
\end{equation}
and
\begin{equation}\label{4-5}
\begin{aligned}
\int_{\Omega_{0}}\int^{\tau}_{0^{+}}\frac{d_{1}}{\rho^{2}(t)}(\Delta \zeta\phi- \Delta\phi\zeta)dtdx=&\int^{\tau}_{0^{+}}\frac{d_{1}}{\rho^{2}(t)}\bigg[\int_{\Omega_{0}}\phi d\nabla \zeta-\int_{\Omega_{0}}\zeta d\nabla \phi\bigg]dt\\
=&\int^{\tau}_{0^{+}}\frac{d_{1}}{\rho^{2}(t)}\bigg[\Big(\phi \nabla \zeta-\zeta\nabla \phi\Big)\Big|_{\partial \Omega_{0}}\bigg]dt=0.
\end{aligned}
\end{equation}
Using \eqref{4-4} and \eqref{4-5}, \eqref{4-3} can be simplified to
\begin{equation}\label{4-6}
\begin{aligned}
(\lambda_{1}-\mu_{1})\int_{\Omega_{0}}\int^{\tau}_{0^{+}}\phi \zeta dtdx=&\int_{\Omega_{0}}\int^{\tau}_{0^{+}}(f'(0)\eta\phi-a_{12}\psi\zeta)dtdx.
\end{aligned}
\end{equation}
A similar discussion yields that
\begin{equation}\label{4-7}
\begin{aligned}
(\lambda_{1}-\mu_{1})\int_{\Omega_{0}}\int^{\tau}_{0^{+}}\psi \eta dtdx=&\int_{\Omega_{0}}\int^{\tau}_{0^{+}}(a_{12}\psi\zeta-f'(0)\eta\phi)dtdx.
\end{aligned}
\end{equation}
Adding \eqref{4-6} and \eqref{4-7} gives that
\begin{equation*}
\begin{aligned}
(\lambda_{1}-\mu_{1})\int_{\Omega_{0}}\int^{\tau}_{0^{+}}\psi \eta dtdx=(\mu_{1}-\lambda_{1})\int_{\Omega_{0}}\int^{\tau}_{0^{+}}\phi \zeta dtdx.
\end{aligned}
\end{equation*}
This implies that $\lambda_{1}=\mu_{1}$. This completes the proof.
\end{proof}
\end{theorem}
\autoref{theorem 4-1} tells us that the principal eigenvalue $\lambda_{1}$ of problem \eqref{2-2} and the principal eigenvalue $\mu_{1}$ of corresponding adjoint problem \eqref{4-1} are equal. Based on this conclusion, we now discuss the effects of $\Omega_{0}$ and $g'(0)$ on the principal eigenvalue $\lambda_{1}$ of problem \eqref{2-2}.
\begin{theorem}\label{theorem 4-2}
The following assertions are valid.
\begin{enumerate}
\item[$(1)$]
$\lambda_{1}$ is strongly decreasing about $g'(0)$.
\item[$(2)$]
$\lambda_{1}$ is strongly decreasing about $\Omega_{0}$, i.e., $\lambda_{1}(\Omega^{2}_{0})<\lambda_{1}(\Omega^{1}_{0})$ holds in that meaning $\Omega^{1}_{0}\subset \Omega^{2}_{0}$ and  $\Omega^{2}_{0}-\Omega^{1}_{0}$ is a nonempty and open set.
\end{enumerate}
\begin{proof}
(1) Obviously, $\phi$, $\psi$, and $\lambda_{1}$ are all $\mathbb{C}^{1}$-functions of $g'(0)$. For the sake of simplicity, denote $\frac{\partial \phi}{\partial g'(0) }$, $\frac{\partial \psi}{\partial g'(0) }$, and $\frac{\partial \lambda_{1}}{\partial g'(0) }$ by $\phi'$, $\psi'$, and $\lambda_{1}'$ respectively. Differentiating problem \eqref{2-2} with respect to $g'(0)$, one can obtain that
\begin{eqnarray}\label{4-8}
\left\{
\begin{array}{ll}
\frac{\partial \phi'}{\partial t}=\frac{d_{1}}{\rho^{2}(t)}\Delta\phi'-\big(\frac{n\dot{\rho}(t)}{\rho(t)}+a_{11}-\lambda_{1}\big)\phi'+a_{12}\psi'+\lambda_{1}' \phi,\; &\, t\in(0^{+}, \tau], x\in\Omega_{0}, \\[2mm]
\frac{\partial \psi'}{\partial t}=\frac{d_{2}}{\rho^{2}(t)}\Delta\psi' -\big(\frac{n\dot{\rho}(t)}{\rho(t)}+a_{22}-\lambda_{1}\big)\psi'+f'(0)\phi'+\lambda_{1}' \psi,\; &\, t\in(0^{+}, \tau], x\in\Omega_{0}, \\[2mm]
\phi'=0, \psi'=0,\; &\, t\in[0, \tau], x\in\partial\Omega_{0},\\[2mm]
\phi'(0,x)=\phi'(\tau,x), \psi'(0,x)=\psi'(\tau,x),\; &\,x\in\xoverline{\Omega}_{0},\\[2mm]
\phi'(0^{+},x)=\phi(0,x)+g'(0)\phi'(0,x), \psi'(0^{+},x)=\psi'(0,x), \; &\, x\in\Omega_{0}.
\end{array} \right.
\end{eqnarray}
By multiplying the first equations of problem \eqref{4-1} and problem \eqref{4-8} by $\phi'$ and $\zeta$ respectively,  it follows that
\begin{eqnarray}\label{4-9}
\left\{
\begin{array}{ll}
\frac{\partial \zeta}{\partial t}\phi'=-\frac{d_{1}}{\rho^{2}(t)}\Delta \zeta\phi'+\frac{n\dot{\rho}(t)}{\rho(t)}\zeta\phi'+a_{11}\zeta\phi'-f'(0)\eta\phi'-\mu_{1} \zeta\phi' ,\\[2mm]
\frac{\partial \phi'}{\partial t}\zeta=\frac{d_{1}}{\rho^{2}(t)}\Delta\phi'\zeta-\big(\frac{n\dot{\rho}(t)}{\rho(t)}+a_{11}-\lambda_{1}\big)\phi'\zeta+a_{12}\psi'\zeta+\lambda_{1}' \phi\zeta.
\end{array}
\right.
\end{eqnarray}
Then, integrating on each side of the first and second equations of \eqref{4-9} on $(0^{+}, \tau]\times\Omega_{0}$ respectively, and adding the obtained results give that
\begin{equation}\label{4-10}
\begin{aligned}
\int_{\Omega_{0}}\int^{\tau}_{0^{+}}\bigg(\frac{\partial \zeta}{\partial t}\phi'+\frac{\partial \phi'}{\partial t}\zeta\bigg)dtdx
=&\int_{\Omega_{0}}\int^{\tau}_{0^{+}}\big(a_{12}\psi'\zeta-f'(0)\eta\phi'\big)dtdx+\lambda_{1}'\int_{\Omega_{0}}\int^{\tau}_{0^{+}}\phi\zeta dtdx\\
&+\int_{\Omega_{0}}\int^{\tau}_{0^{+}}(\lambda_{1}-\mu_{1})\phi'\zeta dtdx.
\end{aligned}
\end{equation}
With the help of \autoref{theorem 4-1}, Fubini's theorem, and the formula of integration by parts, \eqref{4-10} can be further simplified to
\begin{equation}\label{4-11}
\begin{aligned}
\int_{\Omega_{0}}\int^{\tau}_{0^{+}}\big(f'(0)\eta\phi'-a_{12}\psi'\zeta\big)dtdx=\frac{1}{g'(0)}\int_{\Omega_{0}}\phi(0,x)\zeta (0,x)dx+\lambda_{1}'\int_{\Omega_{0}}\int^{\tau}_{0^{+}}\phi\zeta dtdx.
\end{aligned}
\end{equation}
A similar discussion yields that
\begin{equation}\label{4-12}
\begin{aligned}
\int_{\Omega_{0}}\int^{\tau}_{0^{+}}(a_{12}\psi'\zeta-f'(0)\eta\phi')dtdx=\lambda_{1}'\int_{\Omega_{0}}\int^{\tau}_{0^{+}}\psi\eta dtdx.
\end{aligned}
\end{equation}
Adding \eqref{4-11} and \eqref{4-12} gives that
\begin{equation*}
\begin{aligned}
\lambda_{1}'=-\frac{\int_{\Omega_{0}}\phi(0,x)\zeta (0,x)dx}{g'(0)\int_{\Omega_{0}}\int^{\tau}_{0^{+}}\big(\phi\zeta+\psi\eta\big)dtdx}<0.
\end{aligned}
\end{equation*}
This means that $\lambda_{1}$ is strongly decreasing about $g'(0)$. This completes the proof of assertion (1).

(2) Let $\big(\lambda_{1}(\Omega^{1}_{0}), \phi_{1}, \psi_{1}\big)$ and $\big(\lambda_{1}(\Omega^{2}_{0}), \phi_{2}, \psi_{2}\big)$ be the principal eigenpairs of eigenvalue problem \eqref{2-2} with initial regions $\Omega^{1}_{0}$ and $\Omega^{2}_{0}$, respectively, and $\big(\mu(\Omega^{1}_{0}), \zeta_{1}, \eta_{1}\big)$ and $\big(\mu(\Omega^{2}_{0}), \zeta_{2}, \eta_{2}\big)$ be the principal eigenpairs of adjoint problem \eqref{4-1} with initial regions $\Omega^{1}_{0}$ and $\Omega^{2}_{0}$, respectively. By multiplying the first equations of problem \eqref{2-2} with initial region $\Omega^{2}_{0}$ and problem \eqref{4-1} with initial region $\Omega^{1}_{0}$ by $\zeta_{1}$ and $\phi_{2}$ respectively,  it follows that
\begin{eqnarray}\label{4-13}
\left\{
\begin{array}{ll}
\frac{\partial \phi_{2}}{\partial t}\zeta_{1}=\frac{d_{1}}{\rho^{2}(t)}\Delta \phi_{2}\zeta_{1}-\frac{n\dot{\rho}(t)}{\rho(t)}\phi_{2}\zeta_{1}-a_{11}\phi_{2}\zeta_{1}+a_{12}\psi_{2}\zeta_{1}+\lambda_{1}\big(\Omega^{2}_{0}\big) \phi_{2} \zeta_{1}, \\[2mm]
\frac{\partial \zeta_{1}}{\partial t}\phi_{2}=-\frac{d_{1}}{\rho^{2}(t)}\Delta \zeta_{1}\phi_{2}+\frac{n\dot{\rho}(t)}{\rho(t)}\zeta_{1}\phi_{2}+a_{11}\zeta_{1}\phi_{2}-f'(0)\eta_{1}\phi_{2}-\mu\big(\Omega^{1}_{0}\big) \zeta_{1}\phi_{2}.
\end{array}
\right.
\end{eqnarray}
Then, integrating both sides of the first and second equations of \eqref{4-13} on $(0^{+}, \tau]\times\Omega^{1}_{0}$ respectively, and adding the obtained results give that
\begin{equation}\label{4-14}
\begin{aligned}
\int_{\Omega^{1}_{0}}\int^{\tau}_{0^{+}}\bigg(\frac{\partial \zeta_{1}}{\partial t}\phi_{2}+\frac{\partial \phi_{2}}{\partial t}\zeta_{1}\bigg)dtdx
=&\int_{\Omega^{1}_{0}}\int^{\tau}_{0^{+}}\Bigg[\big(a_{12}\psi_{2}\zeta_{1}-f'(0)\eta_{1}\phi_{2}\big)+\frac{d_{1}}{\rho^{2}(t)}(\Delta\phi_{2}\zeta_{1}-\Delta\zeta_{1}\phi_{2})\Bigg]dtdx\\
&+\Big(\lambda_{1}\big(\Omega^{2}_{0}\big) -\mu\big(\Omega^{1}_{0}\big)\Big )\int_{\Omega^{1}_{0}}\int^{\tau}_{0^{+}}\phi_{2} \zeta_{1} dtdx.
\end{aligned}
\end{equation}
With the help of \autoref{theorem 4-1}, Fubini's theorem, and the formula of integration by parts, \eqref{4-14} can be further simplified to
\begin{equation}\label{4-15}
\begin{aligned}
\big(\lambda_{1}\big(\Omega^{2}_{0}\big) -\lambda_{1}\big(\Omega^{1}_{0}\big)\big)\int_{\Omega^{1}_{0}}\int^{\tau}_{0^{+}}\phi_{2} \zeta_{1} dtdx
=&\int_{\Omega^{1}_{0}}\int^{\tau}_{0^{+}}\big(f'(0)\eta_{1}\phi_{2}-a_{12}\psi_{2}\zeta_{1}\big)dtdx+\int^{\tau}_{0^{+}}\int_{\partial\Omega^{1}_{0}}\frac{d_{1}\phi_{2}}{\rho^{2}(t)}\frac{\partial \zeta_{1} }{\partial n}dxdt.\\
\end{aligned}
\end{equation}
A similar discussion yields that
\begin{equation}\label{4-16}
\begin{aligned}
\big(\lambda_{1}\big(\Omega^{2}_{0}\big) -\lambda_{1}\big(\Omega^{1}_{0}\big)\big)\int_{\Omega^{1}_{0}}\int^{\tau}_{0^{+}}\psi_{2} \eta_{1} dtdx
=&\int_{\Omega^{1}_{0}}\int^{\tau}_{0^{+}}\big(a_{12}\zeta_{1}\psi_{2}-f'(0)\phi_{2}\eta_{1})dtdx+\int^{\tau}_{0^{+}}\int_{\partial\Omega^{1}_{0}}\frac{d_{2}\psi_{2}}{\rho^{2}(t)}\frac{\partial \eta_{1} }{\partial n}dxdt.\\
\end{aligned}
\end{equation}
Adding \eqref{4-15} and \eqref{4-16} gives that
\begin{equation}\label{4-17}
\begin{aligned}
\lambda_{1}\big(\Omega^{2}_{0}\big)-\lambda_{1}\big(\Omega^{1}_{0}\big)=\frac{\int^{\tau}_{0^{+}}\int_{\partial\Omega^{1}_{0}}\Big[\frac{d_{1}\phi_{2}}{\rho^{2}(t)}\frac{\partial \zeta_{1} }{\partial n}+\frac{d_{2}\psi_{2}}{\rho^{2}(t)}\frac{\partial \eta_{1} }{\partial n}\Big]dxdt}{\int^{\tau}_{0^{+}}\int_{\partial\Omega^{1}_{0}}\big(\phi_{2} \zeta_{1}+\psi_{2} \eta_{1}\big)dtdx}.
\end{aligned}
\end{equation}
With the help of the strong maximum principle, it follows that $\frac{\partial \zeta_{1} }{\partial n}, \frac{\partial \eta_{1} }{\partial n}<0$ for all $(t,x)\in[0,\tau]\times\partial\Omega^{1}_{0}$. Moreover, the condition $\Omega^{1}_{0}\subset \Omega^{2}_{0}$ provides that $\phi_{2}\geq0(\not \equiv0)$ and $\psi_{2}\geq0(\not \equiv0)$ for $t\in[0,\tau]$ and $x\in\partial\Omega^{1}_{0}$. Therefore, \eqref{4-17} yields that
\begin{equation*}
\lambda_{1}\big(\Omega^{2}_{0}\big)-\lambda_{1}\big(\Omega^{1}_{0}\big)<0.
\end{equation*}
 This means that $\lambda_{1}$ is strictly monotonic decreasing with respect to $\Omega_{0}$. This completes the proof of assertion (2).
\end{proof}
\end{theorem}
\begin{Remark}
From the perspective of epidemic prevention and control, \autoref{theorem 4-2} provides the following statements.
\begin{enumerate}
\item[$(1)$]
The stronger the intensity $1-g'(0)$ of periodic disinfection of bacteria in the environment, the more beneficial to the fight against the epidemics.
\item[$(2)$]
Containing the region in which an epidemic is occurring within a smaller area has a positive influence on containing the disease.
\end{enumerate}
\end{Remark}

\autoref{theorem 4-2} above gives the monotonicity of $\lambda_{1}$ about $g'(0)$ and $\Omega_{0}$ . In order to facilitate further understanding of this monotonicity using numerical simulations, we next give some estimates of the principal eigenvalue $\lambda_{1}$ in some special cases.

\begin{theorem}\label{theorem 4-3}
The following assertions are valid.
\begin{enumerate}
\item[$(1)$]
Assume $g'(0)=1$. Then, we have that
\begin{equation*}
\lambda_{1}\leq\frac{(d_{1}+d_{2})\lambda_{0}\xoverline{\rho^{-2}}+a_{11}+a_{22}-\sqrt{\big[(d_{1}-d_{2})\lambda_{0}\xoverline{\rho^{-2}}+a_{11}-a_{22}\big]^{2}+4a_{12}f'(0)}}{2},
\end{equation*}
where $\xoverline{\rho^{-2}}=\frac{1}{\tau}\int_{0}^{\tau}\frac{1}{\rho^{2}(t)}dt$.
\item[$(2)$]
Assume $g'(0)=1$ and $d_{1}=d_{2}$. Then, we have that
\begin{equation*}
\lambda_{1}\geq\frac{d_{1}\lambda_{0}}{\tau}\int_{0}^{\tau}\frac{1}{\rho^{2}(t)}dt-\max\big(a_{12},f'(0)\big).
\end{equation*}
\item[$(3)$]
Assume $g'(0)=1$ and $\rho(t)=1$. Then, we have that
\begin{equation*}
\lambda_{1}=\frac{(d_{1}+d_{2})\lambda_{0}+a_{11}+a_{22}-\sqrt{(a_{22}+(d_{2}-d_{1})\lambda_{0}-a_{11})^{2}+4a_{12}f'(0)}}{2}.
\end{equation*}
\end{enumerate}
\begin{proof}
(1) Let $\phi(t,x)=\Phi(t)X(x)~\text{and}~\psi(t,x)=\Psi(t)X(x)=\kappa(t)\Phi(t)X(x)$, where $\big(\lambda_{0}, X(x)\big)$ is the principal eigenpair of the eigenvalue problem
\begin{eqnarray*}
\left\{
\begin{array}{ll}
-\Delta X(x)=\lambda X(x),\; &\, x\in \Omega_{0}, \\[2mm]
X(x)=0,\; &\, x\in \partial\Omega_{0}.
\end{array}
\right.
\end{eqnarray*}
With the help of assumption $g'(0)=1$, it follows from problem \eqref{2-2} that
\begin{eqnarray}\label{4-22}
\left\{
\begin{array}{ll}
\Phi'(t)=-\Big[\frac{d_{1}\lambda_{0}}{\rho^{2}(t)}+\frac{n\dot{\rho}(t)}{\rho(t)}+a_{11}-a_{12}\kappa(t)-\lambda_{1}\Big]\Phi(t),\; &\, t\in(0, \tau], \\[2mm]
\Psi'(t)=-\Big[\frac{d_{2}\lambda_{0}}{\rho^{2}(t)}+\frac{n\dot{\rho}(t)}{\rho(t)}+a_{22}-\frac{f'(0)}{\kappa(t)}-\lambda_{1}\Big]\Psi(t),\; &\, t\in(0, \tau], \\[2mm]
\Phi(0)=\Phi(\tau), \Psi(0)=\Psi(\tau).\; &\,
\end{array} \right.
\end{eqnarray}
Integrating on each side of the first and second equations of \eqref{4-22} on $(0^{+}, \tau]$ respectively, it follows that
\begin{equation}\label{4-23}
\begin{cases}
\begin{aligned}
\int_{0}^{\tau}\kappa(t)dt=&\frac{d_{1}\lambda_{0}}{a_{12}}\int_{0}^{\tau}\frac{1}{\rho^{2}(t)}dt+\frac{\tau}{a_{12}}(a_{11}-\lambda_{1}),\\
\int_{0}^{\tau}\frac{1}{\kappa(t)}dt=&\frac{d_{2}\lambda_{0}}{f'(0)}\int_{0}^{\tau}\frac{1}{\rho^{2}(t)}dt+\frac{\tau}{f'(0)}(a_{22}-\lambda_{1}).
\end{aligned}
\end{cases}
\end{equation}
With the help of \eqref{4-23} and the Cauchy-Schwarz inequality, the direct calculation yields that
\begin{equation}\label{4-24}
\Big[d_{1}\lambda_{0}\xoverline{\rho^{-2}}+a_{11}-\lambda_{1}\Big]\Big[d_{2}\lambda_{0}\xoverline{\rho^{-2}}+a_{22}-\lambda_{1}\Big]\geq a_{12}f'(0),
\end{equation}
where $\xoverline{\rho^{-2}}=\frac{1}{\tau}\int_{0}^{\tau}\frac{1}{\rho^{2}(t)}dt$. By making a simple calculation on \eqref{4-24}, $\lambda_{1}$ satisfies that
\begin{equation*}
\lambda_{1}\leq\frac{(d_{1}+d_{2})\lambda_{0}\xoverline{\rho^{-2}}+a_{11}+a_{22}-\sqrt{\big[(d_{1}+d_{2})\lambda_{0}\xoverline{\rho^{-2}}+a_{11}+a_{22}\big]^{2}+4a_{12}f'(0)}}{2}.
\end{equation*}
This completes the proof of assertion (1).

(2) Adding the first two equations of \eqref{4-22} gives
\begin{equation}\label{4-25}
\Big(\Phi(t)+\Psi(t)\Big)'\leq-\Bigg[\frac{d_{1}\lambda_{0}}{\rho^{2}(t)}+\frac{n\dot{\rho}(t)}{\rho(t)}-\max\big(a_{12},f'(0)\big)-\lambda_{1}\Bigg]\Big(\Phi(t)+\Psi(t)\Big)
\end{equation}
Consider the following Cauchy problem
\begin{eqnarray}\label{4-26}
\left\{
\begin{array}{ll}
h'(t)=-\bigg[\frac{d_{1}\lambda_{0}}{\rho^{2}(t)}+\frac{n\dot{\rho}(t)}{\rho(t)}-\max\big(a_{12},f'(0)\big)-\lambda\bigg]h(t),\; &\, t\in(0, \tau], \\[2mm]
h(0)=h(\tau).\; &\,
\end{array} \right.
\end{eqnarray}
Integrating both sides of the first equation of \eqref{4-26} over $[0, \tau]$ yields
\begin{equation*}
\lambda=\frac{d_{1}\lambda_{0}}{\tau}\int_{0}^{\tau}\frac{1}{\rho^{2}(t)}dt-\max\big(a_{12},f'(0)\big).
\end{equation*}
Therefore, it follows from \eqref{4-25} that
\begin{equation*}
\lambda_{1}\geq\frac{d_{1}\lambda_{0}}{\tau}\int_{0}^{\tau}\frac{1}{\rho^{2}(t)}dt-\max\big(a_{12},f'(0)\big).
\end{equation*}
This completes the proof of assertion (2).

(3) By using assumption $\rho(t)=1$, the first two equations in \eqref{4-22} can be written as
\begin{equation*}
\Bigg(
\begin{matrix}
   \Phi'(t)\\
   \Psi'(t) \\
  \end{matrix}
\Bigg)=\Bigg(
\begin{matrix}
   \lambda_{1}-d_{1}\lambda_{0}-a_{11}&a_{12}\\
   f'(0)&   \lambda_{1}-d_{2}\lambda_{0}-a_{22}\\
  \end{matrix}
\Bigg)
\Bigg(
\begin{matrix}
   \Phi(t) \\
   \Psi(t) \\
  \end{matrix}
\Bigg):=A
\Bigg(
\begin{matrix}
   \Phi(t) \\
   \Psi(t) \\
  \end{matrix}
\Bigg).
\end{equation*}
Then, it follows from the characteristic equation $|A-\kappa E|=0$ that
\begin{equation*}
\kappa_{1,2}=\frac{2\lambda_{1}-(d_{1}+d_{2})\lambda_{0}-a_{11}-a_{22}\pm\sqrt{(a_{22}+(d_{2}-d_{1})\lambda_{0}-a_{11})^{2}+4a_{12}f'(0)}}{2}.
\end{equation*}
Without loss of generality, assume that $\kappa_{1}>\kappa_{2}$. The linearly independent eigenvectors $(x_{11}, x_{12})$ and $(x_{21}, x_{22})$ associated with eigenvalues $\kappa_{1}$ and $\kappa_{2}$ satisfy
\begin{equation*}
\big(
\begin{matrix}
   x_{i1}&x_{i2}\\
  \end{matrix}
\big)
\Bigg(
\begin{matrix}
   \lambda_{1}-d_{1}\lambda_{0}-a_{11}-\kappa_{i}&a_{12}\\
   f'(0)&   \lambda_{1}-d_{2}\lambda_{0}-a_{22}-\kappa_{i}\\
  \end{matrix}
\Bigg)=
\big(
\begin{matrix}
   0&0\\
  \end{matrix}
\big)
\end{equation*}
for $i=1,2$, and a simple calculation yields that
\begin{equation*}
(x_{11}, x_{12})=\big(a_{12}, a_{11}+d_{1}\lambda_{0}-\lambda_{1}+\kappa_{1}\big)
\end{equation*}
and
\begin{equation*}
(x_{21}, x_{22})=\big(a_{22}+d_{2}\lambda_{0}-\lambda_{1}+\kappa_{2}, f'(0)\big).
\end{equation*}

Next, we consider the algebraic equations
\begin{equation*}
\Bigg(
\begin{matrix}
   e^{\kappa_{1}t} \\
   ke^{\kappa_{2}t} \\
  \end{matrix}
\Bigg)=
\Bigg(
\begin{matrix}
   x_{11}&x_{12}\\
   x_{21}&x_{22}\\
  \end{matrix}
\Bigg)
\Bigg(
\begin{matrix}
   \Phi(t)\\
   \Psi(t) \\
  \end{matrix}
\Bigg):=B
\Bigg(
\begin{matrix}
   \Phi(t)\\
   \Psi(t) \\
  \end{matrix}
\Bigg),
\end{equation*}
and a simple calculation yields that
\begin{equation*}
(\Phi(t), \Psi(t))=\Bigg(\frac{f'(0)e^{\kappa_{1}t}-\big(a_{11}+d_{1}\lambda_{0}-\lambda_{1}+\kappa_{1}\big)ke^{\kappa_{2}t}}{|B|},
\frac{a_{12}ke^{\kappa_{2}t}-\big(a_{22}+d_{2}\lambda_{0}-\lambda_{1}+\kappa_{2}\big)e^{\kappa_{1}t}}{|B|}\Bigg),
\end{equation*}
where
\begin{equation*}
\begin{aligned}
|B|=&a_{12}f'(0)-\big(a_{11}+d_{1}\lambda_{0}-\lambda_{1}+\kappa_{1}\big)\big(a_{22}+d_{2}\lambda_{0}-\lambda_{1}+\kappa_{2}\big)\\
=&a_{12}f'(0)+\big(a_{22}+d_{2}\lambda_{0}-\lambda_{1}+\kappa_{2}\big)^{2}>0.
\end{aligned}
\end{equation*}
Using the third equation of \eqref{4-22}, it follows that
\begin{equation*}
\begin{cases}
\begin{aligned}
&f'(0)-\big(a_{11}+d_{1}\lambda_{0}-\lambda_{1}+\kappa_{1}\big)k=f'(0)e^{\kappa_{1}\tau}-\big(a_{11}+d_{1}\lambda_{0}-\lambda_{1}+\kappa_{1}\big)ke^{\kappa_{2}\tau},\\
&~~~a_{12}k-\big(a_{22}+d_{2}\lambda_{0}-\lambda_{1}+\kappa_{2}\big)=a_{12}ke^{\kappa_{2}\tau}-\big(a_{22}+d_{2}\lambda_{0}-\lambda_{1}+\kappa_{2}\big)e^{\kappa_{1}\tau},
\end{aligned}
\end{cases}
\end{equation*}
and its solution is
\begin{equation*}
\kappa_{1}=0, k=0.
\end{equation*}
Therefore, the principal eigenvalue $\lambda_{1}$ of problem \eqref{2-2} satisfies
\begin{equation*}
\lambda_{1}=\frac{(d_{1}+d_{2})\lambda_{0}+a_{11}+a_{22}-\sqrt{(a_{22}+(d_{2}-d_{1})\lambda_{0}-a_{11})^{2}+4a_{12}f'(0)}}{2},
\end{equation*}
and the positive eigenfunction corresponding to this eigenvalue is
\begin{equation*}
(\Phi(t), \Psi(t))=\Bigg(\frac{f'(0)}{|B|}, \frac{a_{11}+d_{1}\lambda_{0}-\lambda_{1}}{|B|}\Bigg).
\end{equation*}
This completes the proof of assertion (3).
\end{proof}
\end{theorem}
\section{Numerical simulations}\label{Section-6}
With the assist of numerical simulations, this section will validate the correctness of the theoretical results and help us to further understand the important role of periodic pulses and evolving region in the dynamical behaviors of the model.

In all simulations, we consider model \eqref{1-2-20} with Beverton-Holt function and choose the initial region as $\Omega_{0}=[0, \pi]$. Additionally, the initial values are as follows:
\begin{equation*}
u_{0}(x)=8\sin(x),~v_{0}(x)=0.1\sin(x), ~x\in[0, \pi].
\end{equation*}

Next, we give separately the influences of periodic impulses and evolving region on the dynamical behaviours of model \eqref{1-2-20}.
\subsection{The influence of periodic impulse}
This subsection presents numerical simulations on a fixed region in order to conveniently determine the sign of the principal eigenvalue.
Some of the parameters are chosen as $d_{1}=0.05$, $d_{2}=1$, $a_{11}=0.2$, $a_{22}=0.15$, and $a_{1}=10$. Other of the parameters are given in each of the following Examples to show the different dynamical behaviours. The forthcoming \textcolor[rgb]{0.00,0.00,1.00}{Examples} \ref{exm1} and \ref{exm2} will give the influence of periodic impulses on the extinction and persistence of the faecal-oral transmission diseases, respectively.
\begin{figure}[ht]
\centering
\subfigure{ {
\includegraphics[width=0.30\textwidth]{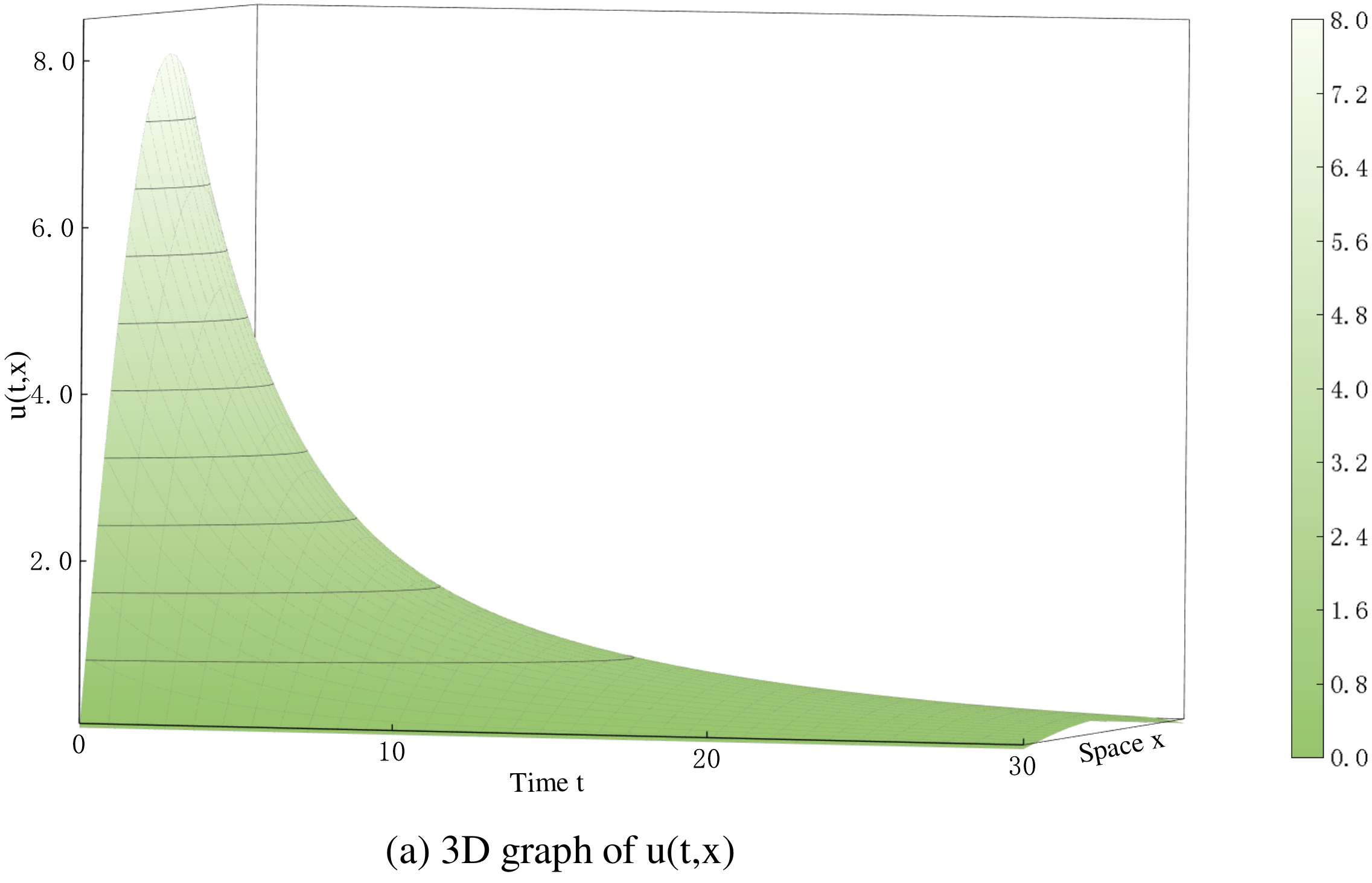}
} }
\subfigure{ {
\includegraphics[width=0.30\textwidth]{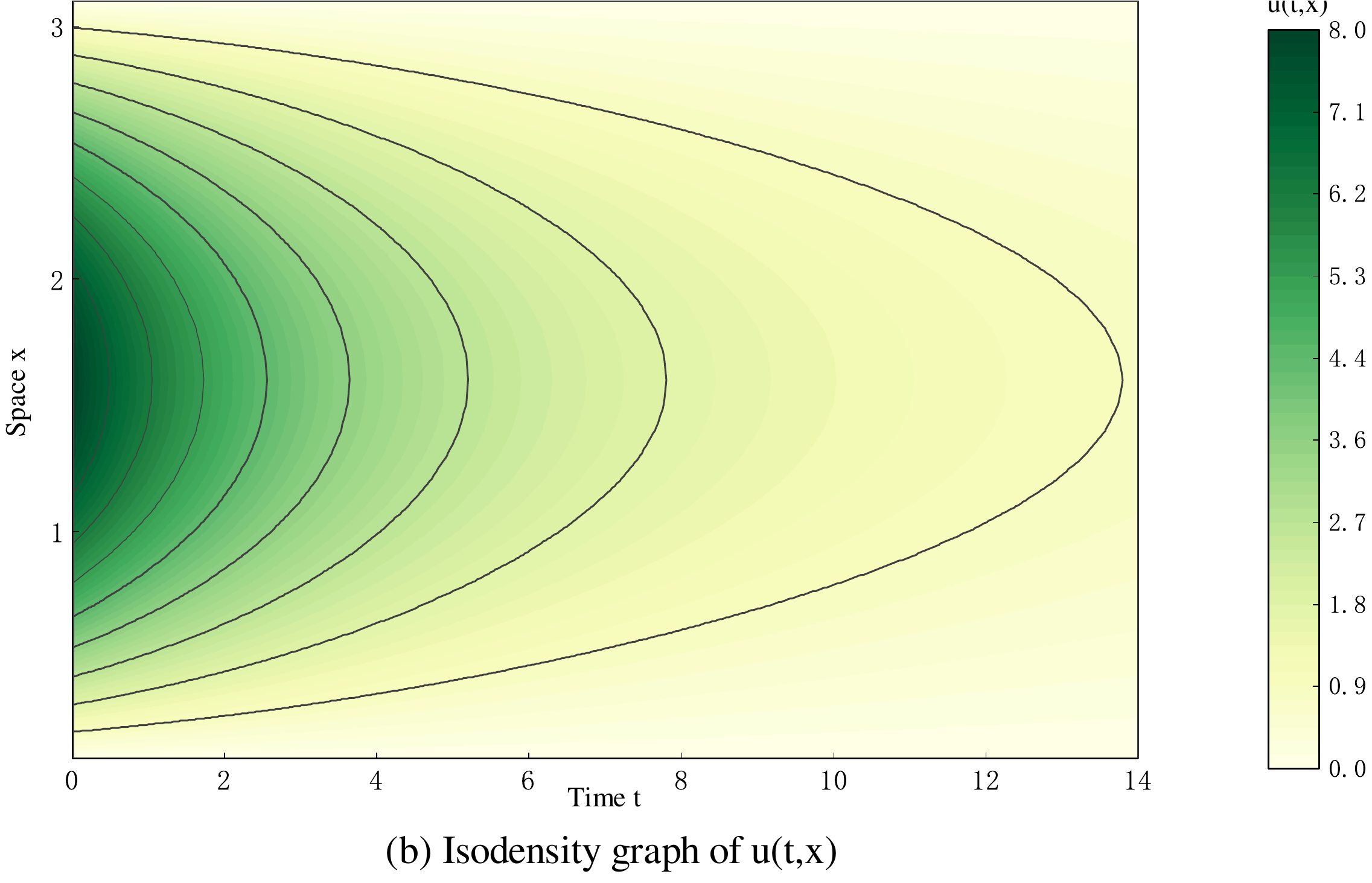}
} }
\subfigure{ {
\includegraphics[width=0.30\textwidth]{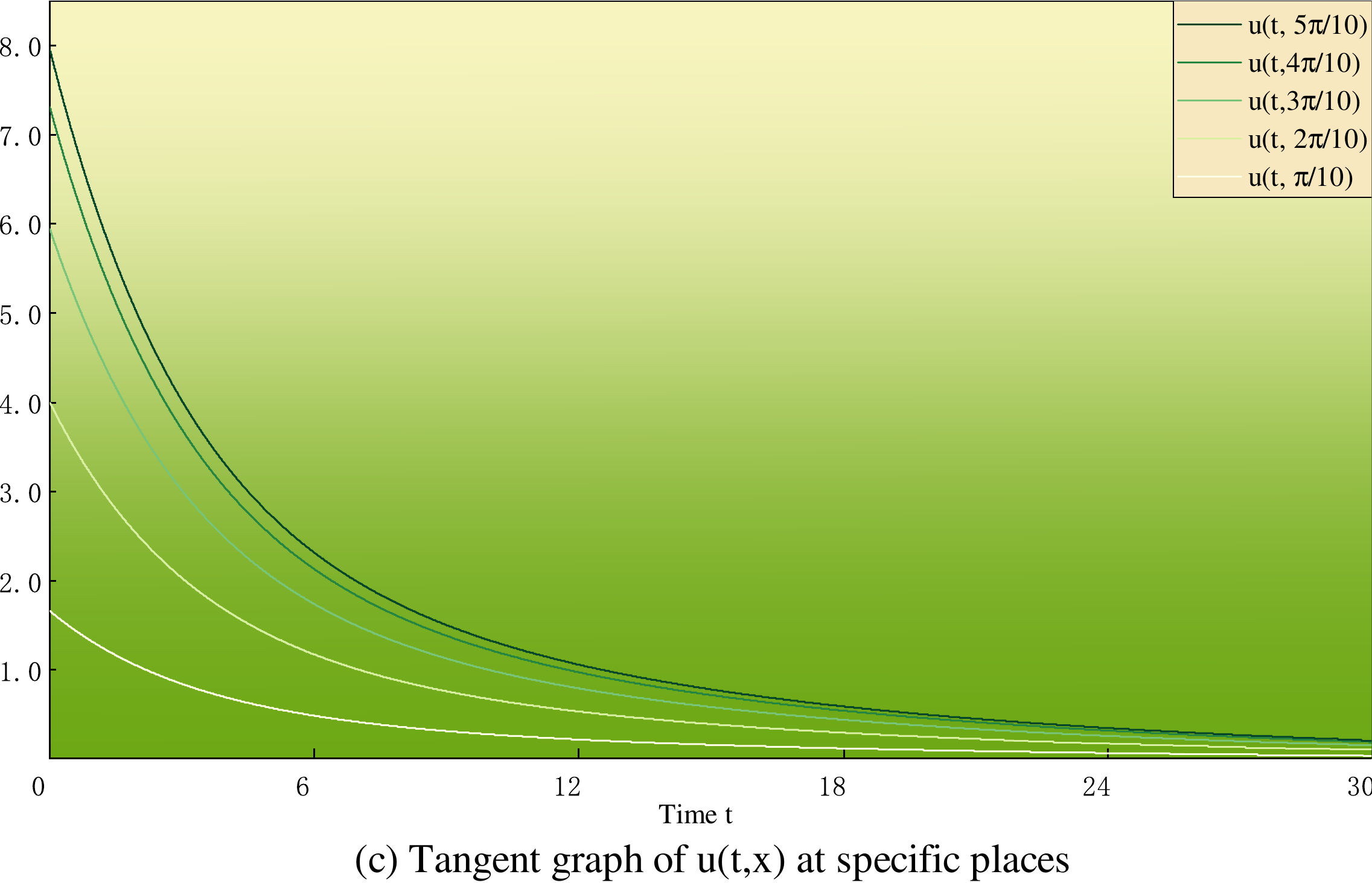}
} }
\subfigure{ {
\includegraphics[width=0.30\textwidth]{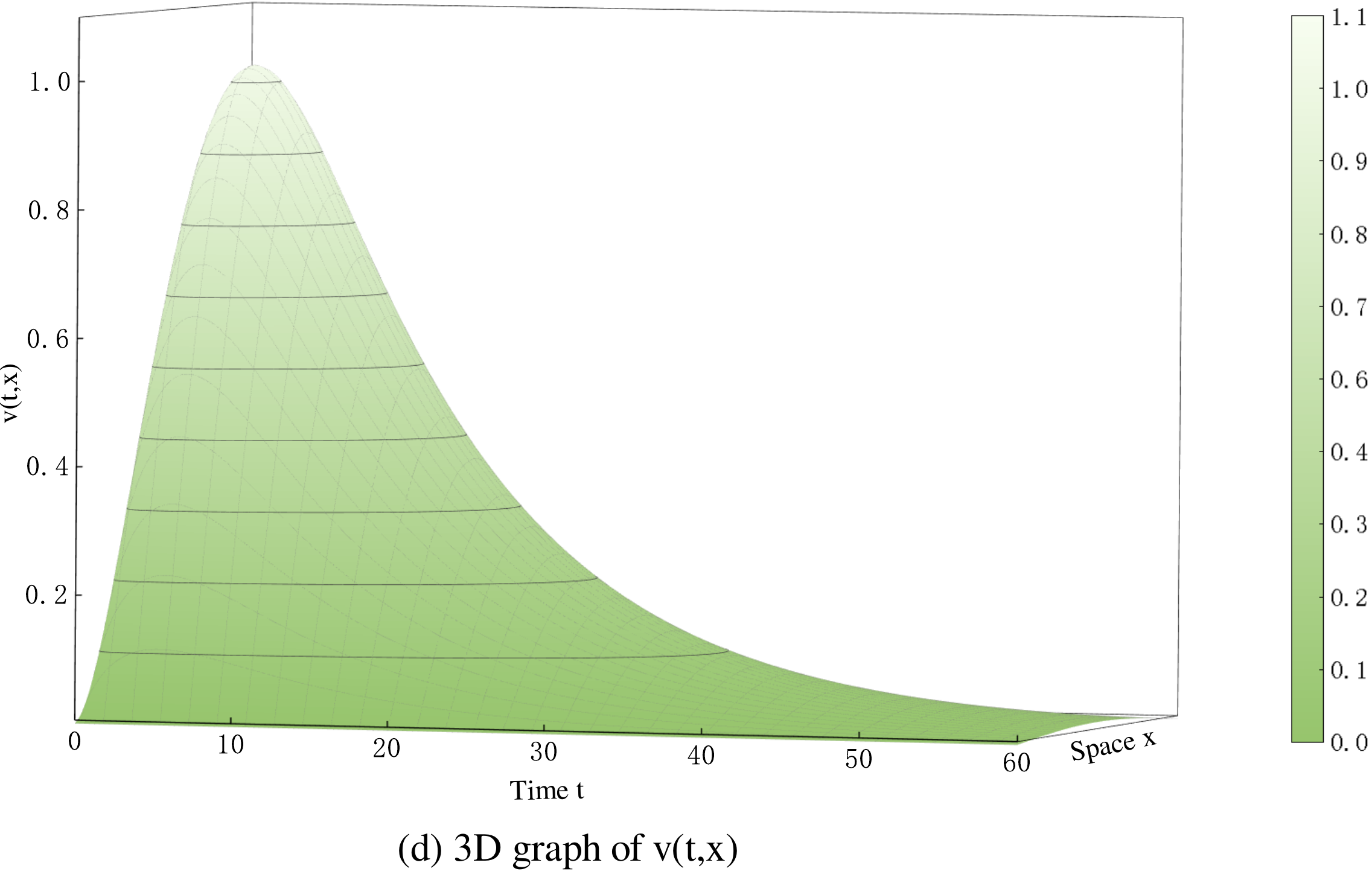}
} }
\subfigure{ {
\includegraphics[width=0.30\textwidth]{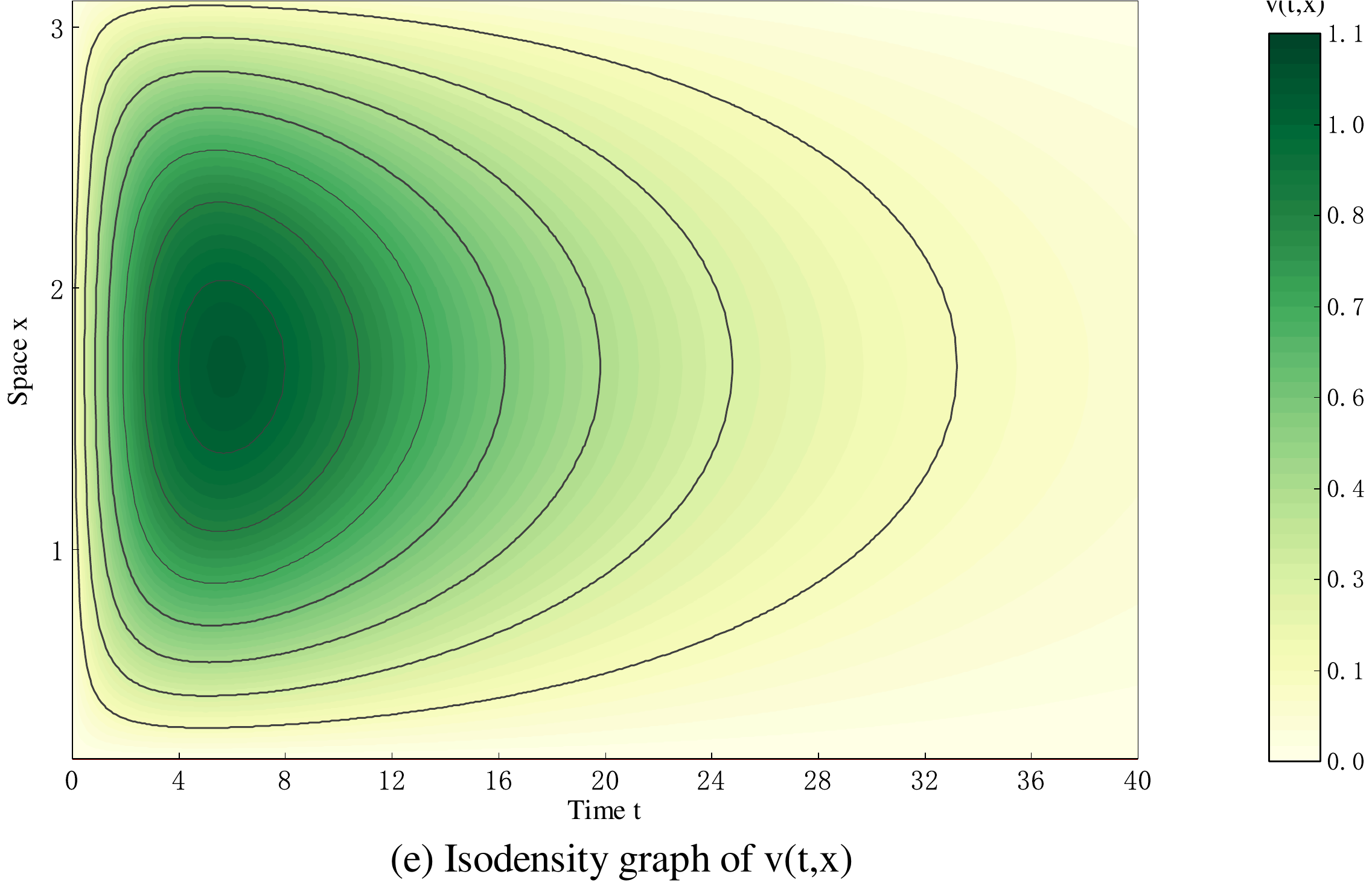}
} }
\subfigure{ {
\includegraphics[width=0.30\textwidth]{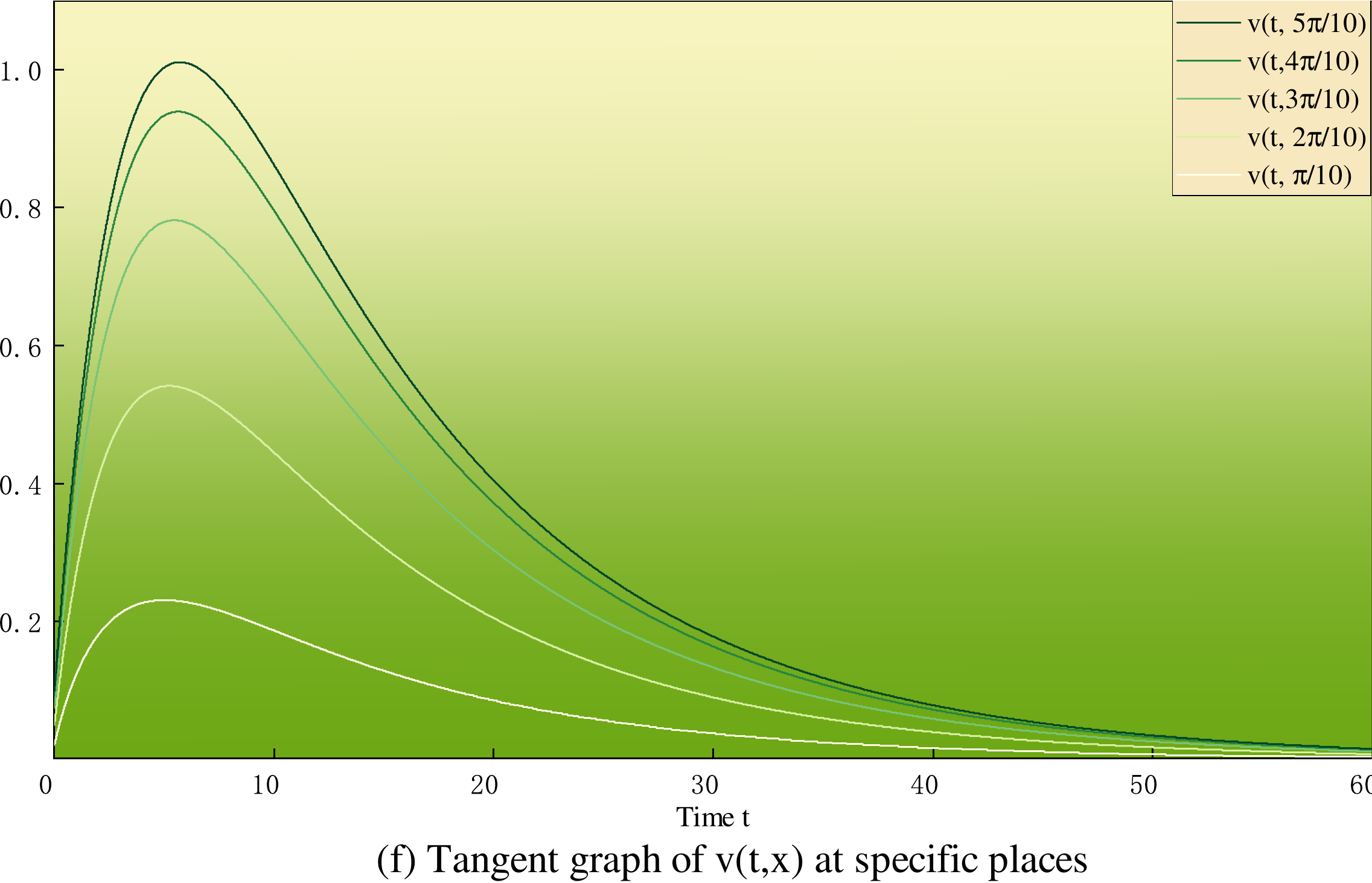}
} }
\caption{When $g(u)=u$(without impulse), graphs (a)-(f) show $(u, v)$ decays to $(0, 0)$.}
\label{A}
\end{figure}
\begin{figure}[!htt]
\centering
\subfigure{ {
\includegraphics[width=0.30\textwidth]{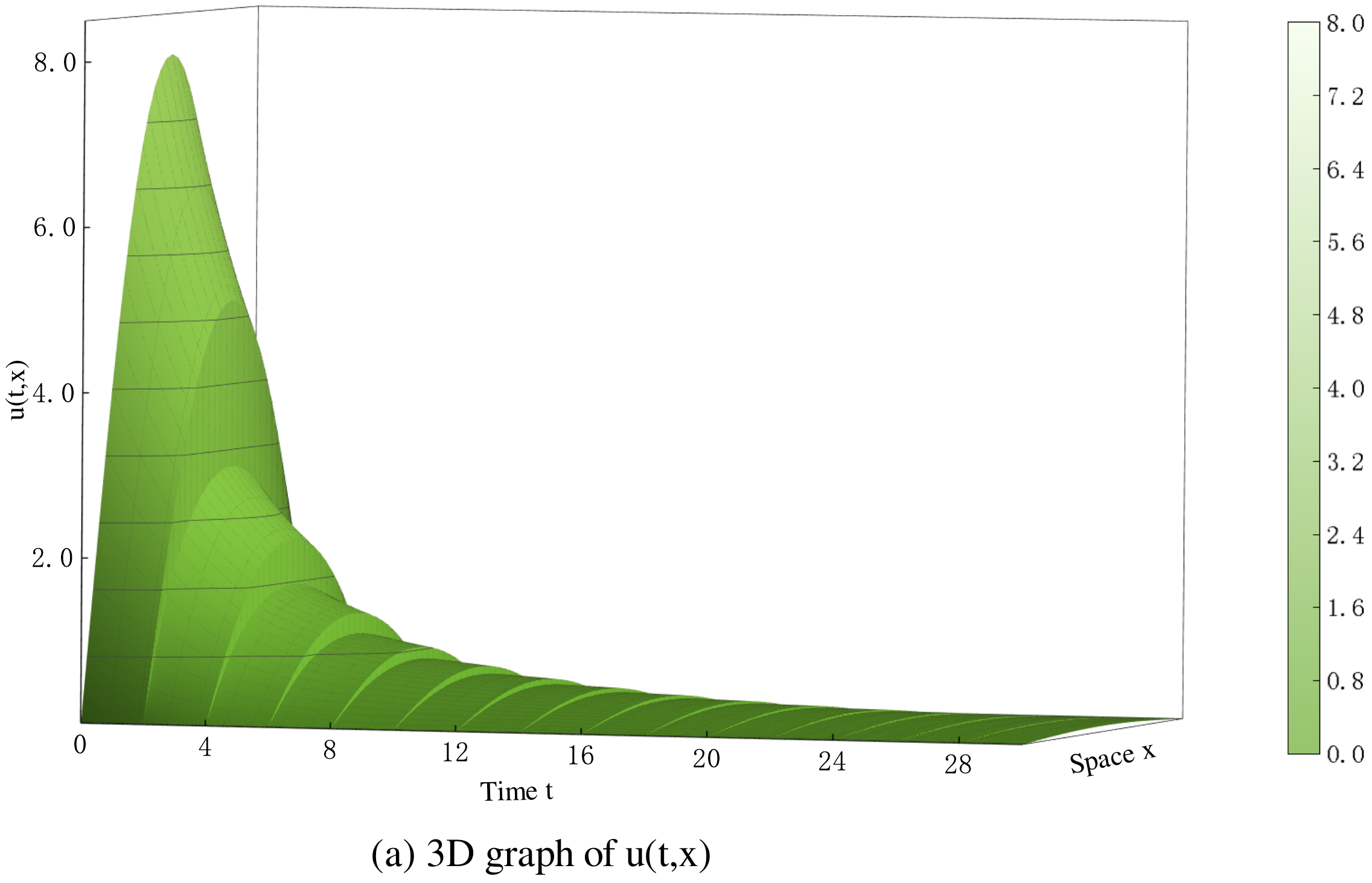}
} }
\subfigure{ {
\includegraphics[width=0.30\textwidth]{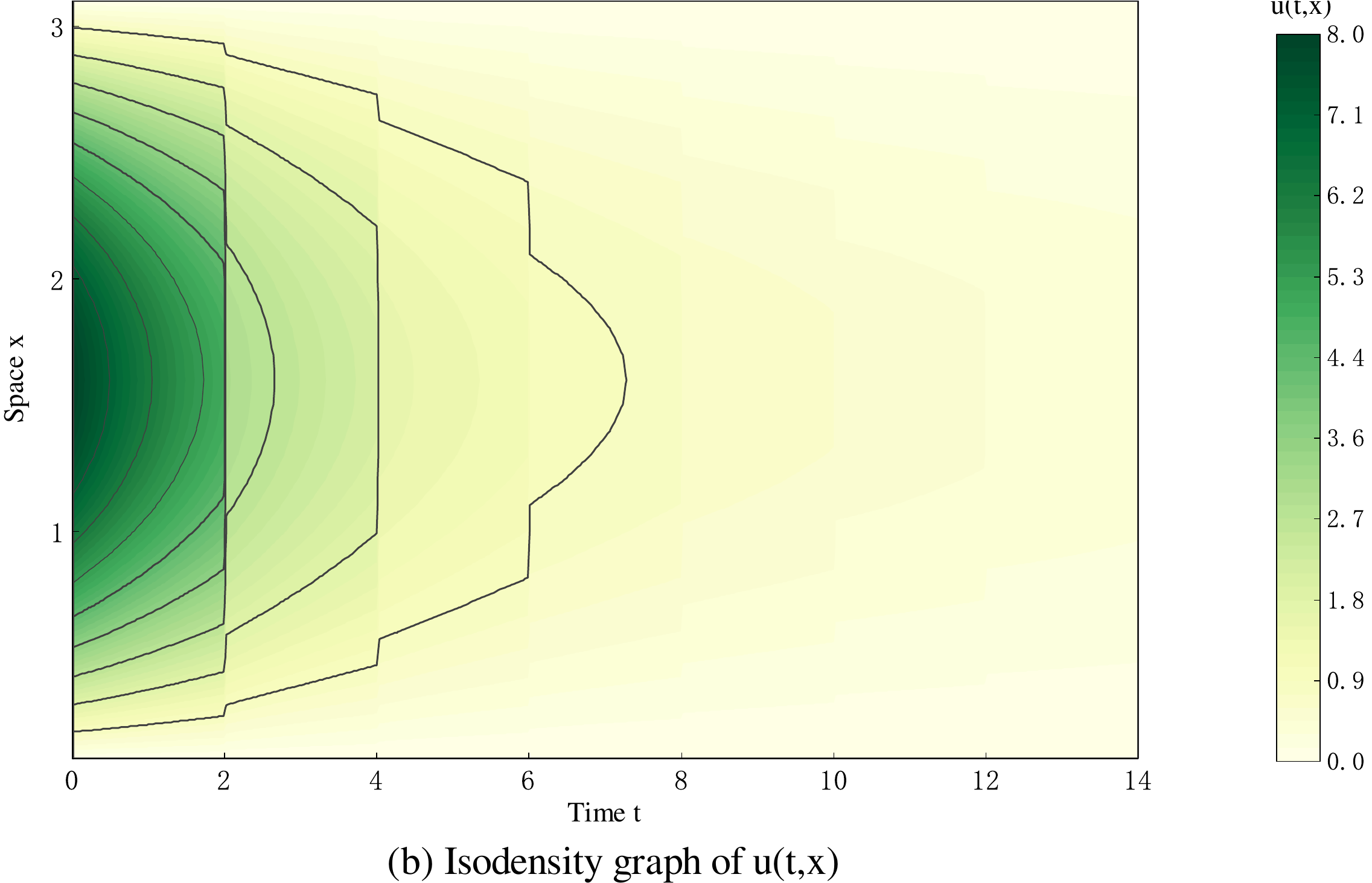}
} }
\subfigure{ {
\includegraphics[width=0.30\textwidth]{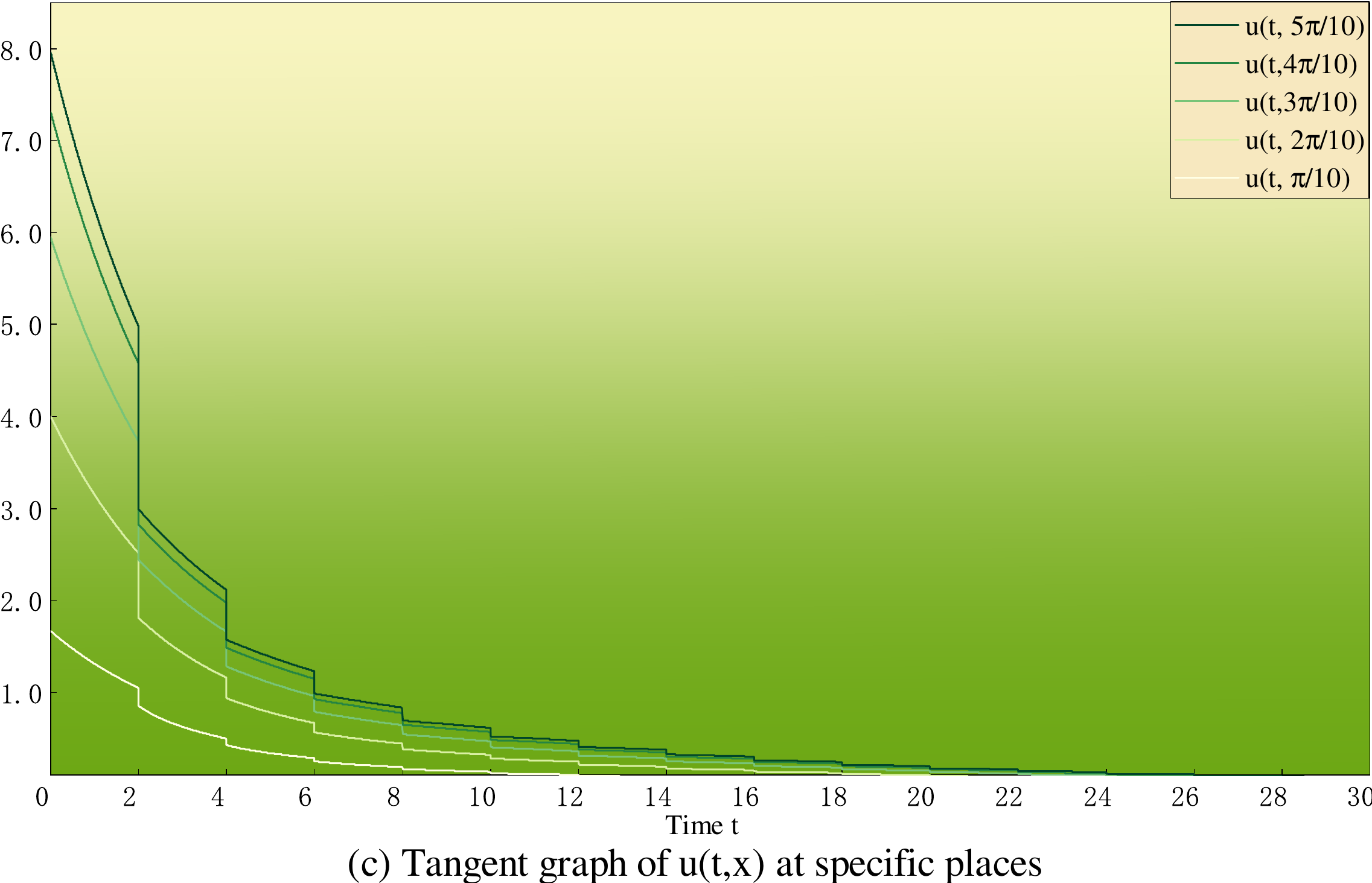}
} }
\subfigure{ {
\includegraphics[width=0.30\textwidth]{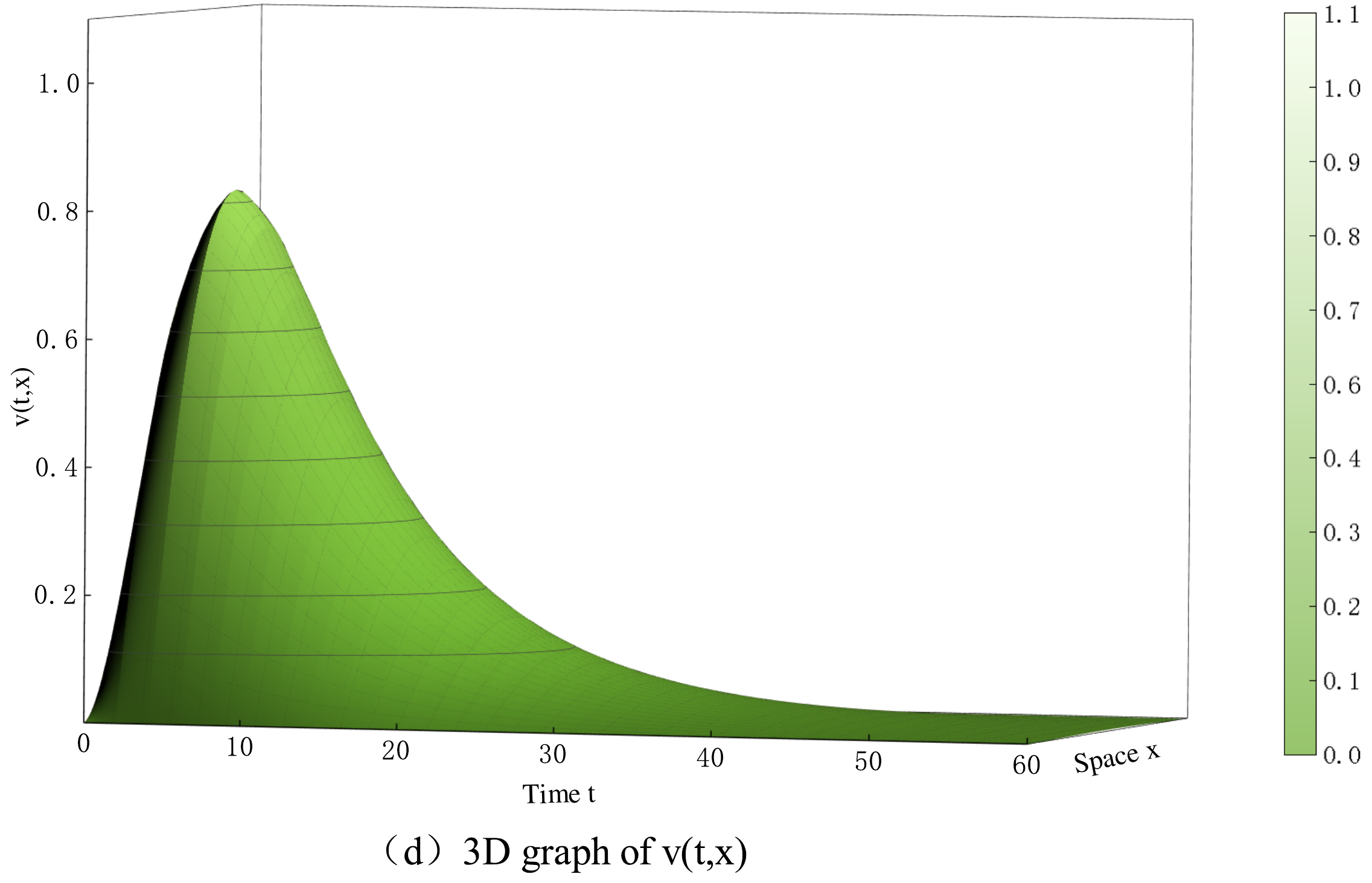}
} }
\subfigure{ {
\includegraphics[width=0.30\textwidth]{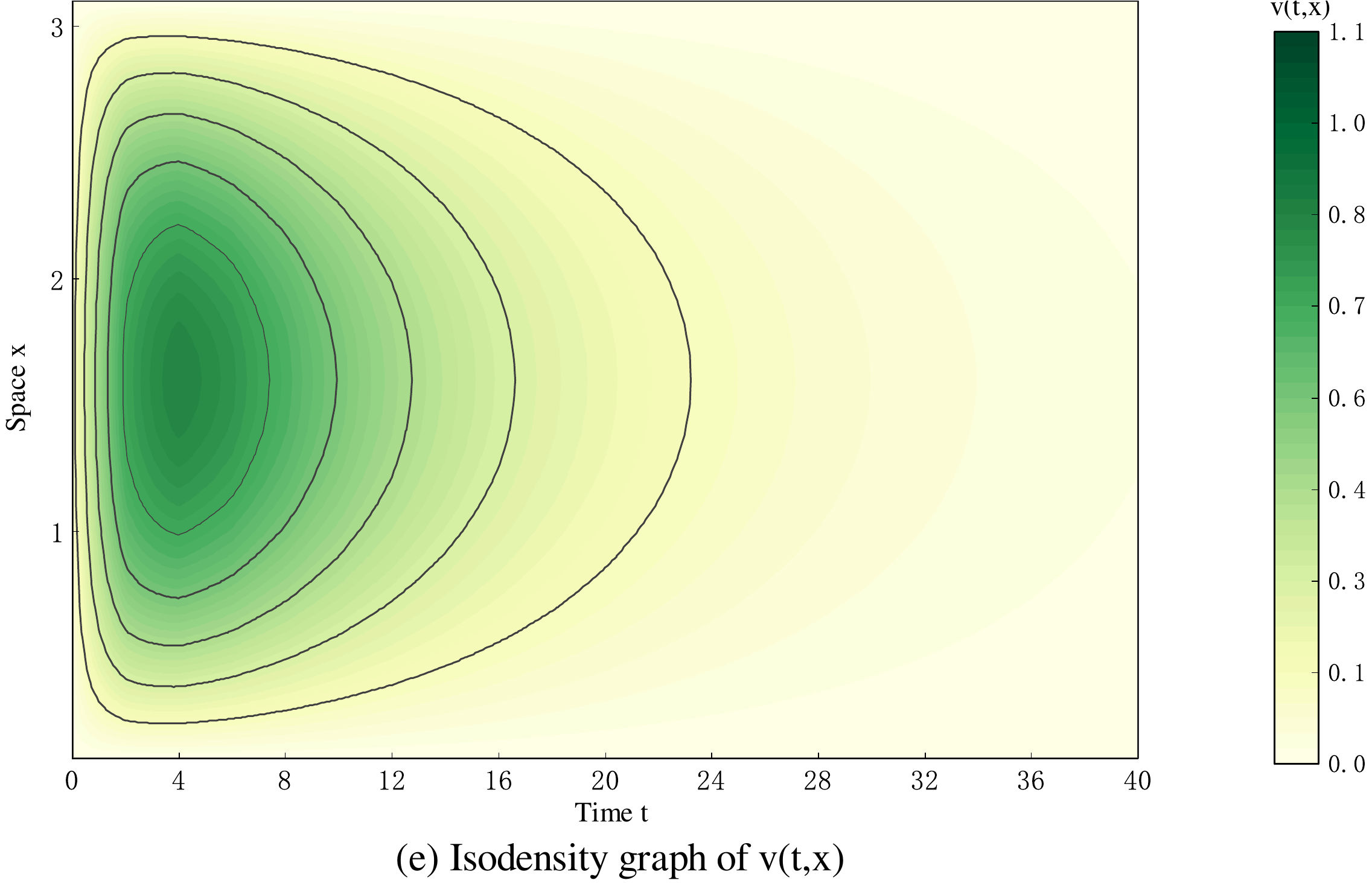}
} }
\subfigure{ {
\includegraphics[width=0.30\textwidth]{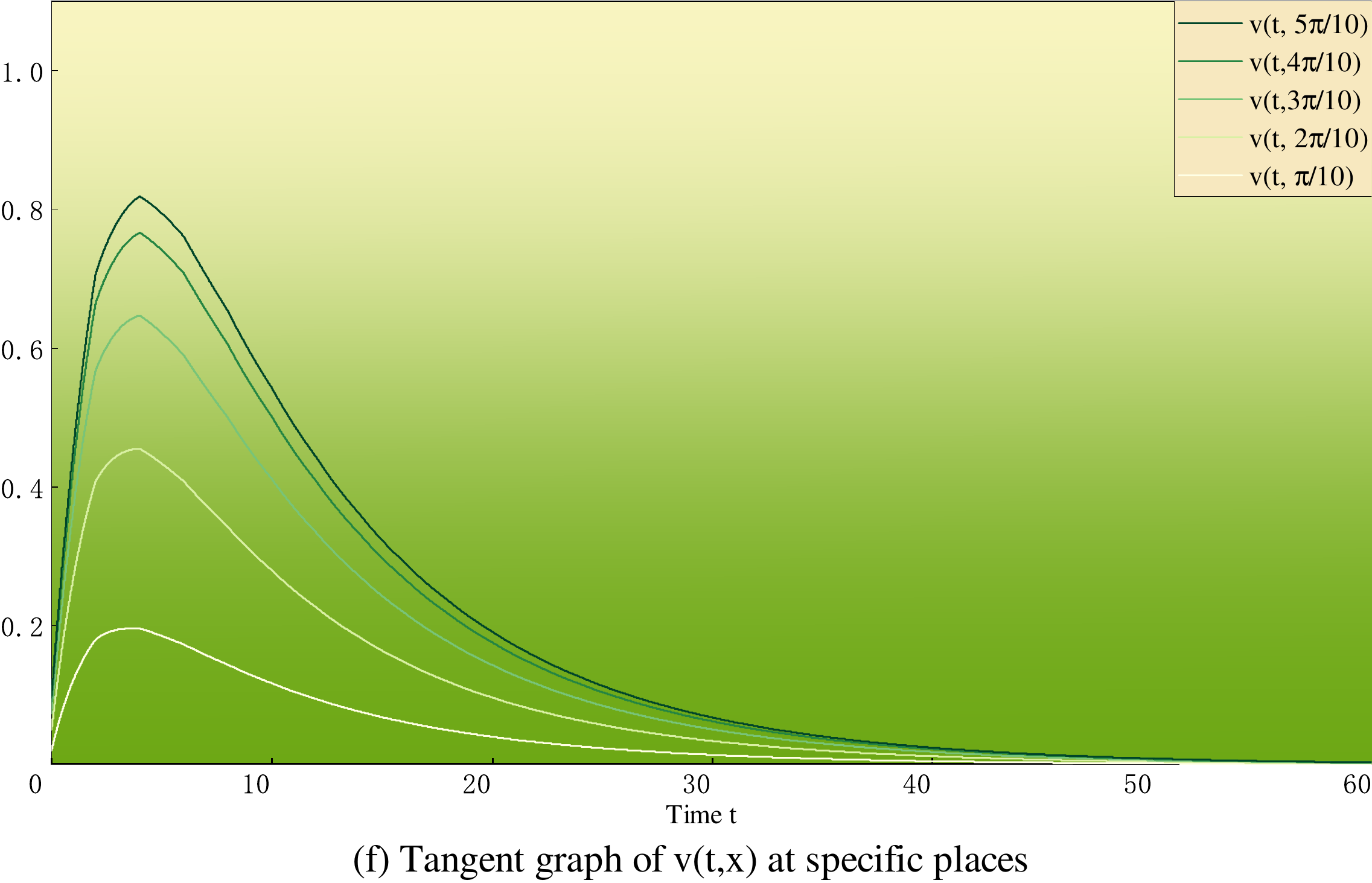}
} }
\caption{When $g(u)=\frac{9u}{10+u}$(with impulse), graphs (a)-(f) show $(u, v)$ more quickly decays to $(0, 0)$.}
\label{B}
\end{figure}
\begin{exm}\label{exm1}
Fix $a_{12}=0.2$, $m_{1}=1$, and $\tau=2$. The periodic impulse function $g$ is taken to be $u$ and $\frac{9u}{10+u}$, respectively.
\end{exm}

Take $g$ to be $u$, i.e., model \eqref{1-2-20} does not have the impulse. Then, it follows from \autoref{theorem 4-3}\textcolor[rgb]{0.00,0.00,1.00}{$(3)$} that $\lambda_{1}\big(g'(0)\big)=\lambda_{1}(1)=0.228$. Thus, \autoref{theorem 3-1} gives that with the progression of time, both $u(t,x)$ and $v(t,x)$ converge uniformly to $0$. In fact, it can be seen from \autoref{A}\textcolor[rgb]{0.00,0.00,1.00}{(a, c, d, f)}  that the bacteria and infected individuals at each location tend together to zero over time. This conclusion is the same as given by \autoref{theorem 3-1}.

Take $g$ to be $\frac{9u}{10+u}$. With the help of \autoref{theorem 4-2}\textcolor[rgb]{0.00,0.00,1.00}{$(1)$}, one can then obtain that $\lambda_{1}\big(g'(0)\big)=\lambda_{1}(0.9)>\lambda_{1}(1)=0.228$. Therefore, \autoref{theorem 3-1} yields that $(0,0)$ is
globally asymptotically stable. Actually, observing \autoref{B} gives us that $u(t,x)$ drops quickly at the pulse point and eventually tends to zero, and $v(t,x)$ also goes to this process, except that it doesn't drop as fast as $u(t,x)$ at the pulse point. In conclusion, the diseases are extinct in the presence of pulses. This conclusion is consistent with that given in \autoref{theorem 3-1}.

A comparison of \autoref{A}\textcolor[rgb]{0.00,0.00,1.00}{(e)} and \autoref{B}\textcolor[rgb]{0.00,0.00,1.00}{(e)} shows that the pulse intervention can reduce the peak in the number of individuals with the infectious diseases and can shorten the time to reach the peak. In addition, the pulse intervention can also shorten the time to extinction of the diseases by observing \autoref{A}\textcolor[rgb]{0.00,0.00,1.00}{(f)} and \autoref{B}\textcolor[rgb]{0.00,0.00,1.00}{(f)}.

\begin{exm}\label{exm2}
Fix $a_{12}=0.9$, $m_{1}=9$, and $\tau=5$. The periodic impulse function $g$ is taken to be $u$ and $\frac{9u}{10+u}$, respectively.
\end{exm}
\begin{figure}[ht]
\centering
\subfigure{ {
\includegraphics[width=0.30\textwidth]{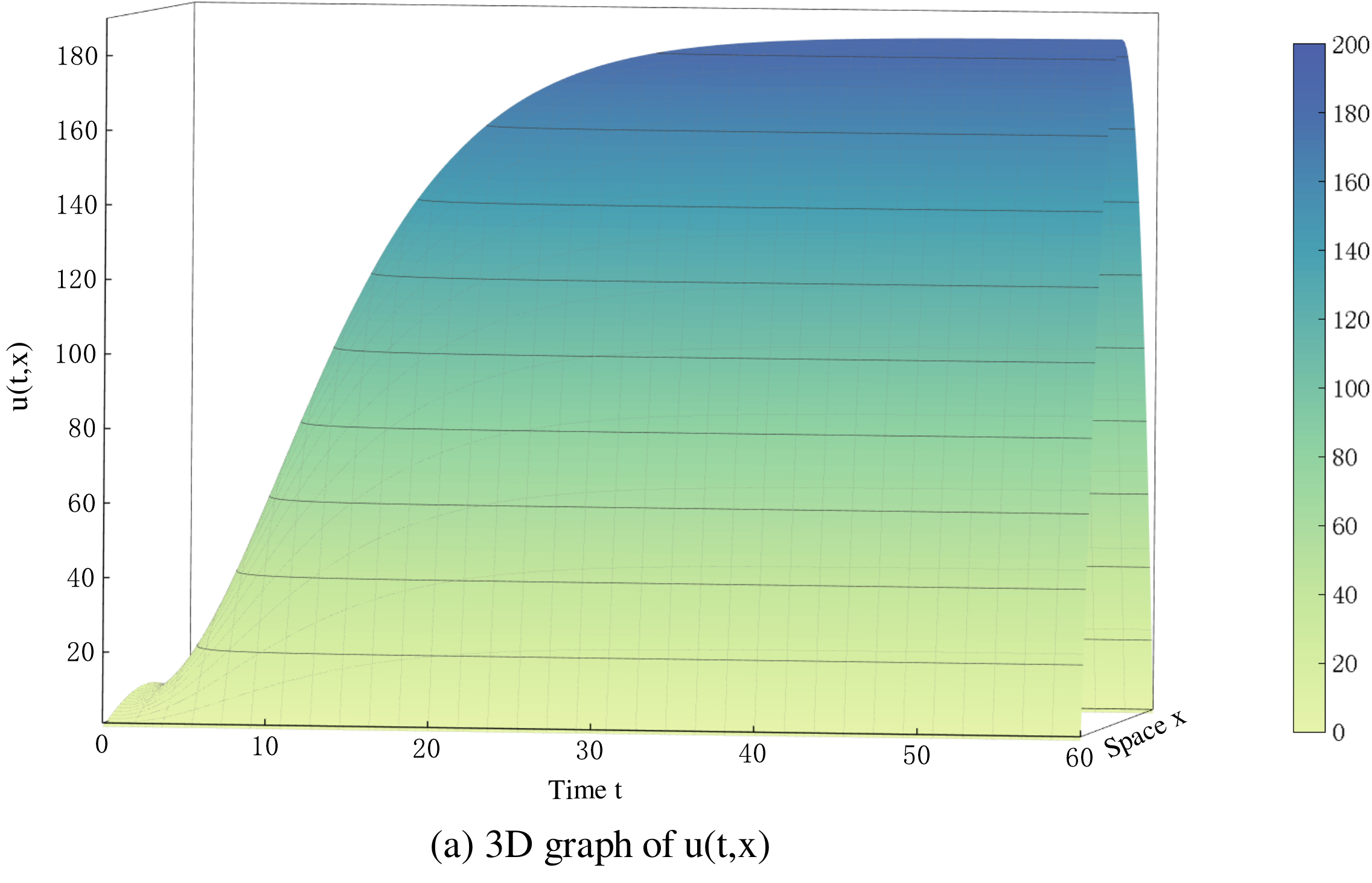}
} }
\subfigure{ {
\includegraphics[width=0.30\textwidth]{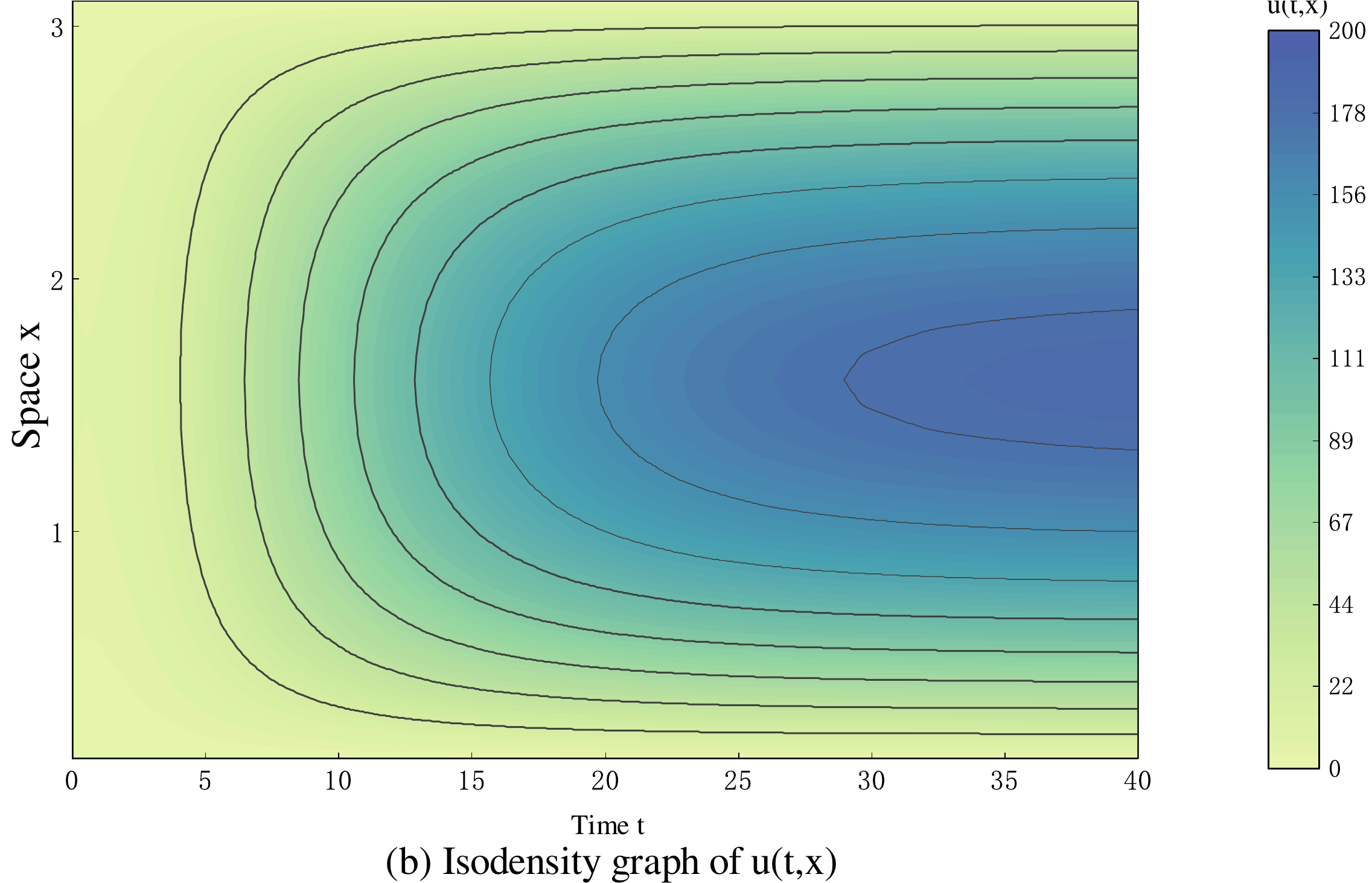}
} }
\subfigure{ {
\includegraphics[width=0.30\textwidth]{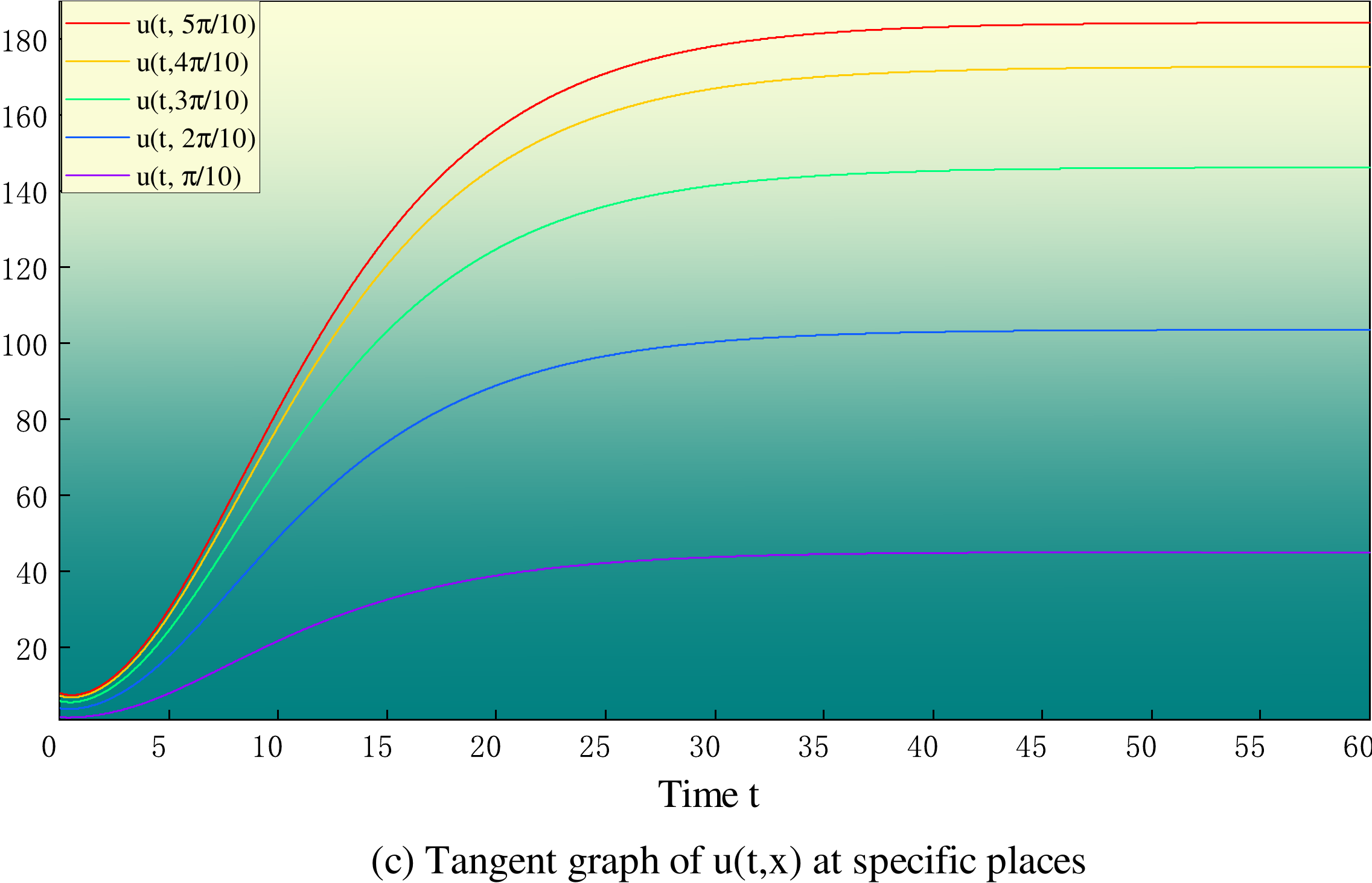}
} }
\subfigure{ {
\includegraphics[width=0.30\textwidth]{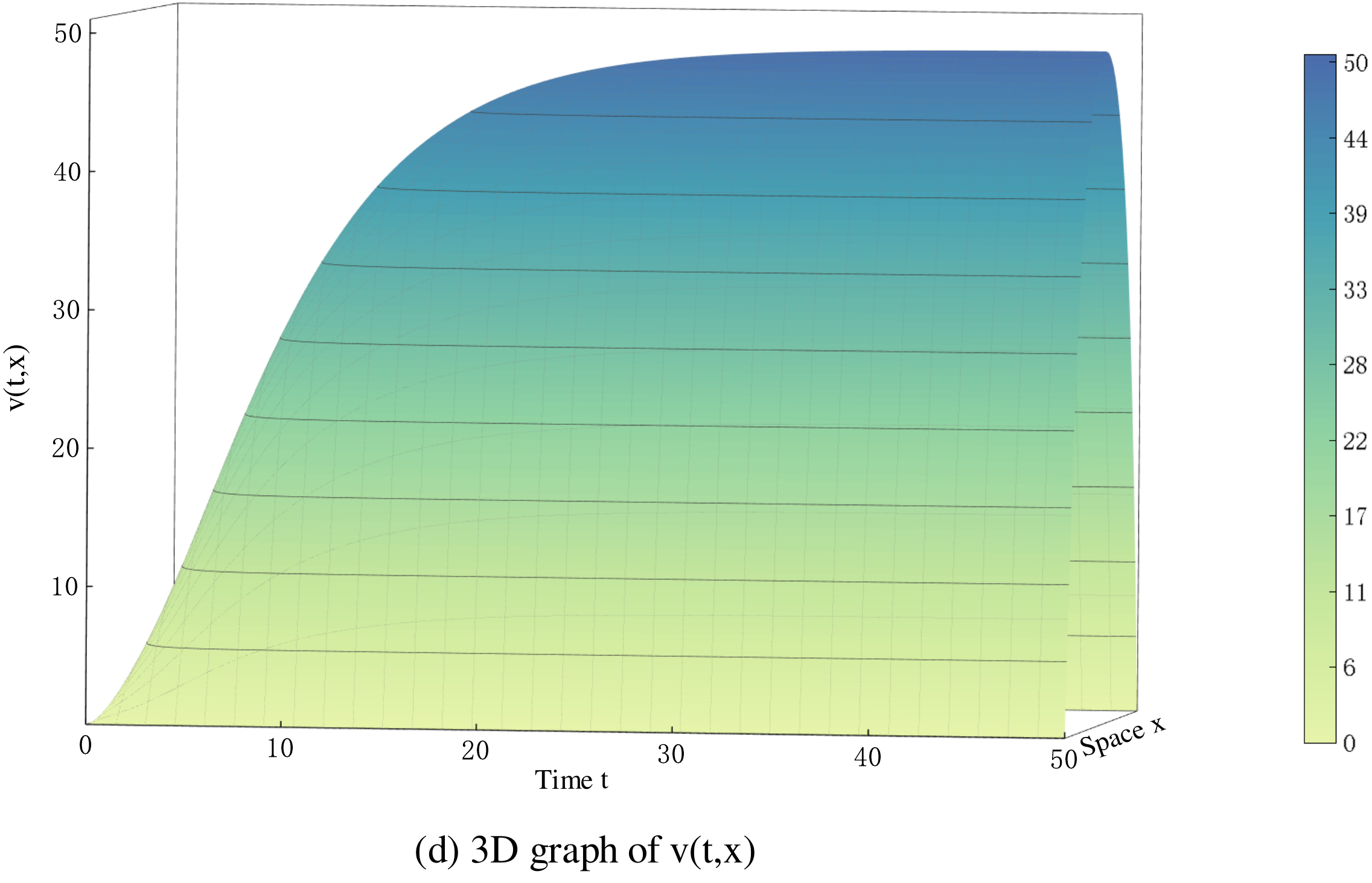}
} }
\subfigure{ {
\includegraphics[width=0.30\textwidth]{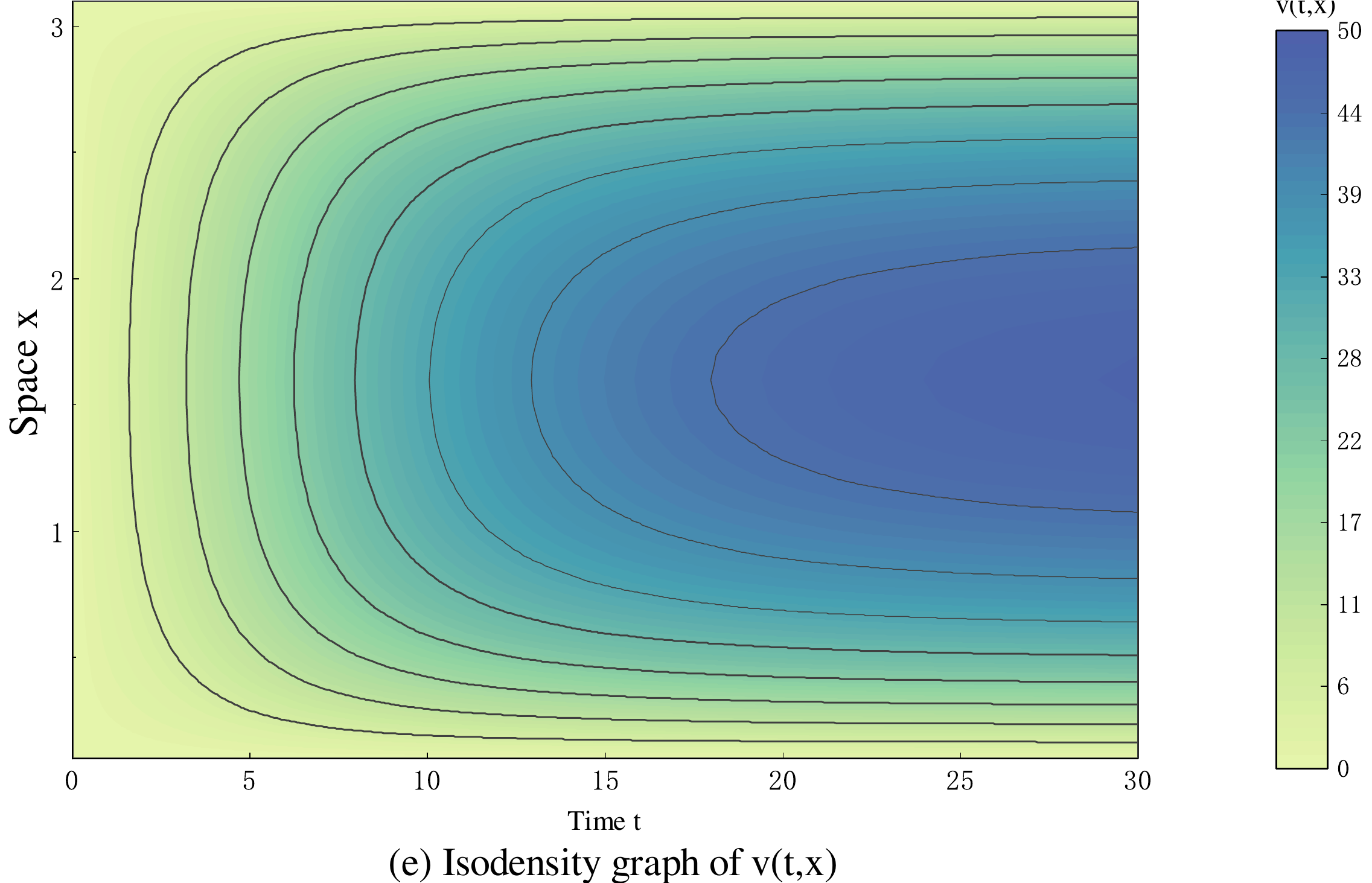}
} }
\subfigure{ {
\includegraphics[width=0.30\textwidth]{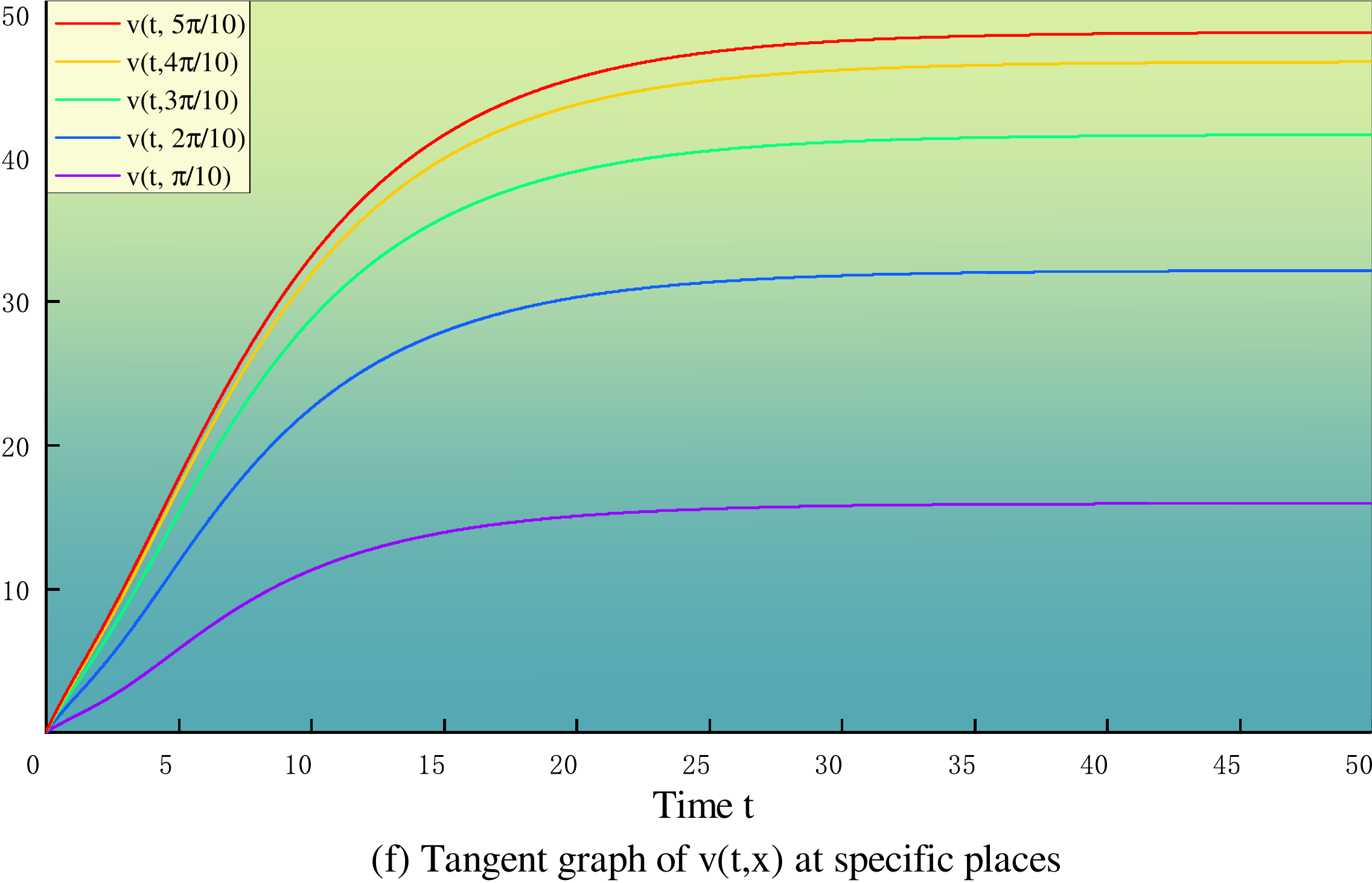}
} }
\caption{When $g(u)=u$(without impulse), graphs (a)-(f) exhibit $u$ and $v$ converge to steady states.}
\label{C}
\end{figure}
\begin{figure}[!ht]
\centering
\subfigure{ {
\includegraphics[width=0.30\textwidth]{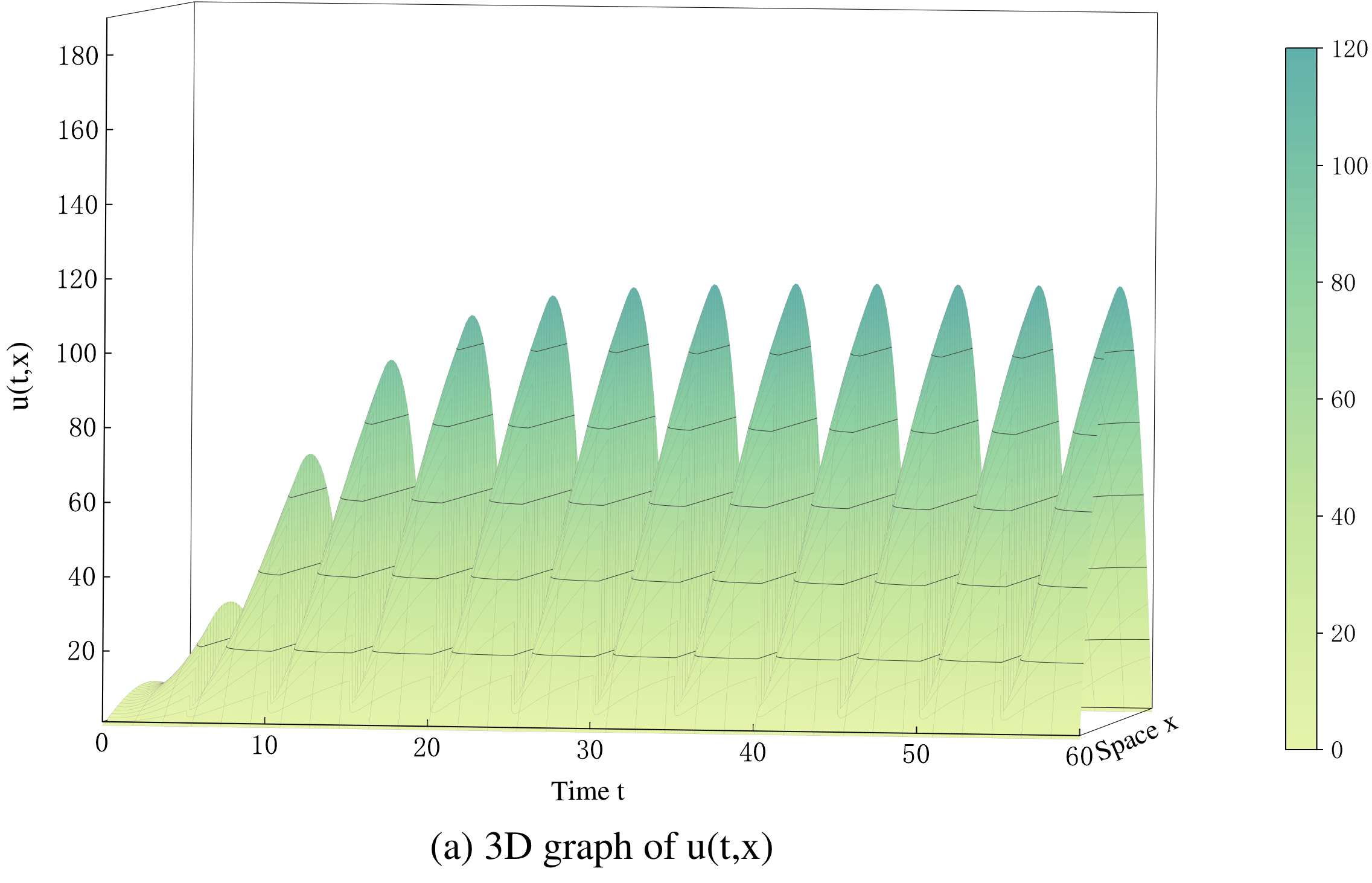}
} }
\subfigure{ {
\includegraphics[width=0.30\textwidth]{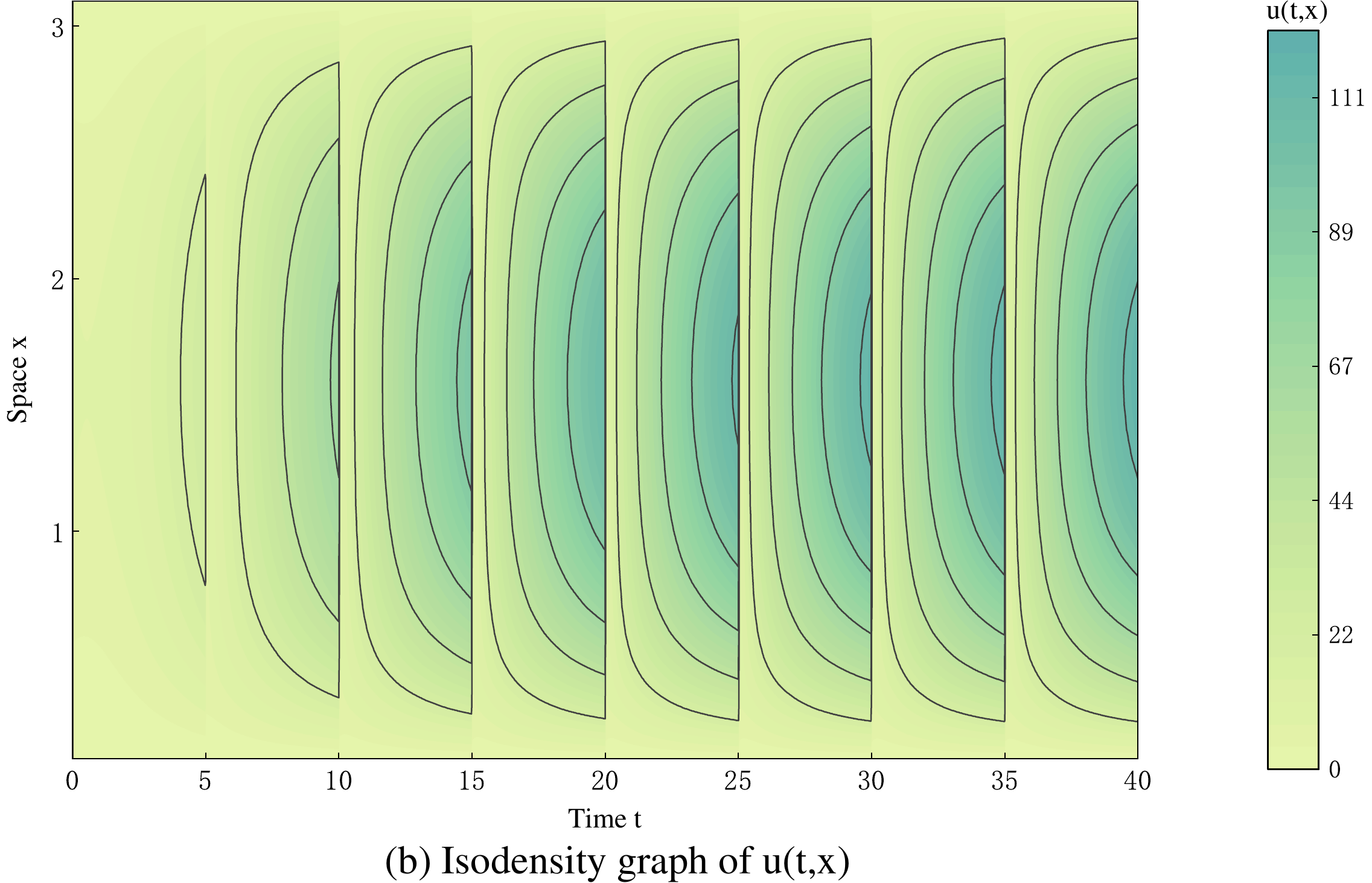}
} }
\subfigure{ {
\includegraphics[width=0.30\textwidth]{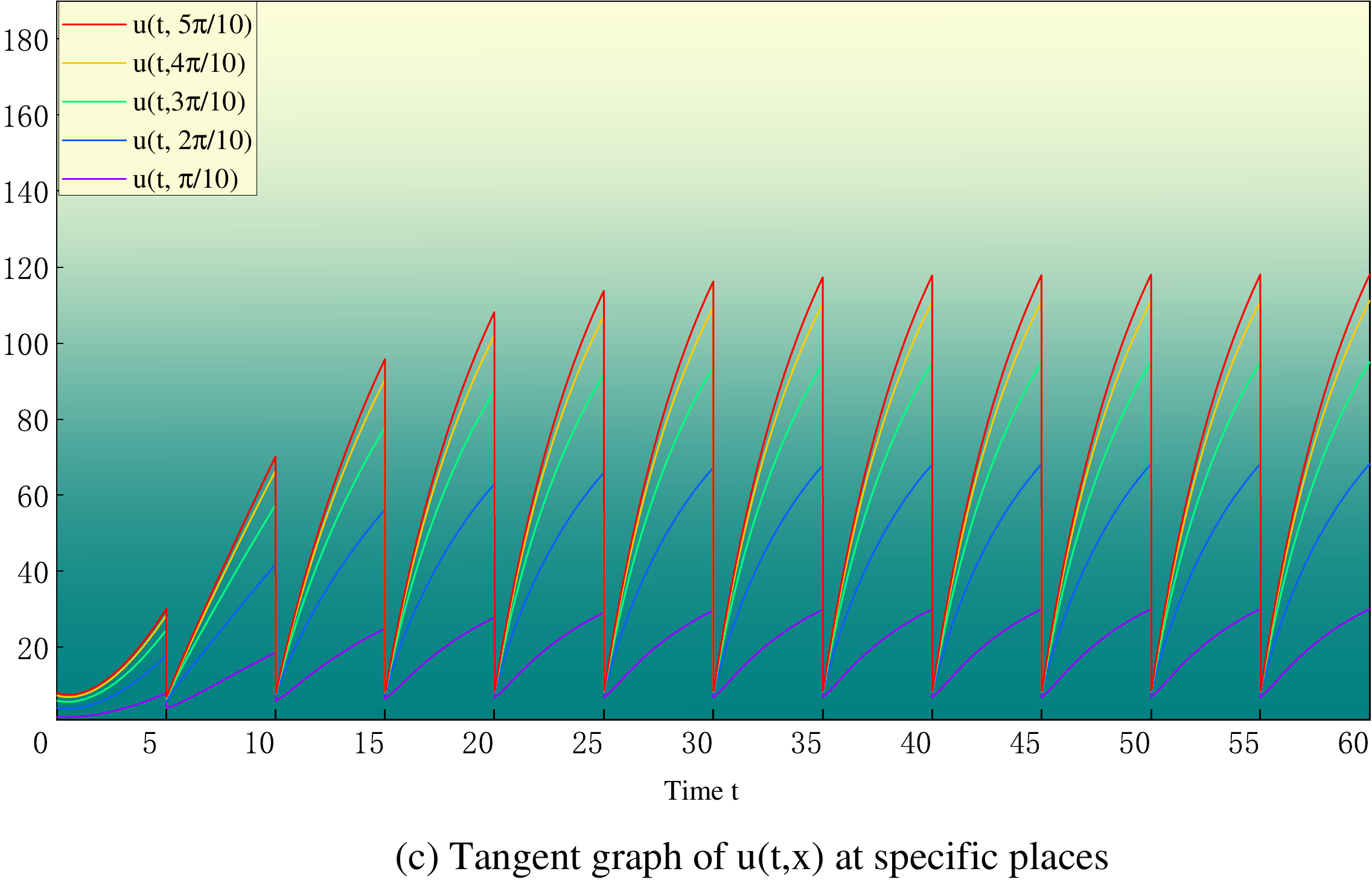}
} }
\subfigure{ {
\includegraphics[width=0.30\textwidth]{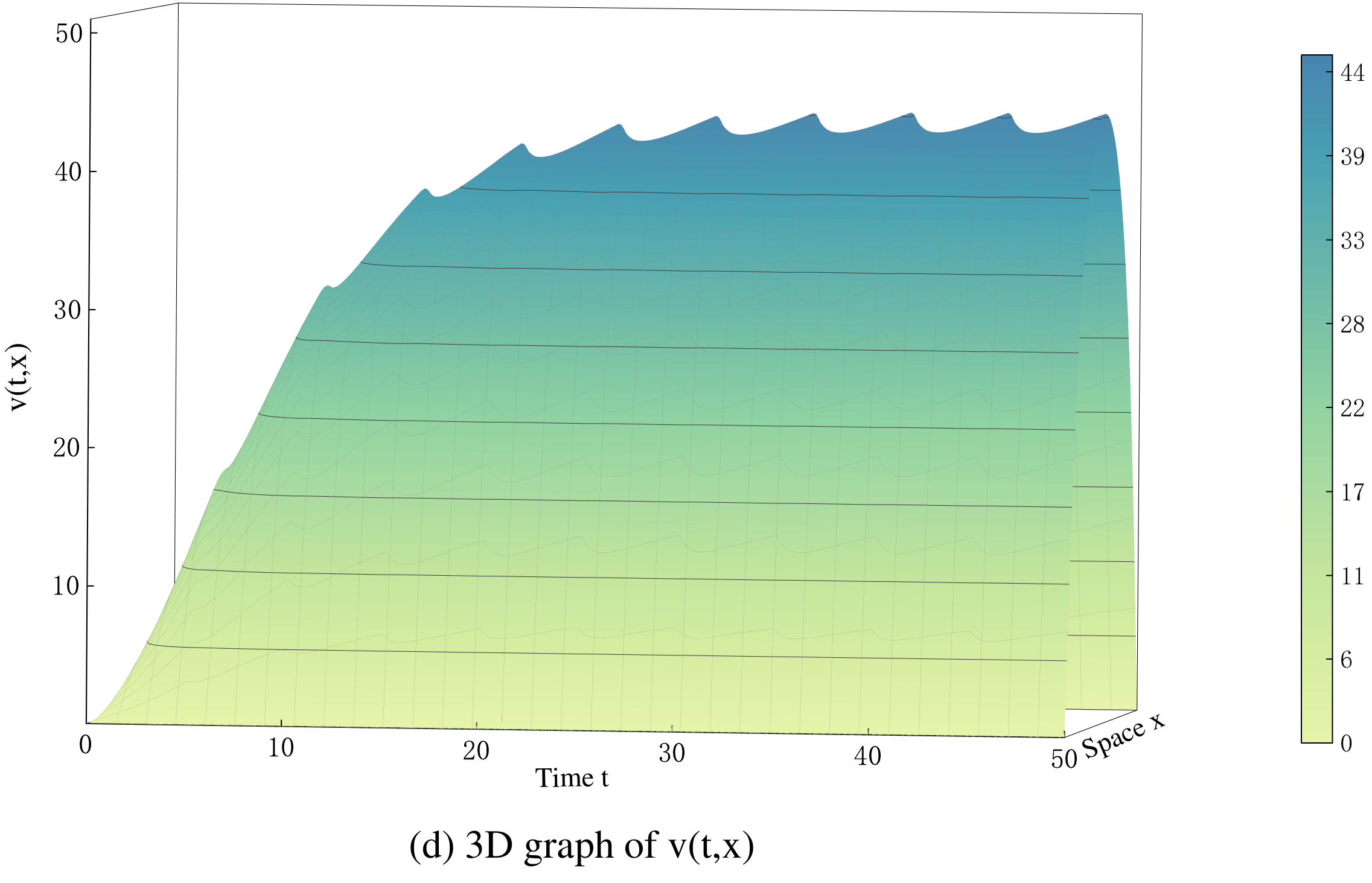}
} }
\subfigure{ {
\includegraphics[width=0.30\textwidth]{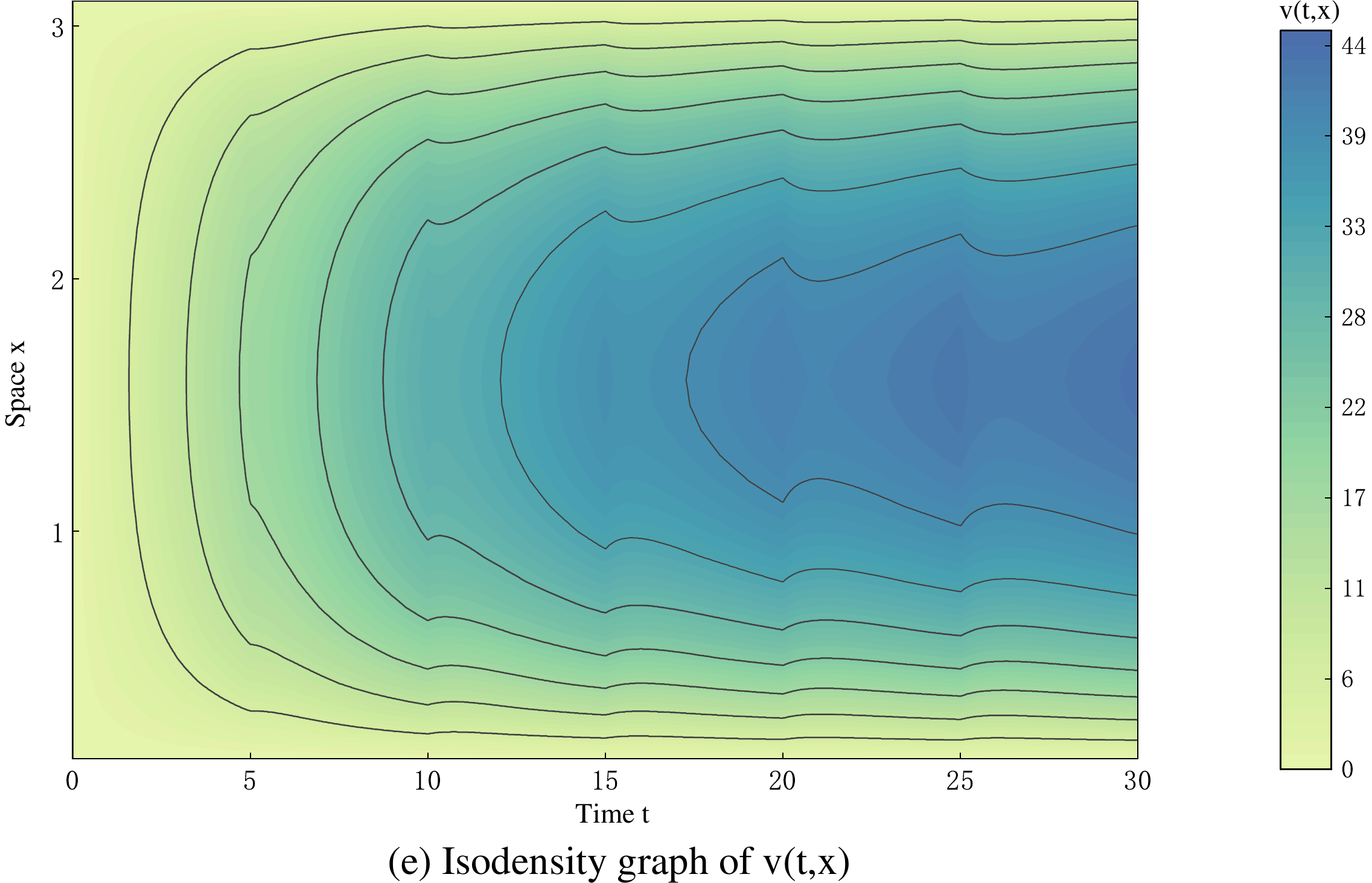}
} }
\subfigure{ {
\includegraphics[width=0.30\textwidth]{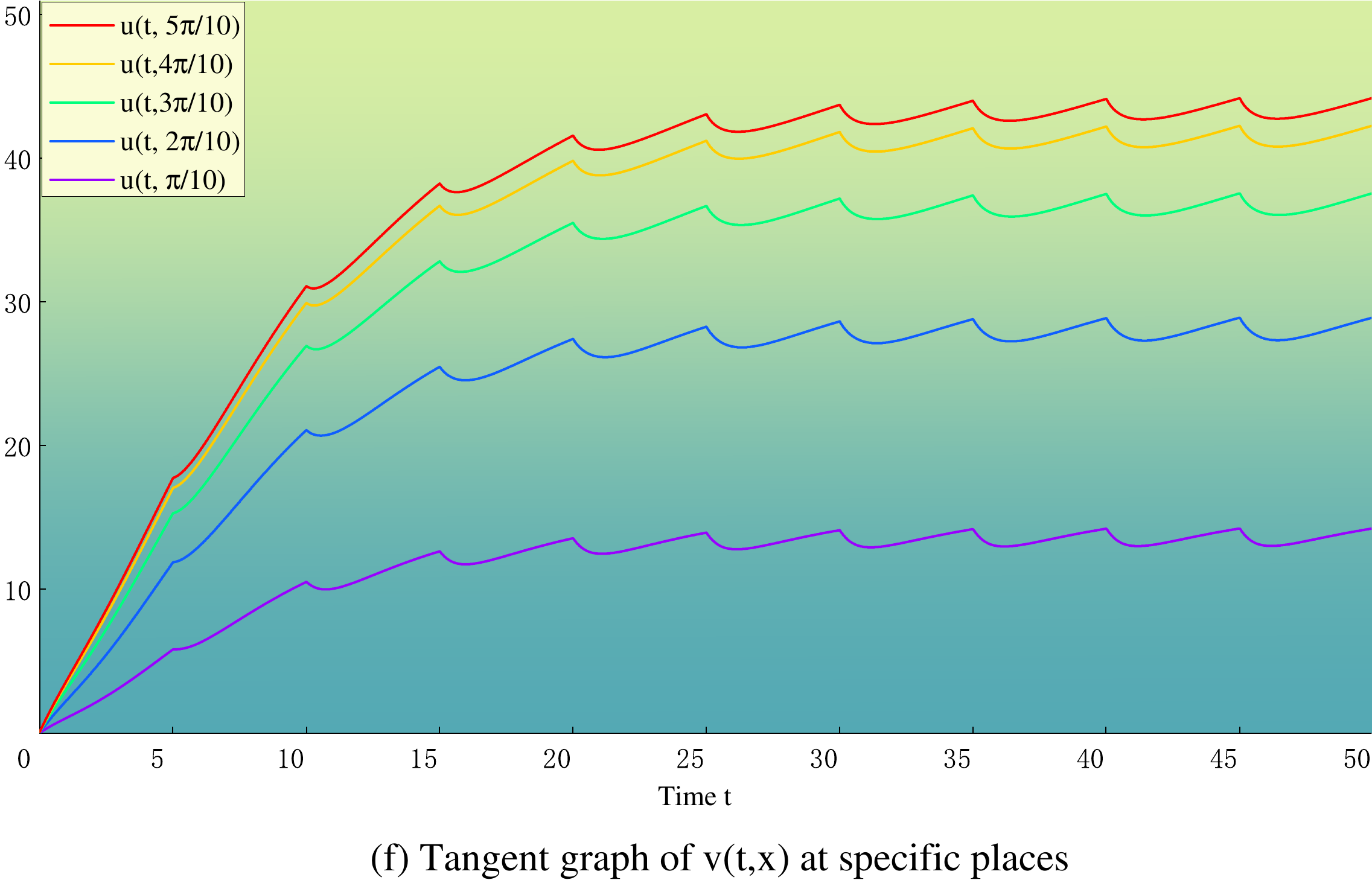}
} }
\caption{\small{When $g(u)=\frac{9u}{10+u}$(with impulse), graphs (a)-(f) show $u$ and $v$ converge to smaller steady states. }}
\label{D}
\end{figure}
In the absence of impulse intervention in model \eqref{1-2-20}, \autoref{theorem 4-3}\textcolor[rgb]{0.00,0.00,1.00}{$(3)$} gives that $\lambda_{1}\big(g'(0)\big)=\lambda_{1}(1)=-0.306$. Then, \textcolor[rgb]{0.00,0.00,1.00}{Theorems} \ref{theorem 3-2} and \ref{theorem 3-3} yield that model \eqref{1-2-20} has a unique positive steady state solution, and it is globally asymptotically stable. Indeed, it can be noticed from \autoref{C}\textcolor[rgb]{0.00,0.00,1.00}{(a, b, d, e)}  that the bacteria and infected individuals converge to positive distributions over time, respectively. This conclusion is the same as given by \autoref{theorem 3-3}.

Now, consider model \eqref{1-2-20} with periodic pulse intervention $g(u)=\frac{9u}{10+u}$. In contrast to \autoref{C}, \autoref{D} shows that the bacteria and infected individuals tend to positive distributions of time period $5$ over time, respectively, and that both bacteria and infected individuals possess lower steady states. Hence, when infectious diseases of faecal-oral transmission are long-standing, periodic impulse intervention can help to reduce the number of infected persons.

\textcolor[rgb]{0.00,0.00,1.00}{Examples} \ref{exm1} and \ref{exm2} show that periodic impulses play a positive role in the control of faecally-orally transmitted diseases.
\subsection{The influence of regional evolution rate}
To further understand the influence of regional evolution rate on the spread of the diseases, all numerical simulations presented in this subsection are done in case of no pulse. Some of the parameters are chosen as $a_{22}=0.1$, $\tau=2$, and $a_{1}=10$. Other of the parameters are given in each of the following Examples to show the different dynamical behaviours. The forthcoming \textcolor[rgb]{0.00,0.00,1.00}{Examples} \ref{exm3} and \ref{exm4} will give the influence of regional evolution rate on the extinction and persistence of the faecal-oral transmission diseases, respectively.
\begin{exm}\label{exm3}
Fix $d_{1}=0.05$, $d_{2}=0.1$, $a_{11}=0.17$, $a_{12}=0.2$, and $m_{1}=1.5$. The regional evolution rate $\rho(t)$ is taken to be $1$ and $0.7e^{2(1-cos(\pi t))}$, respectively.
\end{exm}
\begin{figure}[!ht]
\centering
\subfigure{ {
\includegraphics[width=0.30\textwidth]{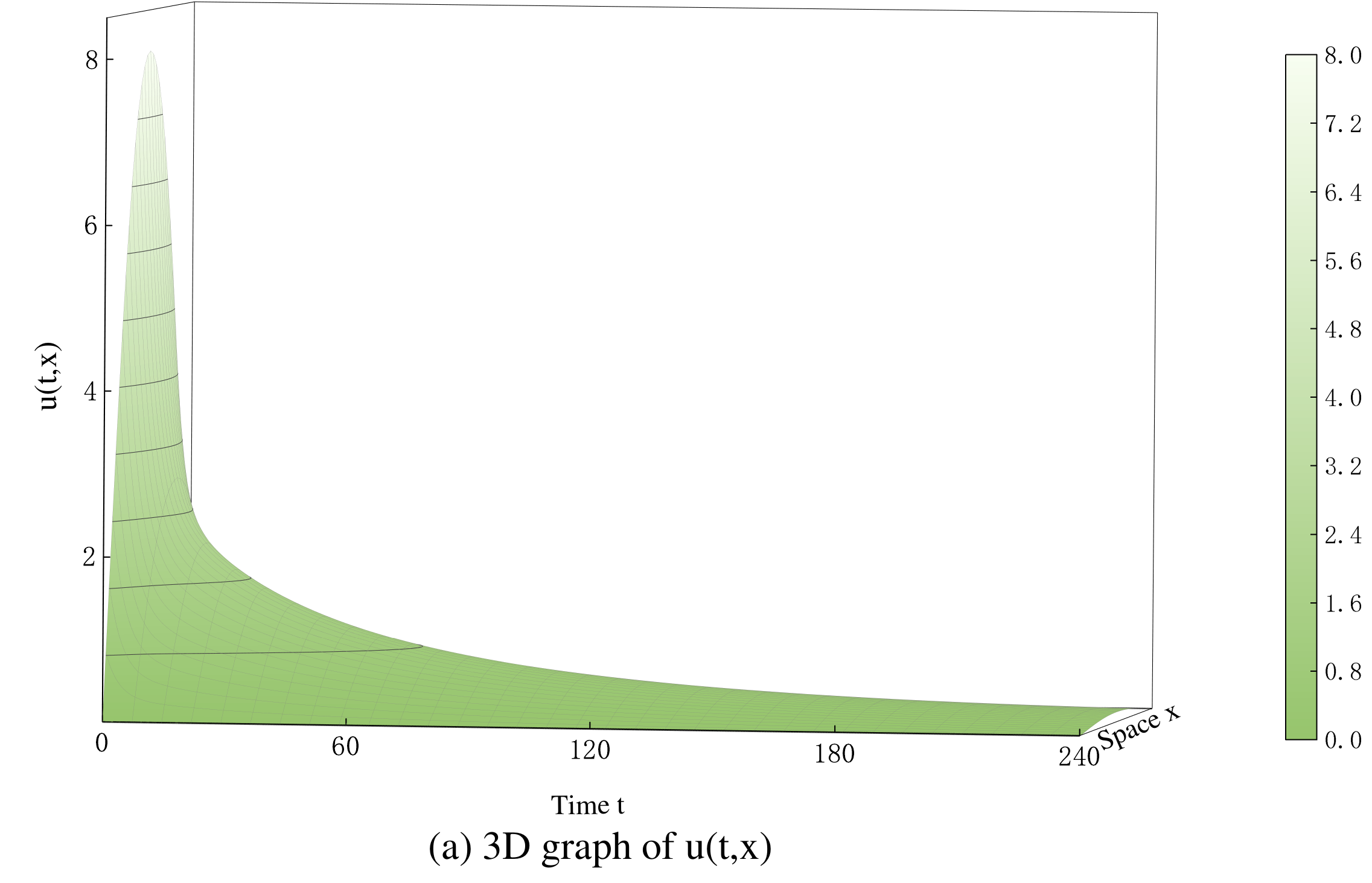}
} }
\subfigure{ {
\includegraphics[width=0.30\textwidth]{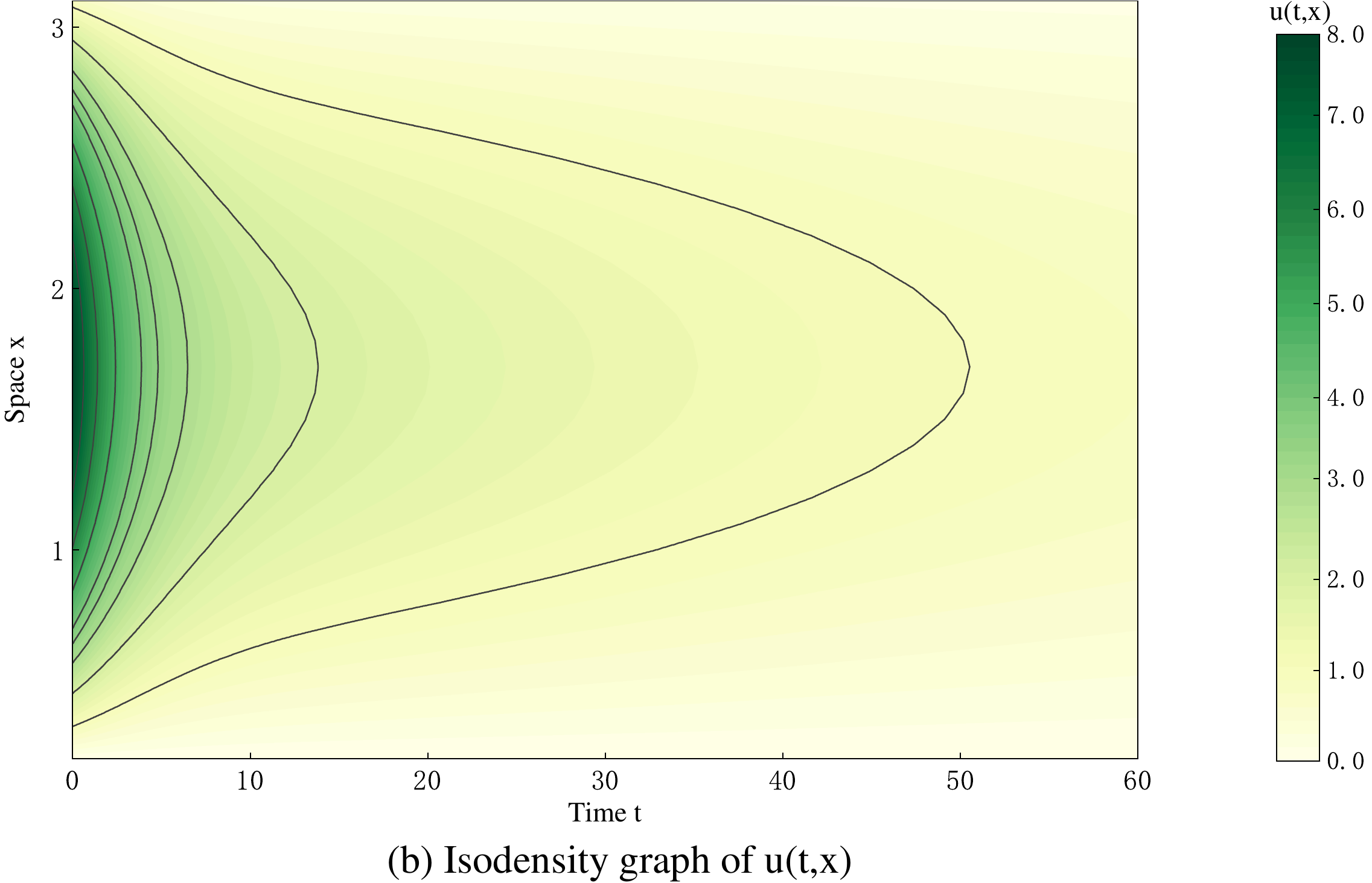}
} }
\subfigure{ {
\includegraphics[width=0.30\textwidth]{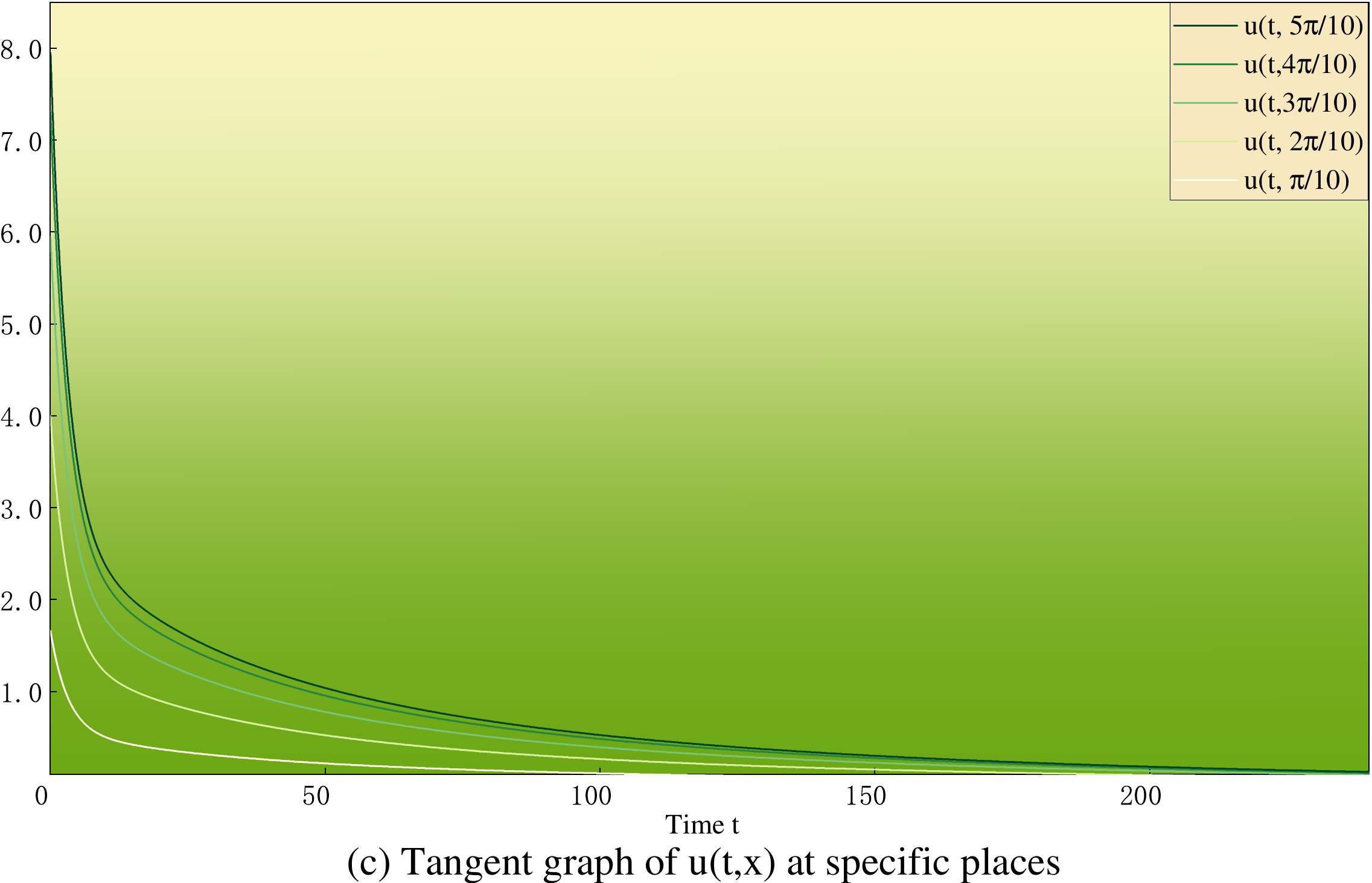}
} }
\subfigure{ {
\includegraphics[width=0.30\textwidth]{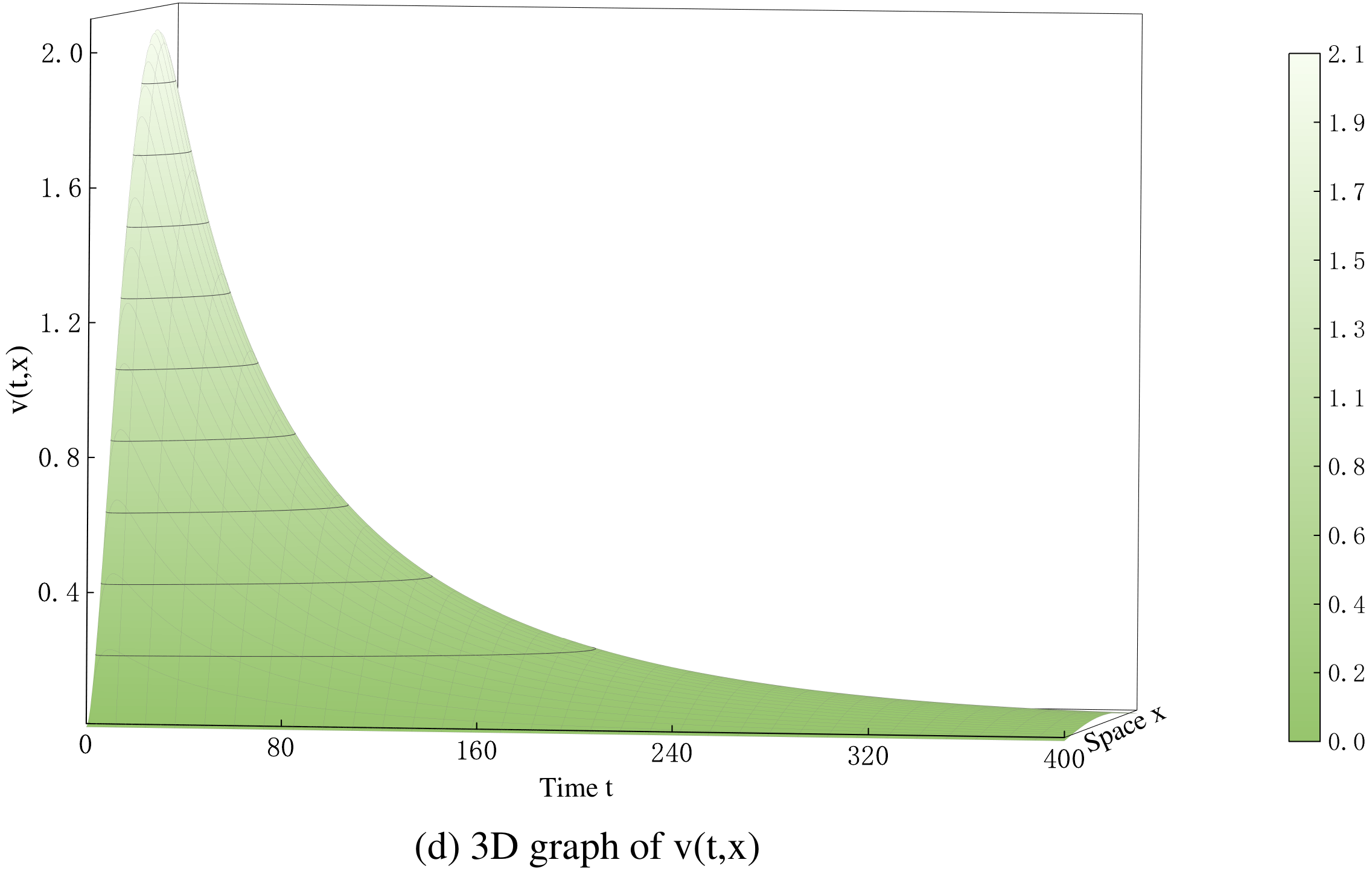}
} }
\subfigure{ {
\includegraphics[width=0.30\textwidth]{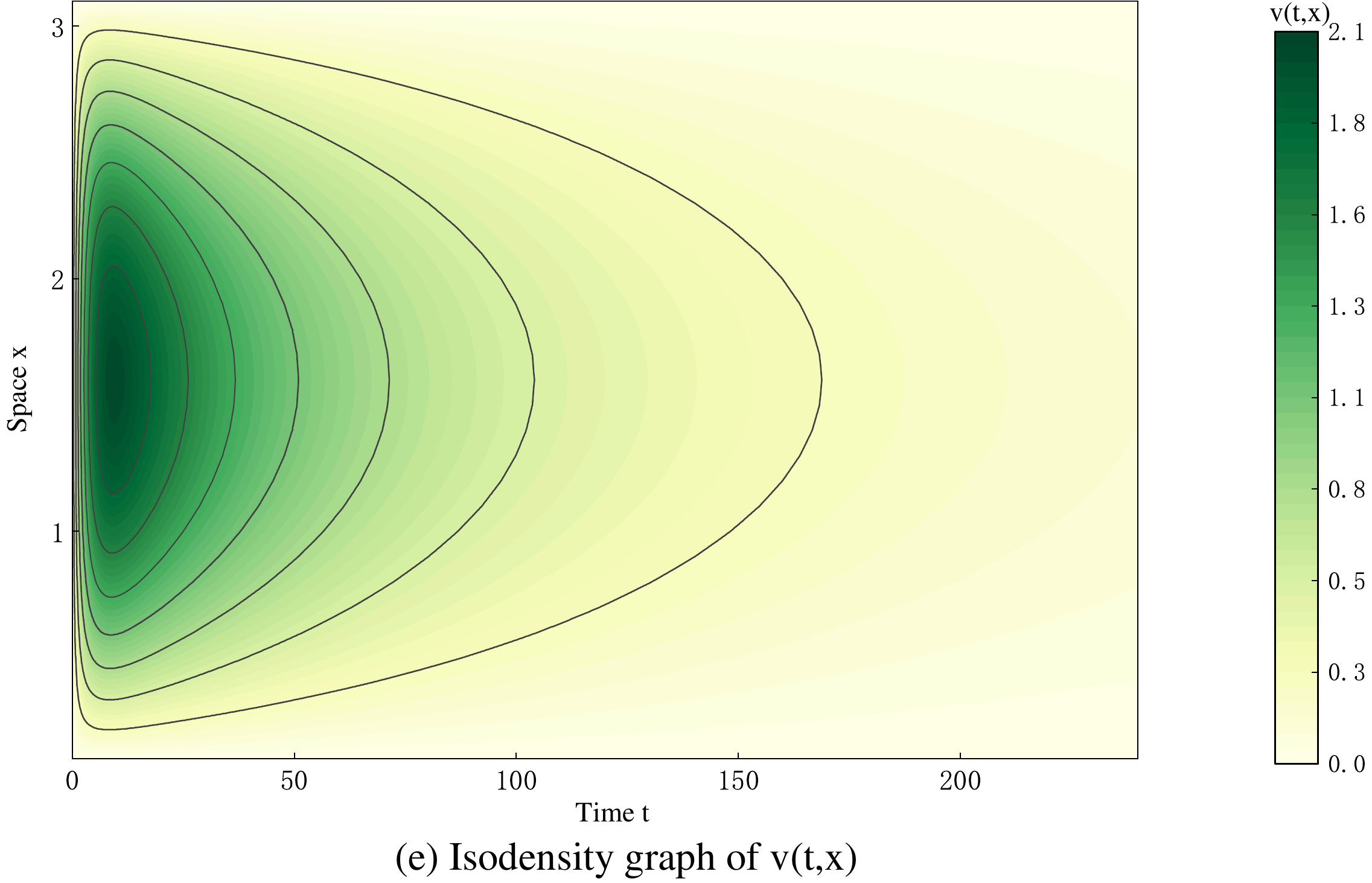}
} }
\subfigure{ {
\includegraphics[width=0.30\textwidth]{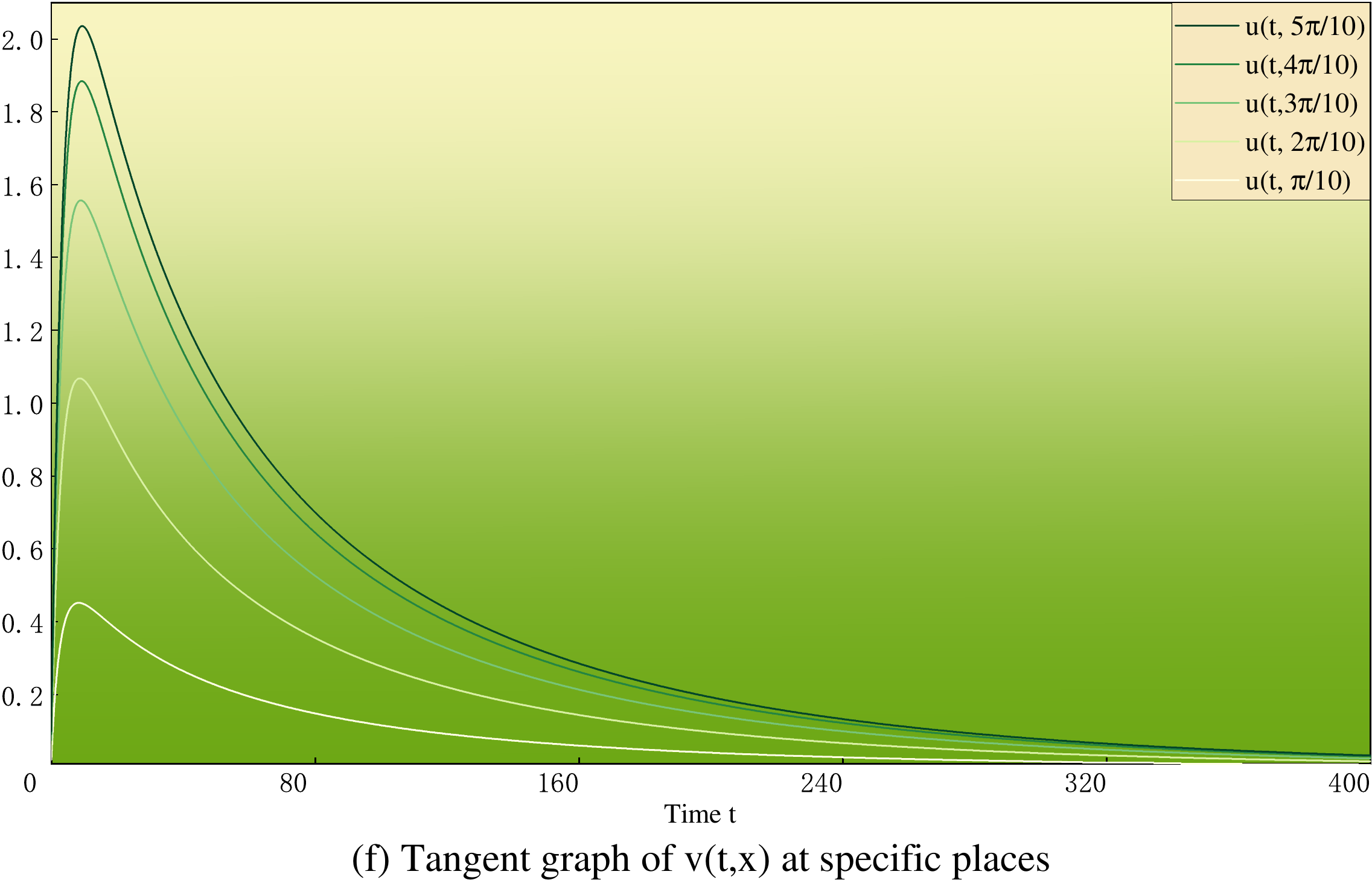}
} }
\caption{When $\rho(t)=1$(fixed region), graphs (a)-(f) show that $(u, v)$ decays to $(0, 0)$. }
\label{E}
\end{figure}
\begin{figure}[!ht]
\centering
\subfigure{ {
\includegraphics[width=0.30\textwidth]{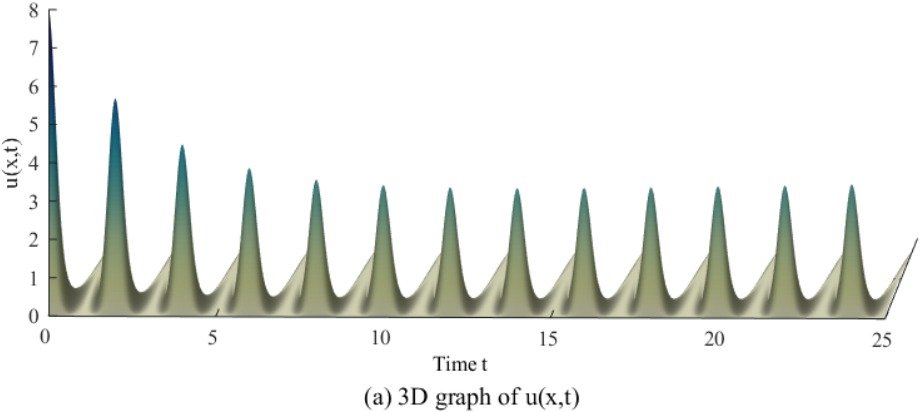}
} }
\subfigure{ {
\includegraphics[width=0.30\textwidth]{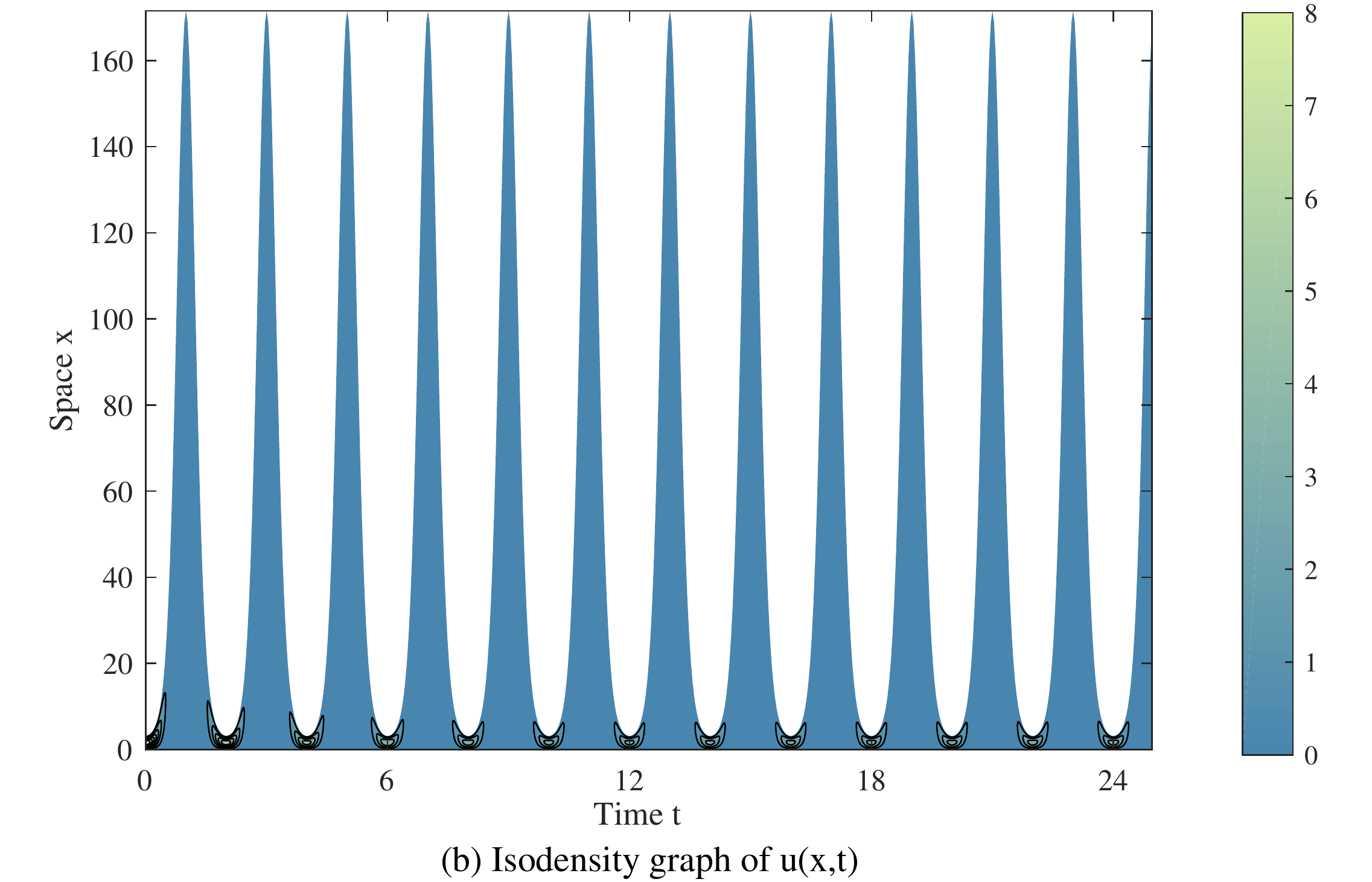}
} }
\subfigure{ {
\includegraphics[width=0.30\textwidth]{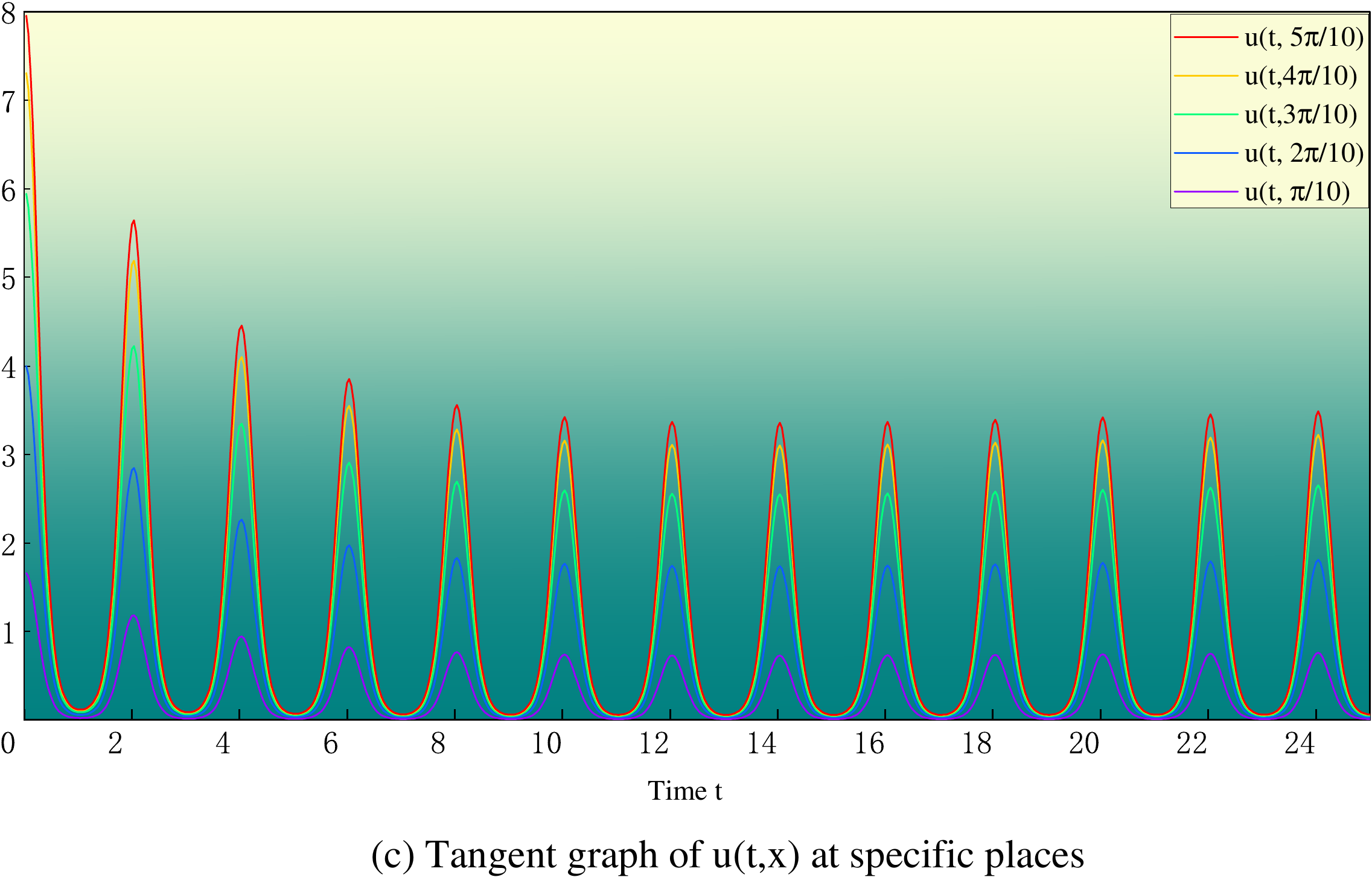}
} }
\subfigure{ {
\includegraphics[width=0.30\textwidth]{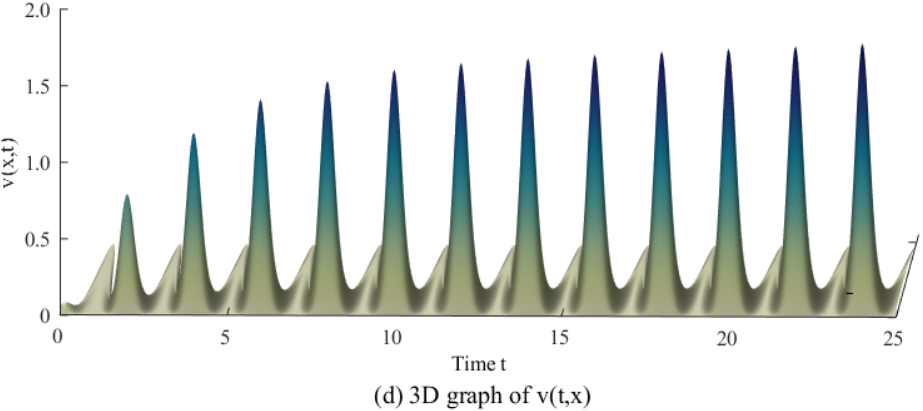}
} }
\subfigure{ {
\includegraphics[width=0.30\textwidth]{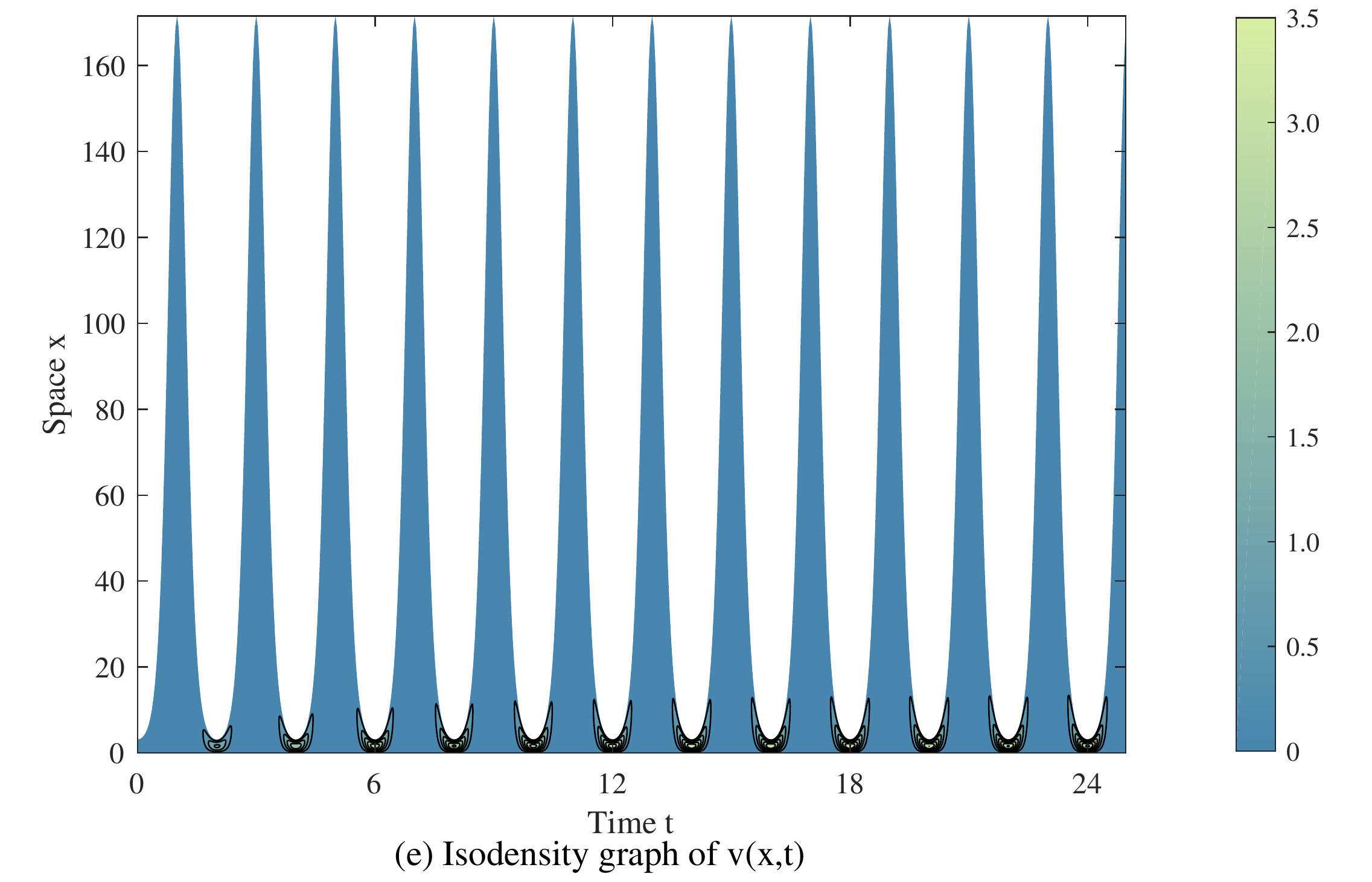}
} }
\subfigure{ {
\includegraphics[width=0.30\textwidth]{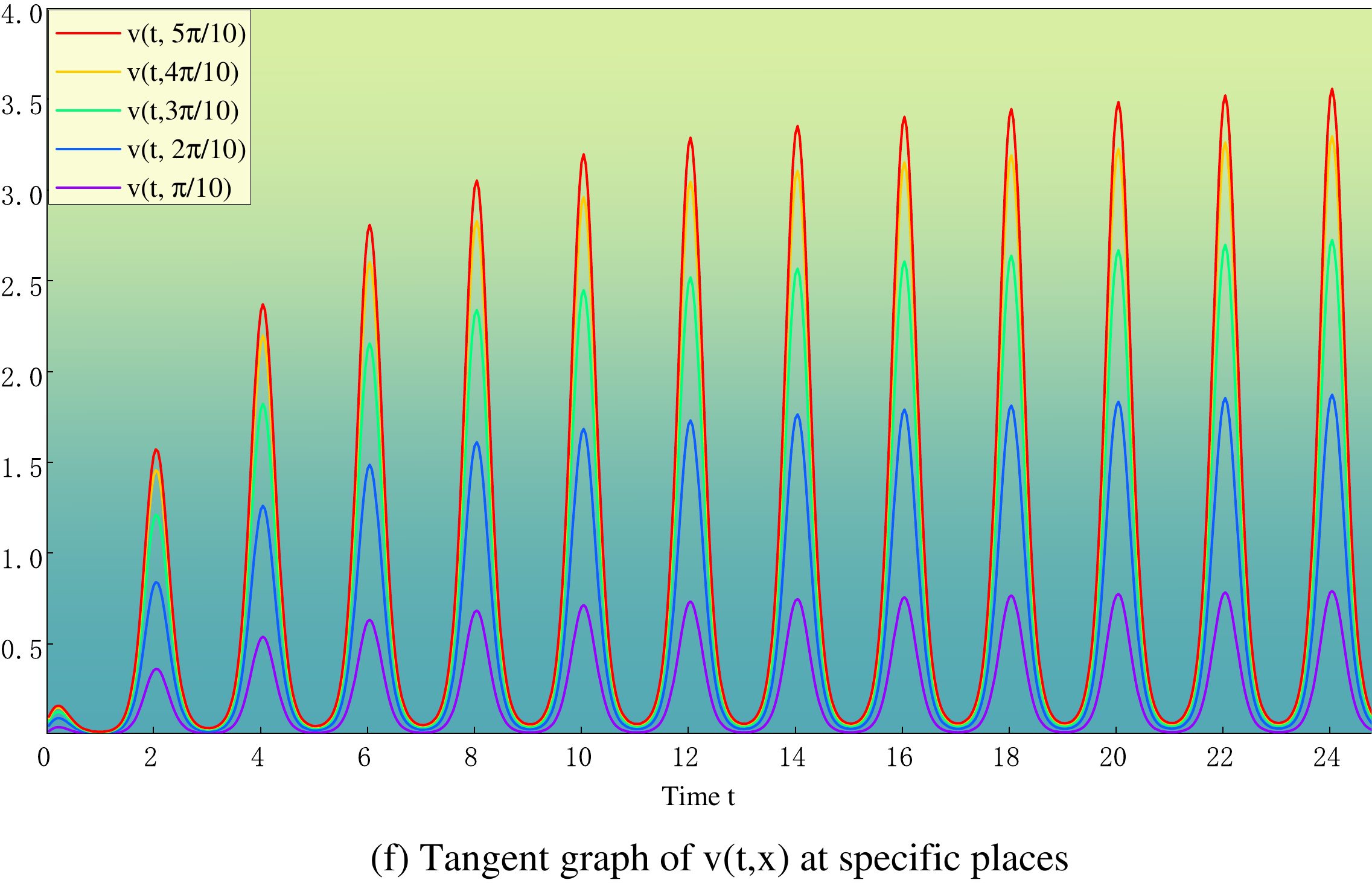}
} }
\caption{\small{When $\rho(t)=e^{2(1-\cos(\pi t))}$(evolving region), graphs (a)-(f) show $u$ and $v$ converge to steady states.}}
\label{F}
\end{figure}

Select regional evolving rate $\rho(t)$ as $1$. Then, it follows from \autoref{theorem 4-3}\textcolor[rgb]{0.00,0.00,1.00}{$(3)$} that $\lambda_{1}\big(\rho(t)\big)=0.037$. This combined with \autoref{theorem 3-1} yields that with the progression of time, both $u(t,x)$ and $v(t,x)$ converge uniformly to $0$. Actually, it can be seen from \autoref{E}\textcolor[rgb]{0.00,0.00,1.00}{(a, c, d, f)}  that the bacteria and infected individuals at each location tend together to zero over time. This conclusion is the same as given by \autoref{theorem 3-1}.

Select regional evolving rate $\rho(t)$ as $e^{2(1-\cos(\pi t))}$. With the help of \autoref{theorem 4-3}\textcolor[rgb]{0.00,0.00,1.00}{$(1)$} that $\lambda_{1}\big(\rho(t)\big)<-0.025$. Therefore, \textcolor[rgb]{0.00,0.00,1.00}{Theorems} \ref{theorem 3-2} and \ref{theorem 3-3} yield that model \eqref{1-2-20} has a unique positive periodic steady state solution, and it is globally asymptotically stable. Actually, it can be noticed from \autoref{F}\textcolor[rgb]{0.00,0.00,1.00}{(a, b, d, e)}  that the bacteria and infected individuals respectively converge to positive periodic distributions over time. This coincides with the conclusion of \autoref{theorem 3-3}.

Fecal-oral diseases become extinct after a period of time, when the populations' habitat is fixed. However, this diseases are persistent, when the habitat of the populations changes with a large evolutionary rate $\rho(t)$($\xoverline{\rho^{-2}}<1$). Thus, larger evolutionary rate is not beneficial for the control and prevention of the diseases.
\begin{exm}\label{exm4}
Fix $d_{1}=0.75$, $d_{2}=0.75$, $a_{11}=0.1$, $a_{12}=0.5$, and $m_{1}=20$. The regional evolution rate $\rho(t)$ is taken to be $1$ and $0.7e^{-0.15(1-cos(\pi t))}$, respectively.
\end{exm}
\begin{figure}[!ht]
\centering
\subfigure{ {
\includegraphics[width=0.30\textwidth]{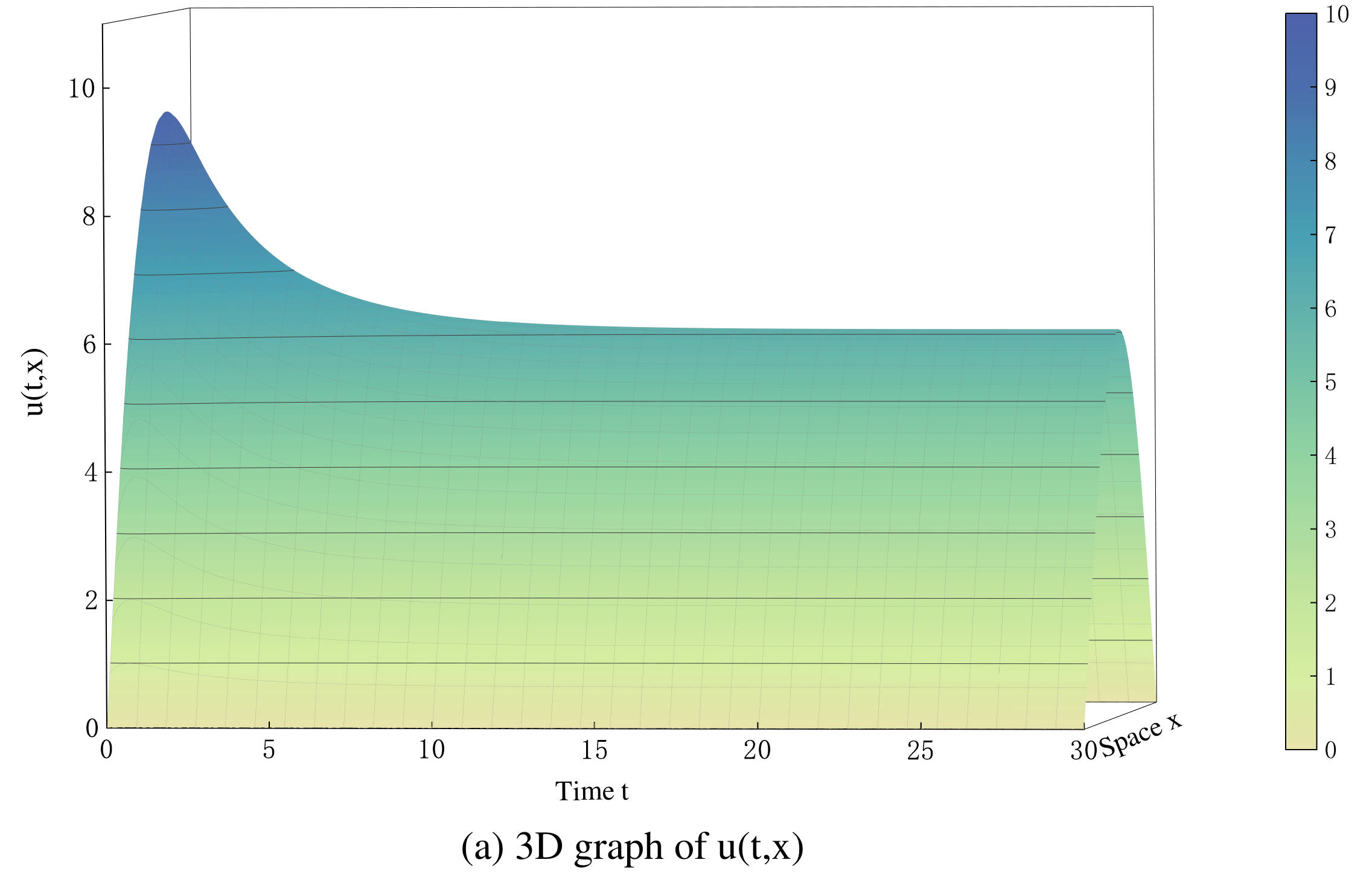}
} }
\subfigure{ {
\includegraphics[width=0.30\textwidth]{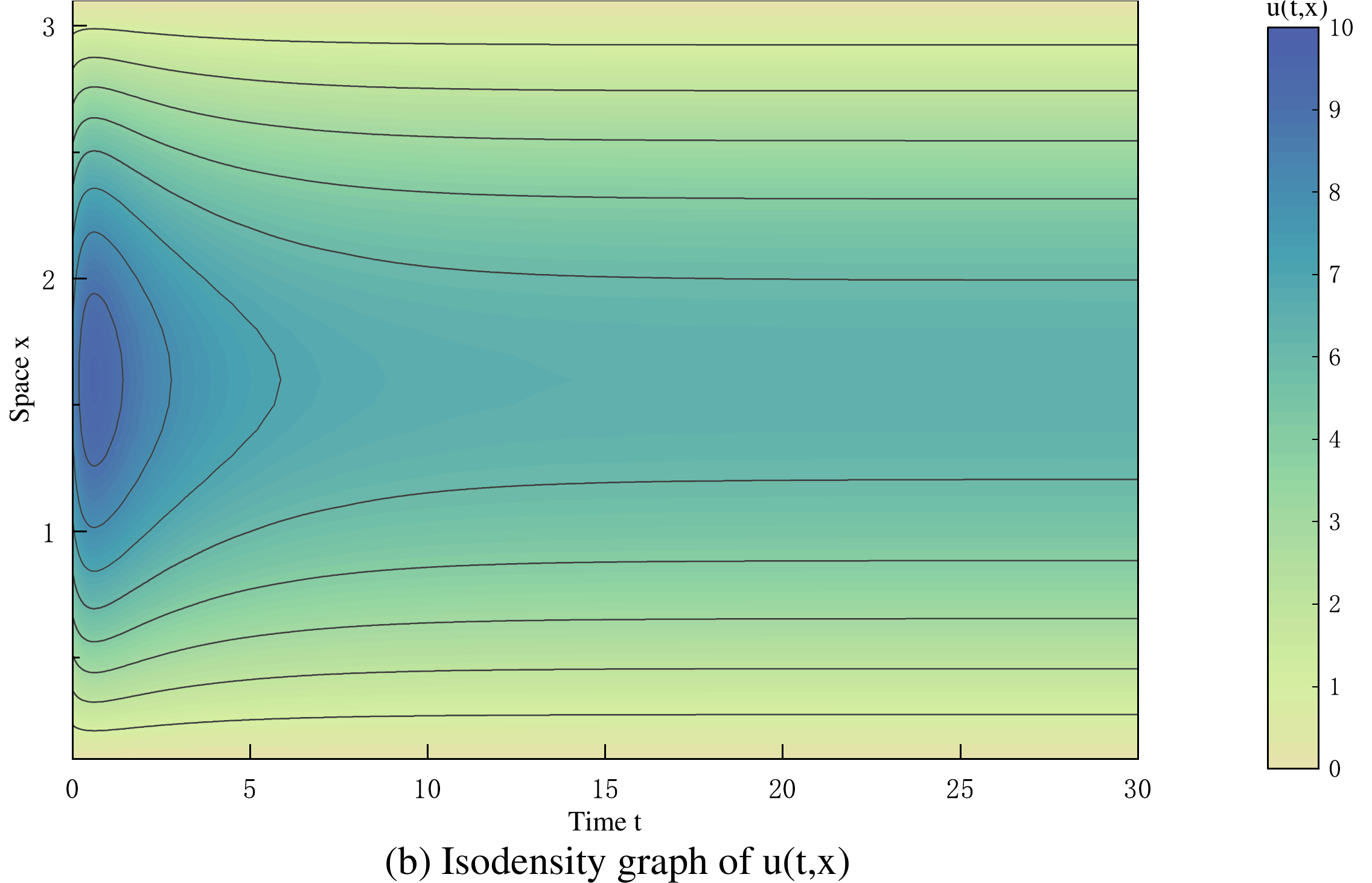}
} }
\subfigure{ {
\includegraphics[width=0.30\textwidth]{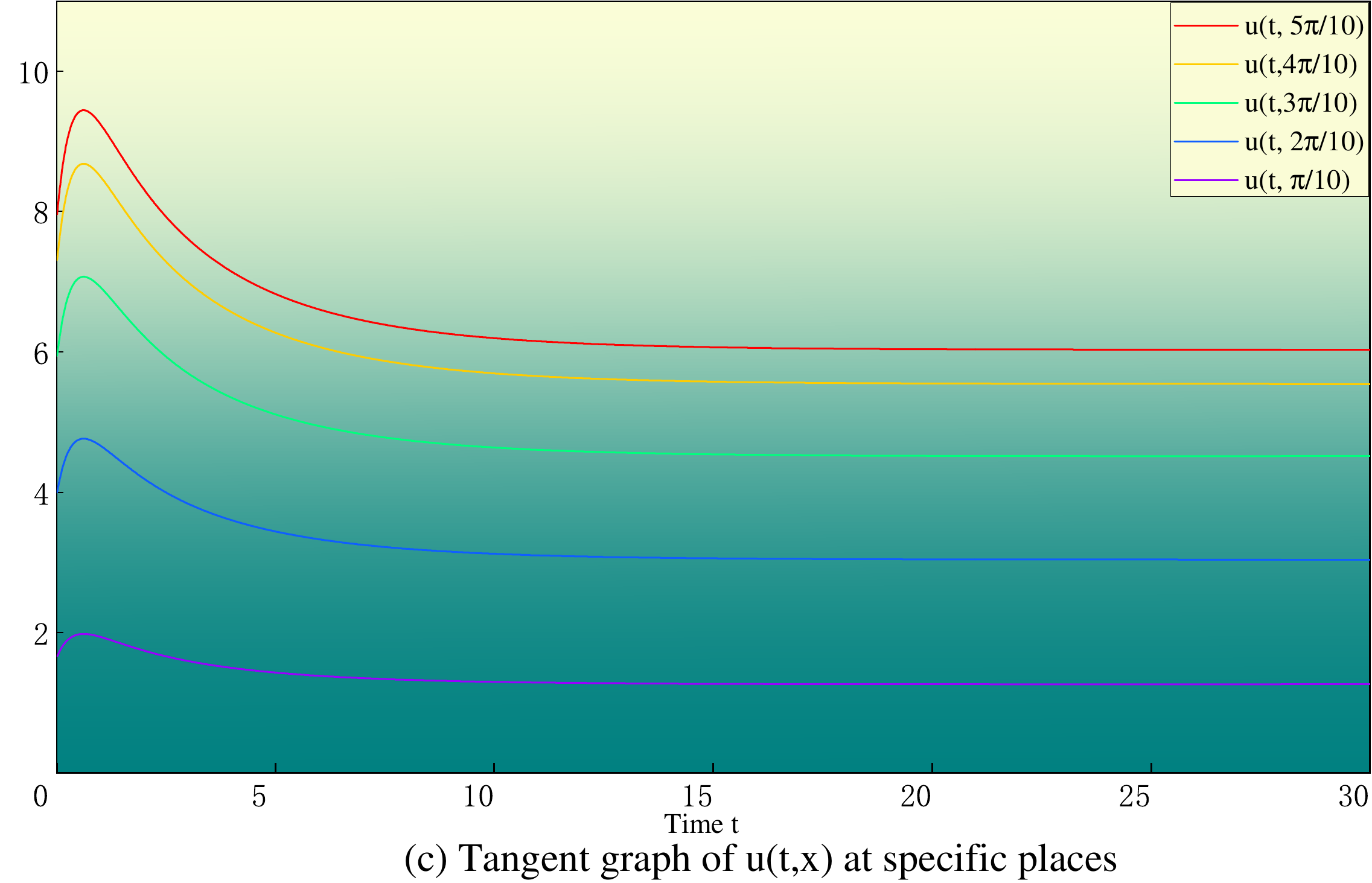}
} }
\subfigure{ {
\includegraphics[width=0.30\textwidth]{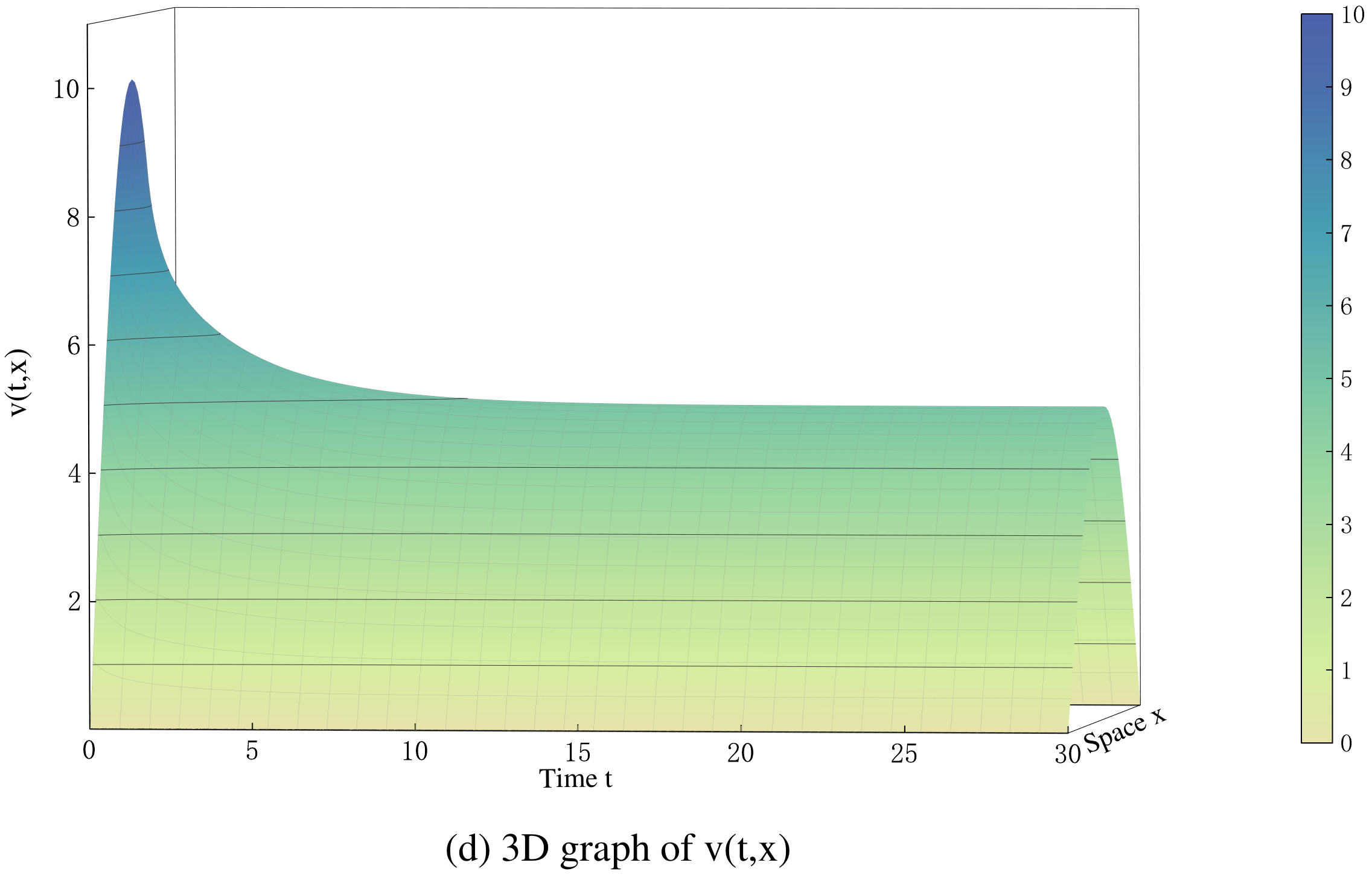}
} }
\subfigure{ {
\includegraphics[width=0.30\textwidth]{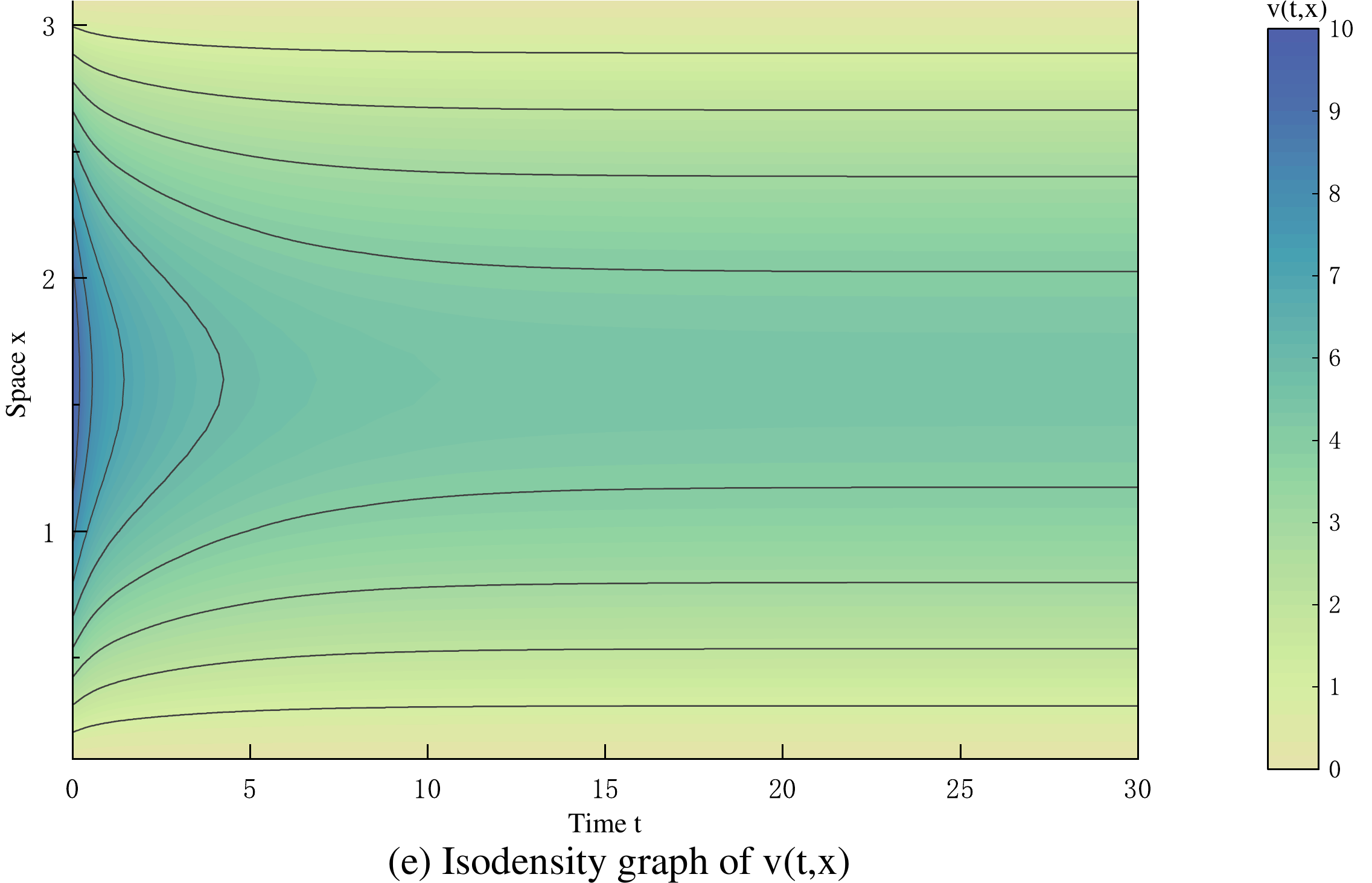}
} }
\subfigure{ {
\includegraphics[width=0.30\textwidth]{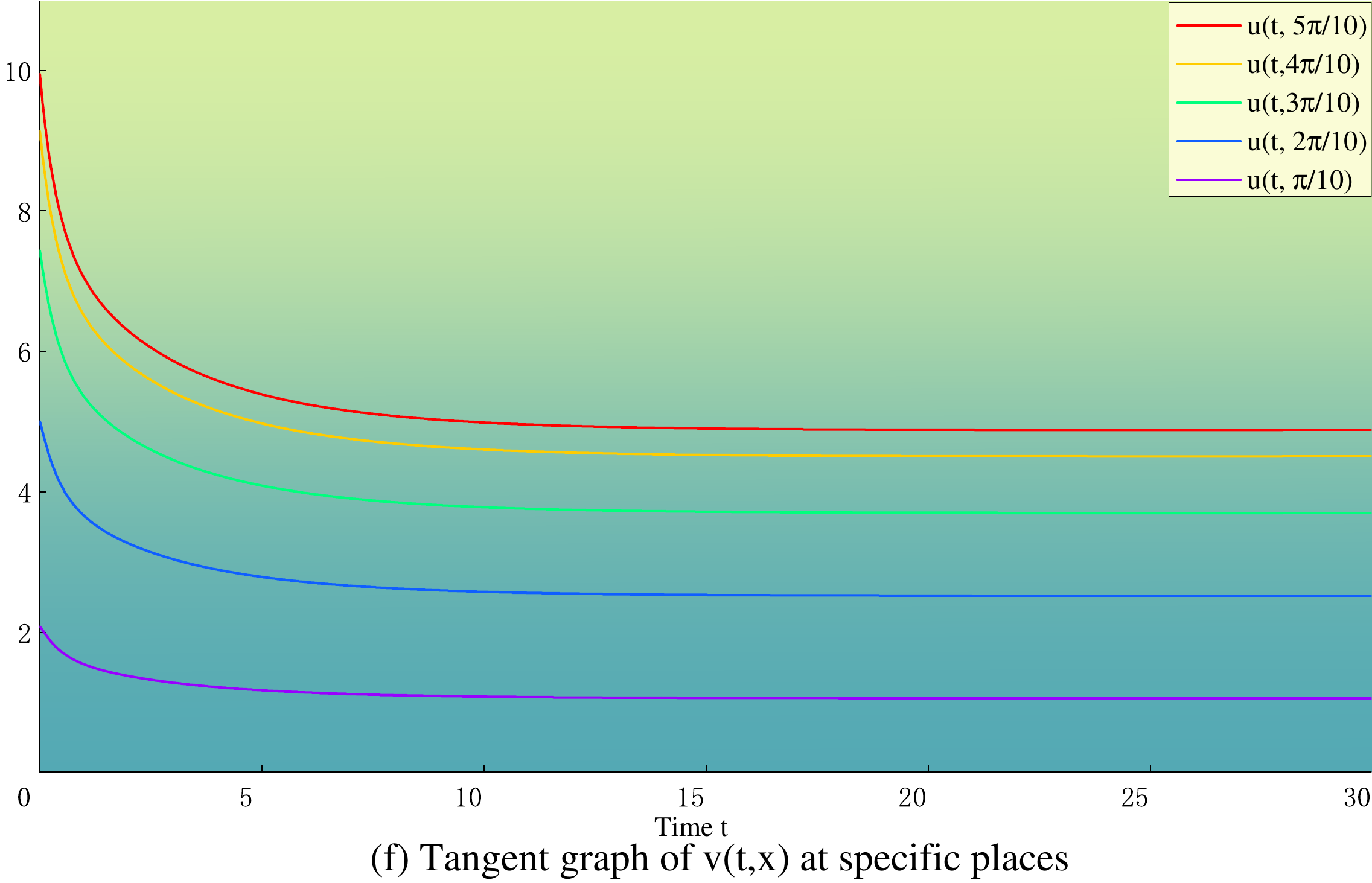}
} }
\caption{When $\rho(t)=1$(fixed region), graphs (a)-(f) show that $u$ and $v$ converge to steady states.}
\label{G}
\end{figure}
\begin{figure}[!ht]
\centering
\subfigure{ {
\includegraphics[width=0.30\textwidth]{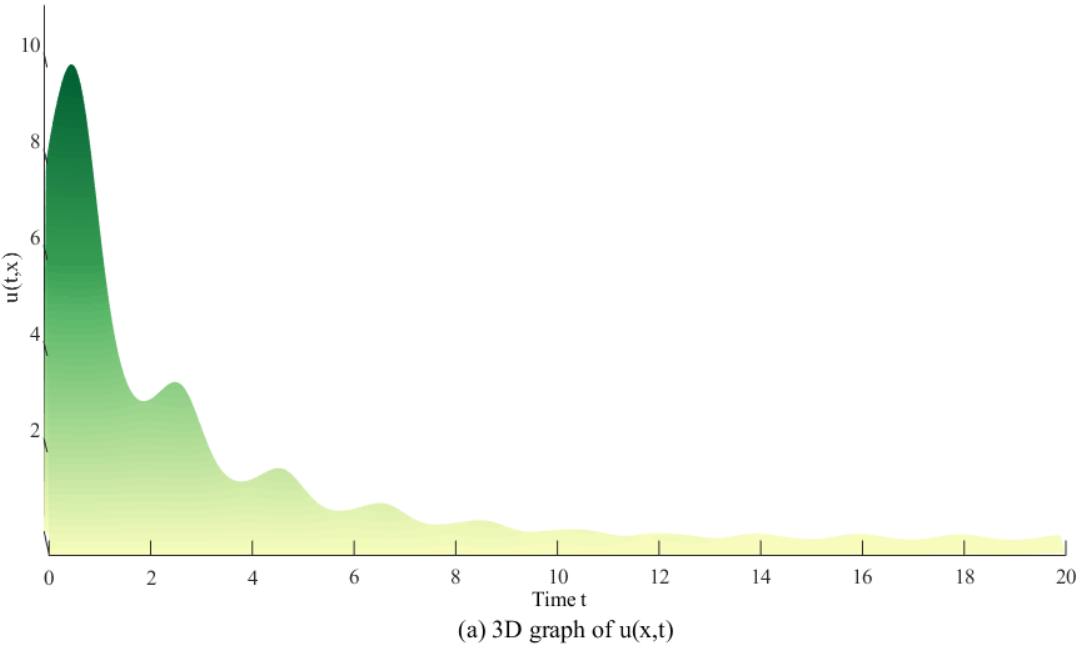}
} }
\subfigure{ {
\includegraphics[width=0.30\textwidth]{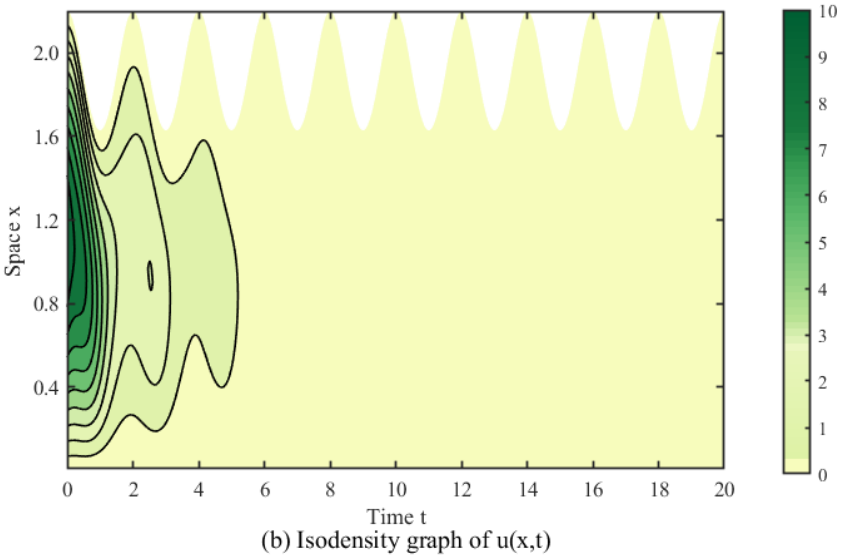}
} }
\subfigure{ {
\includegraphics[width=0.30\textwidth]{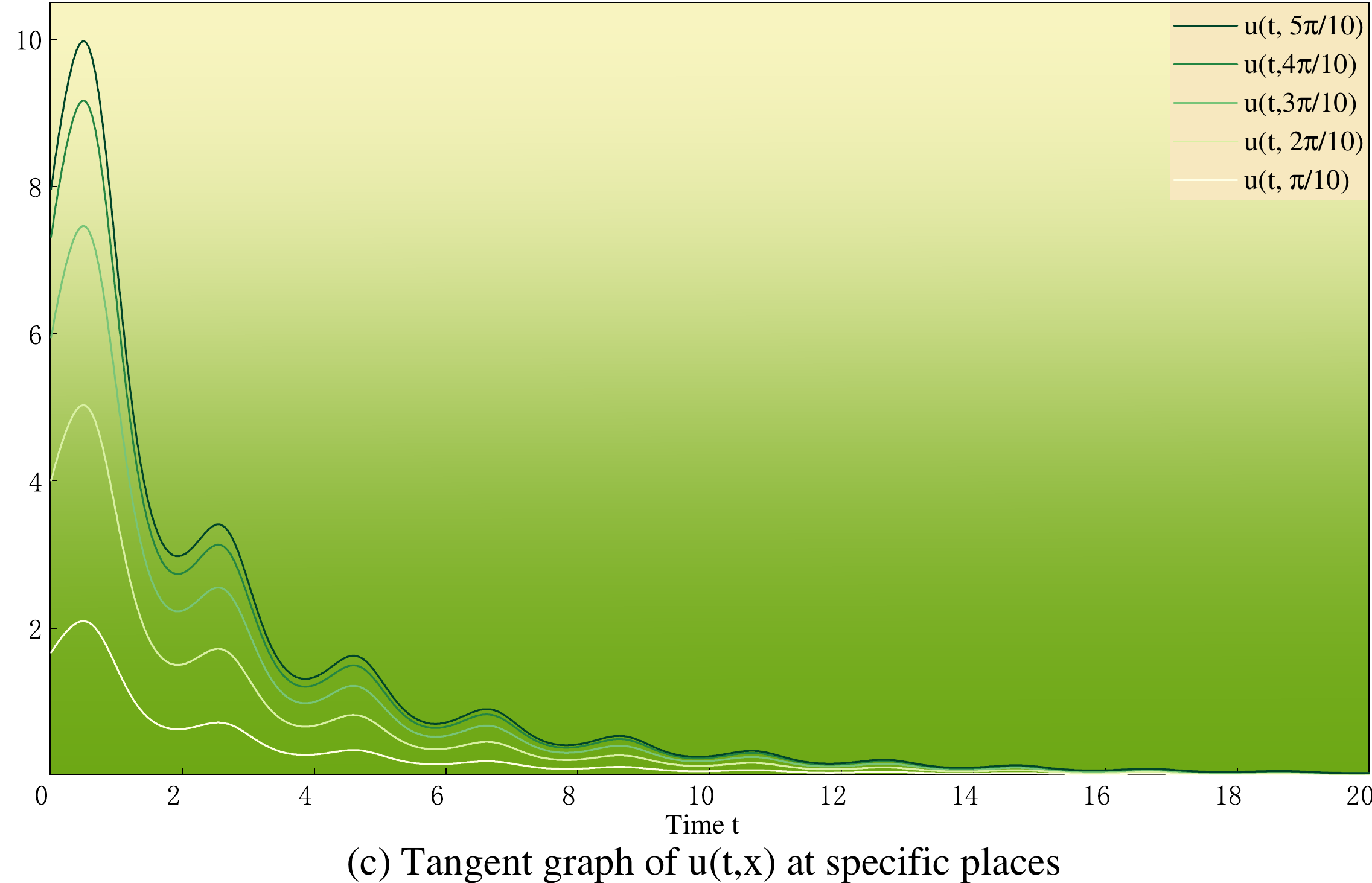}
} }
\subfigure{ {
\includegraphics[width=0.30\textwidth]{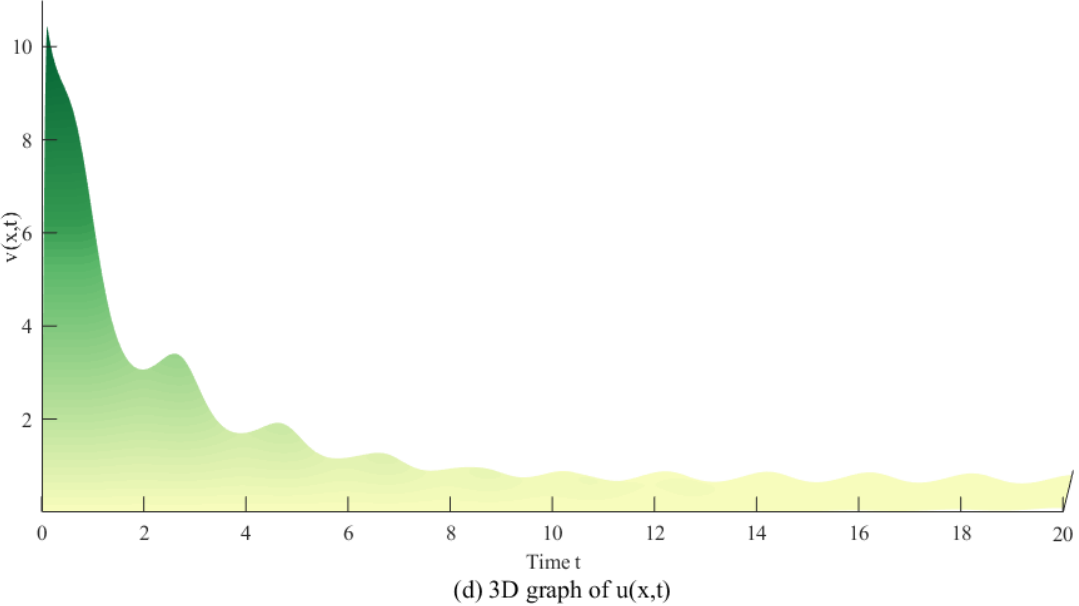}
} }
\subfigure{ {
\includegraphics[width=0.30\textwidth]{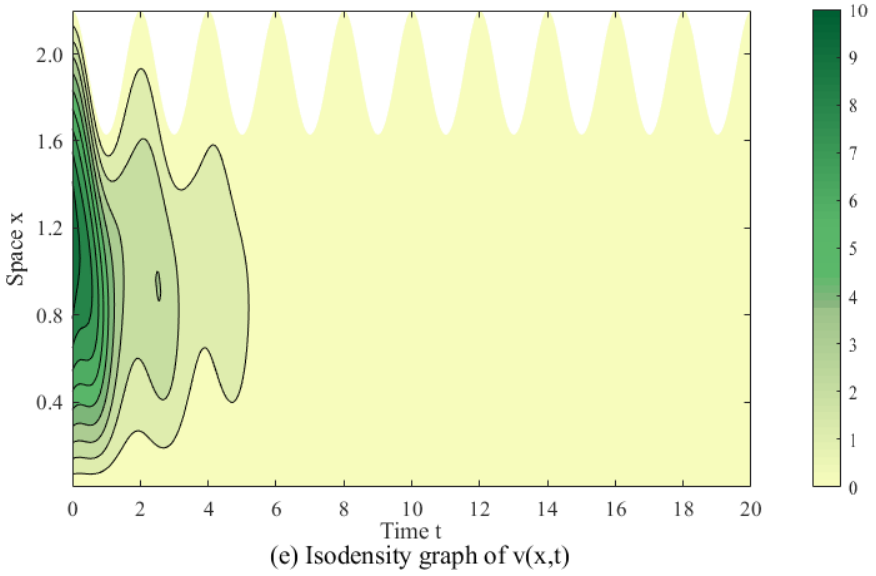}
} }
\subfigure{ {
\includegraphics[width=0.30\textwidth]{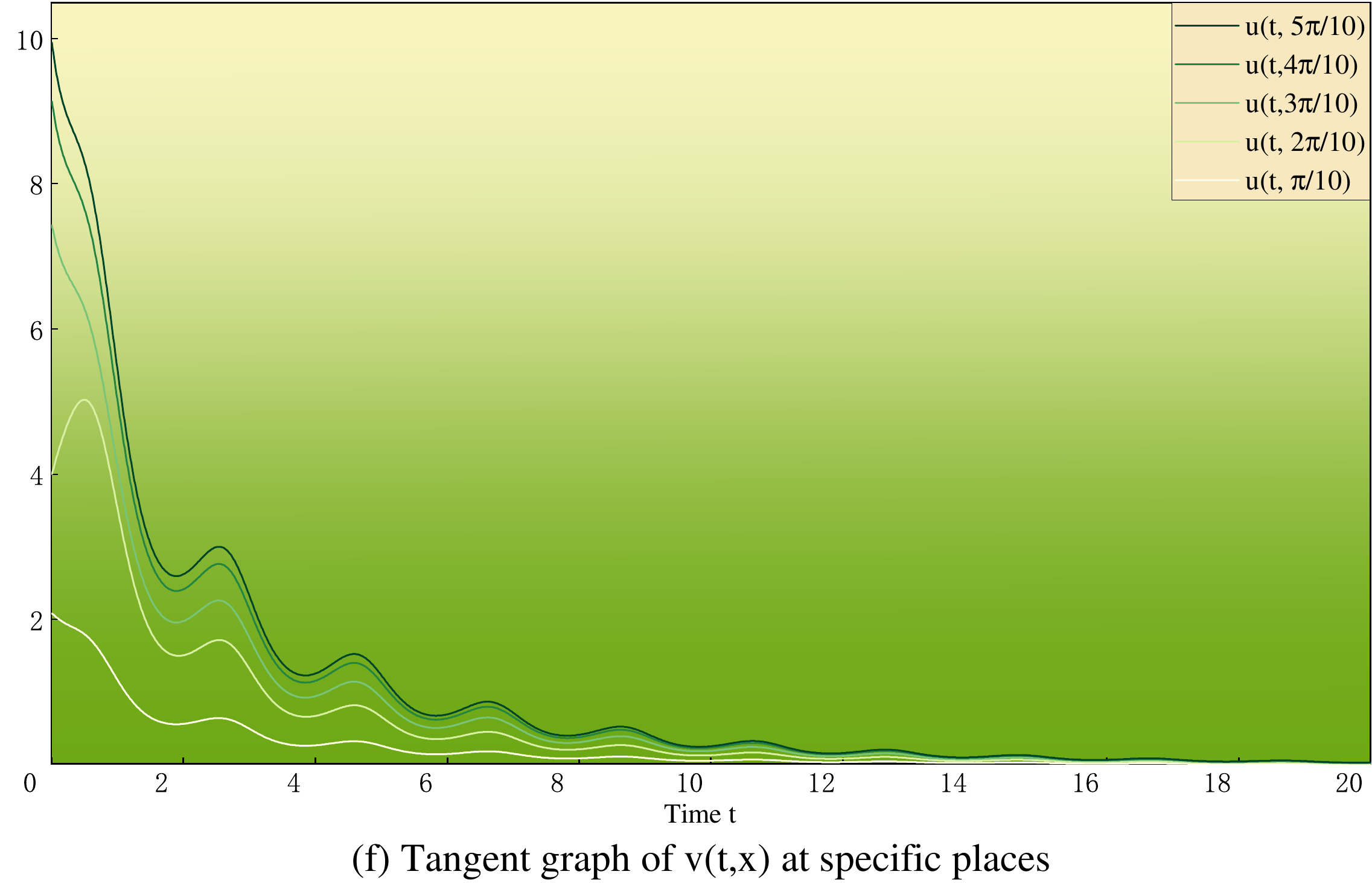}
} }
\caption{\small{When $\rho(t)=0.7e^{-0.15(1-\cos(\pi t))}$(evolving region), graphs (a)-(f) show $(u(t,x), v(t,x))$ decays to $(0, 0)$.}}
\label{H}
\end{figure}
Take $\rho(t)$ to be $1$. Then, it follows from \autoref{theorem 4-3}\textcolor[rgb]{0.00,0.00,1.00}{$(3)$} that $\lambda_{1}\big(\rho(t)\big)=-0.469$. Therefore, \textcolor[rgb]{0.00,0.00,1.00}{Theorems} \ref{theorem 3-2} and \ref{theorem 3-3} yield that model \eqref{1-2-20} has a unique positive steady state solution, and it is globally asymptotically stable. In fact, it can be seen from \autoref{G}\textcolor[rgb]{0.00,0.00,1.00}{(a, c, d, f)}  that the bacteria and infected individuals respectively converge to positive distributions over time. This conclusion is the same as given by \autoref{theorem 3-3}.

Take $\rho(t)$ to be $0.7e^{-0.15(1-\cos(\pi t))}$. With the help of \autoref{theorem 4-3}\textcolor[rgb]{0.00,0.00,1.00}{$(2)$}, one can then obtain that $\lambda_{1}\big(\rho(t)\big)>0.113$. Therefore, \autoref{theorem 3-1} yields that $(0,0)$ is globally asymptotically stable. Actually, observing \autoref{H} gives us that the bacteria and infected individuals at each location tend together to zero over time. This coincides with the conclusion of \autoref{theorem 3-1}.

When the habitat of populations does not change with time, the diseases persist. However, when the habitat of populations changes with a small evolutionary rate $\rho(t)$($\xoverline{\rho^{-2}}>1$), the disease dies out. Therefore, small evolving rate play a positive role in the prevention and control of the diseases.

\textcolor[rgb]{0.00,0.00,1.00}{Examples} \ref{exm3} and \ref{exm4} demonstrate that regional evolving rate plays an important role in the spread of the diseases. It can influence-or even change-the dynamical behaviours of the diseases. Moreover, the principal eigenvalue may decrease with increasing evolutionary rate. Most importantly, small evolutionary rate favours the diseases extinction. This has guiding significance in disease prevention and control.
\section{\bf Conclusion and discussion}\label{Section-7}
This paper has developed an impulsive faecal-oral model in a periodically evolving environment in order to understand how pulsed intervention and regional evolution together influence the spread of faecal-oral diseases. The existence, uniqueness, and positivity of the solution of the model are first given in \autoref{theorem 1-1}. To overcome the difficulties caused by pulse, we transform the related impulse period eigenvalue problem. Then, the existence of principal eigenvalue $\lambda_{1}$ is presented in \autoref{theorem 2-1}. Based on the principal eigenvalue, \textcolor[rgb]{0.00,0.00,1.00}{Theorems} \ref{theorem 3-1}-\ref{theorem 3-3} have provided the threshold dynamical behaviours of the model. More specifically, when $\lambda_{1}>0$, disease-free stable state is global asymptotic stability; when $\lambda_{1}<0$, the model has a unique positive periodic steady state, and this stable state is also global asymptotic stability.

With the help of the adjoint problem, this paper also has provided the monotonicity of the principal eigenvalue with respect to initial region and impulse intensity, see \autoref{theorem 4-1}. In addition, \autoref{theorem 4-2} has given the estimation of the principal eigenvalue in some special cases. For the monotonicity of the principal eigenvalue with respect to regional evolution rate, there has been some works on a equation, see, for instance, \cite{19, 2-5}. However, when the model is being represented by two-or even more-equations, this monotonicity has so far not been completely resolved. In this paper, we merely have used numerical simulations to obtain the conclusion that the principal eigenvalue is decreasing about regional evolutionary rate. This will be our future research direction.

\section*{\bf Declaration of competing interest}
 The authors have no conflict of interest.

\end{document}